\input amstex
\documentstyle {amsppt}

\pagewidth{32pc} 
\pageheight{45pc} 
\mag=1200
\baselineskip=15 pt

\hfuzz=5pt
\topmatter
\NoRunningHeads 

\title Discrete Localities II 
\endtitle
\author Andrew Chermak 
\endauthor
\date
March 2022 
\enddate

\endtopmatter

\define\w{\widetilde}

\redefine\norm{\trianglelefteq}

\redefine\bar{\overline}

\redefine\maps{\mapsto}
\redefine\i{^{-1}}

\redefine\l{\lambda}
\redefine\s{\sigma}
\redefine\a{\alpha}
\redefine\b{\beta}
\redefine\d{\delta}
\redefine\g{\gamma}
\redefine\e{\epsilon} 

\redefine\t{\tau}
\redefine\r{\rho}

\redefine\G{\Gamma}

\redefine\S{\Sigma}
\redefine\L{\Lambda}

\redefine\<{\langle}
\redefine\>{\rangle}

\redefine\ca{\Cal}

\redefine\D{\Delta}

\redefine\sub{\subseteq} 

\redefine\nsub{\nsubseteq}

\redefine\bX{\bold X}

\redefine\nset{\emptyset} 

\redefine\1{\bold 1} 

\redefine\up{\uparrow} 
\redefine\down{\downarrow}

\redefine\bull{\bullet}

\redefine\bw{\bold w} 
\redefine\bu{\bold u} 
\redefine\bv{\bold v} 
\redefine\bX{\bold X} 
\redefine\bY{\bold Y}

\redefine\ccirc{\circledcirc} 

\document 

\vskip .2in 
\noindent 
{\bf Introduction} 
\vskip .1in 

This paper continues the project begun in [CG], to be referred to henceforth as Part I. 
References to Part I will be made by prefixing an ``I" to the referenced item. Thus I.1.1 provides the 
definition of partial group, I.3.1 provides the definition of locality $(\ca L,\D,S)$, and I.3.4 
introduces the fusion system $\ca F=\ca F_S(\ca L)$. We shall say  
that $(\ca L,\D,S)$ is a locality {\it on} $\ca F$. 

This Part II develops a pair of themes, to be tied together by the notion of a ``proper locality" $\ca L$. 
\roster 

\item "{(1)}" We shall determine the ways in which the set $\D$ of objects can be either expanded or 
contracted, while leaving invariant both the fusion system $\ca F$ and the structure of the lattice of 
partial normal subgroups of $\ca L$.  

\item "{(2)}" We will show that if $\ca L$ is ``proper" (a notion which will be formally introduced 
in 6.7), then $\ca F$ has the various useful properties that one associates with saturated 
fusion systems over finite $p$-groups. 

\endroster 

Before discussing (1) and (2) it may be helpful to briefly review some notions from Part I. Thus, 
let $(\ca L,\D,S)$ be a locality, and let $\bold W(\ca L)$ be the free monoid in the alphabet 
$\ca L$. For each $w\in\bold W(\ca L)$ one has the subset $S_w\sub S$, consisting of all the 
elements $x\in S$ such that $x$ is conjugated sequentially into $S$ by the sequence of entries of $w$,  
and where $S_w=S$ if $w$ is the empty word. We write 
$$ 
\text{$x\maps x^w$ \ \ and\ \ $c_w:S_w to S$}
$$  
for the ``conjugation" mapping of $S_w$ into $S$ given in this way. Here $S_w$ is in fact a subgroup of 
$S$ and $c_w$ is a homomorphism, by I.2.9. For any subgroup $X\leq S_w$ we write $X^w$ for the image of $X$ 
under $c_w$, and if $Y$ is a subgroup of $S$ containing $X^w$ then we shall often abuse the notation by 
allowing ``$c_w$" to denote also the homomorphism $X\to Y$ given by $x\maps x^w$ for $x\in X$. With this 
convention, the fusion system $\ca F_S(\ca L)$ is then by definition the category whose objects are the 
subgroups of $S$ and whose morphisms are the homomorphisms $c_w$ for $w\in\bold W(\ca L)$. 

Now let $\D_0$ be a subset of $\D$, and assume that $\D_0$ is $\ca F$-closed. That is, assume that $\D_0$ is 
non-empty, and that 
$Y\in\D_0$ whenever there exists an $X\in\D_0$ and an $\ca F$-homomorphism from $X$ into $Y$. 
There is then a locality $(\ca L_0,\D_0,S)$ such that $\ca L_0$ is a subset of $\ca L$, and such that 
the inclusion map $\ca L_0\to\ca L$ is a homomorphism of partial groups. These conditions can be  
fulfilled in only one way: by taking the domain $\bold D(\ca L_0)$ for 
the product in $\ca L_0$ to be the set of words $w\in\bold W(\ca L)$ with $S_w\in\D_0$. We then 
say that $\ca L_0$ is the {\it restriction} of $\ca L$ to $\D_0$. 

In the other direction, if $\D^+$ is an $\ca F$-closed set of subgroups of $S$ containing $\ca D$, then 
an {\it expansion} of $\ca L$ to $\D^+$ is a locality $(\ca L^+,\D^+,S)$ such that $\ca L$ is a restriction 
of $\ca L^+$. Sections 4 and 5 of this paper show how to construct expansions of $\ca L$, under  
suitable hypotheses on $\D^+$ and $\ca L$. Moreover, under those hypotheses we shall find that 
$\ca F_S(\ca L^+)=\ca F$ and that $\ca L^+$ is proper if $\ca L$ is proper. We state the main results 
relating to expansions, as follows, postponing for the moment the definition of proper locality and of the 
set $\ca F^s$ of $\ca F$-subcentric subgroups of $S$. 

\proclaim {Theorem A1} Let $(\ca L,\D,S)$ be a proper locality on $\ca F$ and let $\D^+$ be an 
$\ca F$-closed collection of subgroups of $S$ such that $\D\sub\D^+\sub\ca F^s$. Then the following hold. 
\roster 

\item "{(a)}" There exists a proper locality $(\ca L^+,\D^+,S)$, such that $\ca L^+$ is an expansion of 
$\ca L$. Moreover the inclusion map $\ca L\to\ca L^+$ is a homomorphism of partial groups, 
$\ca L^+$ is generated by $\ca L$ as a partial group, and $\ca F_S(\ca L^+)=\ca F$. 

\item "{(b)}" If $(\w{\ca L},\D^+,S)$ is any other proper locality which is an expansion of $\ca L$  
then there is a unique homomorphism $\b:\ca L^+\to\w{\ca L}$ such that $\b$ restricts to the identity 
map on $\ca L$, and $\b$ is then an isomorphism. 

\endroster 
\endproclaim 

\proclaim {Theorem A2} Let $(\ca L,\D,S)$ be a proper locality on $\ca F$, and let $(\ca L^+,\ca D^+,S)$ 
be an expansion of $\ca L$ to an $\ca F$-closed subset $\D^+$ of $F^s$. Let $\frak N$ be the set 
of partial normal subgroups of $\ca L$, and let $\frak N^+$ be the set of partial normal subgroups of 
$\ca L^+$. For each $\ca N\in\frak N$ there is then a unique $\ca N^+\in\frak N^+$ such that 
$\ca N=\ca N^+\cap\ca L$. In particular, we have $S\cap\ca N=S\cap\ca N^+$, and the mapping 
$\ca N\maps\ca N^+$ is an inclusion-preserving bijection $\frak N\to\frak N^+$.  
\endproclaim 

The union of Theorems A1 and A2 may be thought of as a single Theorem A. Theorem A is motivated in 
large part by the goal (to be attained in Part III) of finding an $\ca F$-closed set $\d(\ca F)$ 
of subgroups of $\ca F$ such that there is a proper locality $(\ca L,\d(\ca F),S)$ on $\ca F$, 
and such that each partial normal subgroups $\ca N\norm\ca L$ is itself a proper locality 
$(\ca N,\d(\ca E),T)$ on $\ca E$, where $T=S\cap\ca N$ and where $\ca E$ is the fusion system 
$\ca F_T(\ca N)$.

\vskip .1in 
As mentioned, our second main theme concerns the fusion system $\ca F$ of a proper locality $(\ca L,\D,S)$. 
In this introduction we will draw analogies with some of the basic definitions concerning fusion 
systems over finite $p$-groups (and with which we assume that the reader has some familiarity).  
Readers for whom fusion systems are entirely new may wish to visit section 1 below, and to see the 
definition 8.2 of saturation, before proceeding any further. In order to introduce 
this theme, we first review the notion of finite-dimensionality from Part I. 

Set 
$$
\Omega=\Omega_S(\ca L)=\{S_w\mid w\in\bold W(\ca L)\},  
$$ 
and regard $\Omega$ as a poset via inclusion. Then, by the definition of locality, $\Omega$ is 
finite-dimensional: i.e. there is an upper bound to the lengths of strictly monotone chains in $\Omega$. As 
$\Omega$ is closed under finite intersections (I.2.8(a)), finite-dimensionality implies that $\Omega$ 
is closed under arbitrary intersections, and there is then a mapping $X\maps X^\star$ from the set 
of subgroups of $S$ into $\Omega$, defined by 
$$ 
X^\star=\bigcap\{S_w\mid X\leq S_w\}. 
$$ 
The basic properties of $(^{\star})$ are given by 
I.2.11 and I.2.12, and these include the fact that $(^{\star})$ is equivariant with respect to fusion 
in $\ca F_S(\ca L)$. That is, for any $w\in\bold W(\ca L)$ and any subgroup $X\leq S_w$ we have 
$$ 
(X^w)^\star=(X^\star)^w.  
$$ 
We say that $(\star,\Omega)$ is the {\it stratification} on $\ca F_S(\ca L)$ induced from $\ca L$.  
Such structures will be studied in detail in sections 2 and 3. 

Some basic notions from the theory of fusion systems over a finite $p$-group can be formulated in the 
present context by means of the stratification $(\star,\Omega)$. First, for $X\leq S$ define the 
{\it $\Omega$-dimension} 
$dim_{\Omega}(X)$ of $X$ to be the supremum of the lengths of monotonically ascending chains in $\Omega$ 
that terminate in $X^*$. 

\definition {Definition 0.1} Let $(\ca L,\D,S)$ be a locality, set $\ca F=\ca F_S(\ca L)$, let 
$(\star,\Omega)$ be the stratification on $\ca F$ induced from $\ca L$, and let $X\leq S$ be a 
subgroup of $S$. Then $X$ is {\it fully $\star$-normalized} in $\ca F$ if 
$$ 
dim_\Omega(X)\geq dim_\Omega(Y) 
$$ 
for each $\ca F$-conjugate $Y$ of $X$, and $X$ is {\it fully $\star$-centralized} in $\ca F$ if 
$$ 
dim_\Omega(C_S(X)X)\geq dim_\Omega(C_S(Y)Y) 
$$ 
for each $\ca F$-conjugate $Y$ of $X$. Further, given a subgroup $Y\leq S$, we say that $Y$ is: 
\roster 

\item "{$\cdot$}" {\it centric} in $\ca F$ if $C_S(X)\leq X$ for every  $\ca F$-conjugate $X$ of $Y$.   

\item "{$\cdot$}" {\it $\star$-radical} in $\ca F$ if there exists a fully $\star$-normalized 
$\ca F$-conjugate $X$ of $Y$ such that $O_p(N_{\ca F}(X))=X$. 

\item "{$\cdot$}" {\it $\star$-subcentric} in $\ca F$ if there exists a fully $\star$-normalized 
$\ca F$-conjugate $X$ of $Y$ such that $O_p(N_{\ca F}(X))$ is centric in $\ca F$. 

\endroster 
Write $\ca F^c$ for the set of all $\ca F$-centric subgroups of $S$, $\ca F^{cr}$ for the set of 
all $Y\in\ca F^c$ such that $Y$ is $\star$-radical in $\ca F$, and $\ca F^s$ for the set of all 
$\star$-subcentric subgroups of $S$. 
\enddefinition 

Here the definitions of $N_{\ca F}(X)$ (and of $C_{\ca F}(X)$) are the same as for fusion systems 
over a finite $p$-group (for which see 1.5 below). Also, as in the finite case, if $\ca E$ is a fusion 
system on a $p$-group $T$ then $O_p(\ca E)$ is by definition the largest subgroup $X$ of $T$ such that 
$\ca E=N_{\ca E}(X)$. 

Two remarks are in order at this point. First, it would seem that the sets $\ca F^{cr}$ and $\ca F^s$ 
depend not only on $\ca F$ but on the stratification $(\Omega,\star)$, and hence on $\ca L$. But in 
fact, if $\ca L'$ is any other proper locality on $\ca F$, with the stratification $(\Omega',\star')$ 
induced on $\ca F$, then it will turn out (cf. 2.9) that a subgroup $X$ of $S$ is fully $\star$-normalized if 
and only if $X$ is fully $\star'$-normalized. Thus $(\Omega',\star')$ yields the 
same result for $\ca F^{cr}$ and $\ca F^s$ as does $(\Omega,\star)$, and so the notation (in which we 
have omitted any reference to $\ca L$) is justified. Second, it should be mentioned that 
even in the case of fusion systems over finite $p$-groups, the above definition of 
$\ca F$-radical subgroup is different from the standard 
one - by which a subgroup $Y$ of $S$ is said to be $\ca F$-radical if $Inn(Y)$ is the largest normal 
$p$-subgroup of $Aut_{\ca F}(Y)$. The two definitions will be seen to be equivalent 
(cf. 8.13) when $\ca F$ is the fusion system of a proper locality.

\definition {Definition 0.2} The locality $(\ca L,\D,S)$ on $\ca F$ is {\it proper} if: 
\roster 

\item "{(PL1)}" $\ca F^{cr}\sub\D$.  

\item "{(PL2)}" For each $P\in\D$ the group $\ca L_P:=N_{\ca L}(P)$ is of characteristic $p$. 
That is, $C_{\ca L_P}(O_p(\ca L_P))\leq O_p(\ca L_P)$. 

\item "{(PL3)}" $S$ has the normalizer-increasing property. That is, if $P$ and $Q$ are subgroups 
of $S$ with $P<Q$ ($P$ a proper subgroup of $Q$) then $P<N_Q(P)$. 

\endroster 
\enddefinition 

We note that by the Appendix A to Part I, the class of groups $S$ having the normalizer-increasing property 
includes all homomorphic images of $p$-groups having a faithful finite-dimensional representation over 
the algebraic closure of a finite field. 
\vskip .1in

\proclaim {Theorem B} Let $(\ca L,\D,S)$ be a proper locality on $\ca F$. let $A$ be fully normalized 
in $\ca F$, and let $B$ be fully centralized in $\ca F$, with respect to the stratification $(\Omega,\star)$ 
induced from $\ca L$. Then the following hold. 
\roster 

\item "{(a)}" $\ca F$ is inductive. That is, for each $\ca F$-conjugate $A'$ of $A$ there exists 
an $\ca F$-homomorphism $\phi:N_S(A')\to N_S(A)$ with $A'\phi=A$. 

\item "{(b)}" There exists a proper locality $(\ca L_A,\D_A,N_S(A))$ on $N_{\ca F}(A)$, and a proper 
locality $(\ca C_B,\S_B,C_S(B))$ on $C_{\ca F}(B)$. 

\item "{(c)}" $\ca F$ is saturated. 

\endroster 
\endproclaim 

The proofs of Theorem B and of the related results 3.8 and 8.1, are based largely on ideas taken from the 
proof of [Theorem 2.2 in BCGLO]. Further information on $\ca F$ of the kind carried by Theorem B may be found 
in results 8.5 and 8.6. In 8.10 we show via arguments taken from [He2] that if $\ca F$ is the 
fusion system of a proper locality $(\ca L,\D,S)$ then the set $\ca F^s$ of $\ca F$-subcentric subgroups 
of $S$ is $\ca F$-closed. Combining this result with Theorem A1 yields the notion of the {\it subcentric 
closure} $(\ca L^s,\ca F^s,S)$ of $\ca L$. This will be an important tool in the remaining parts of this 
series. 

Finally, the author would like to express his gratitude to Alex Gonzales, for suggesting the 
possibility of this project, for fruitful discussions, and for his assistance with Appendix A 
where an equivalence is established between the $p$-local compact groups of [BLO2] and certain 
proper localities. 

\vskip .2in 
\noindent 
{\bf Section 1: Abstract fusion systems} 
\vskip .1in

This initial section consists of definitions, some of their basic consequences, and a few 
illustrative examples.

\definition {Definition 1.1} Let $S$ be a group. A {\it fusion system} on $S$ is a category $\ca F$ 
whose objects are the subgroups of $S$, and whose morphisms satisfy the following 
conditions.  
\roster 

\item "{(1)}" Each $\ca F$-morphism $X\to Y$ is an injective homomorphism of groups. 

\item "{(2)}" If $g\in S$, and $X$ and $Y$ are subgroups of $S$ such that 
$X^g\leq Y$ then the conjugation map $c_g:X\to Y$ is an $\ca F$-morphism. 

\item "{(3)}" If $\phi:X\to Y$ is an $\ca F$-morphism then the bijection $X\to Im(\phi)$ defined by  
$\phi$ is an $\ca F$-isomorphism. 

\endroster 
\enddefinition  

One most often refers to $\ca F$-morphisms as $\ca F$-homomorphisms, in order to emphasize that they are 
bona fide homomorphisms of groups. 
Notice that (2) implies that all inclusion maps between subgroups of $S$ are $\ca F$-homomorphisms. 
It follows from this observation, and from (3), that if $\phi:X\to Y$ is an $\ca F$-homomorphism which maps 
a subgroup $X_0$ of $X$ into a subgroup $Y_0$ of $Y$, then the map $\phi_0:X_0\to Y_0$ induced by $\phi$ 
is again an $\ca F$-homomorphism. Such a homomorphism $\phi_0$ will be referred to as a {\it restriction} 
of $\phi$, and we say also that $\phi$ is an {\it extension} of $\phi_0$. 
\vskip .1in  

\definition {Examples} 
\roster 

\item "{(1)}" Let $G$ be a group and let $S$ be a subgroup of $G$. There is then a fusion system 
$\ca F:=\ca F_S(G)$ on $S$ in which the set $Hom_{\ca F}(X,Y)$ 
of $\ca F$-homomorphisms $X\to Y$ is the set of maps 
$c_g:X\to Y$ given by conjugation by those elements $g\in G$ for which $X^g\leq Y$. 

\item "{(2)}" The {\it total fusion system} $\bar{\ca F}_S$ on $S$, is characterized by 
$Hom_{\bar{\ca F}_S}(X,Y)$ being the set of all injective homomomorphisms $X\to Y$. 

\item "{(3)}" Let $(\ca L,\D,S)$ be a locality. As discussed in I.3.4, there is then a fusion system 
$\ca F=\ca F_S(\ca L)$, in which $Hom_{\ca F}(X,Y)$ is the set of all conjugation maps 
$c_w:X\to Y$ where $w\in\bold W(\ca L)$, $X\leq S_w$, and $X^w\leq Y$. In the special case that $\ca L$ is 
a group we have $c_w=c_{\Pi(w)}$ by I.2.3(c), and so this definition of $\ca F$ agrees with the definition 
given in (1) in the case that $\ca L$ is a group. 

\endroster 
\enddefinition 

\definition {Definition 1.2} Let $\ca F$ be a fusion system on $S$, and let $\ca F'$ be a 
fusion system on $S'$. A homomorphism $\a:S\to S'$ is {\it fusion-preseving} (relative to 
$\ca F$ and $\ca F'$) if, for each $\ca F$-homomorphism $\phi:P\to Q$, there exists an $\ca F'$-homomorphism 
$\psi:P\a\to Q\a$ such that $\a\mid_P\circ\psi=\phi\circ\a\mid_Q$. 
\enddefinition 

Notice that since the $\ca F'$-homomorphism $\psi$ in the preceding definition is injective, it is uniquely 
determined. Thus, if $\a:S\to S'$ is a fusion-preserving homomorphism then $\a$ induces a mapping 
$\a_{P,Q}:Hom_{\ca F}(P,Q)\to Hom_{\ca F'}(P\a,Q\a)$, for each pair $(P,Q)$ of 
subgroups of $S$. The following result follows directly from this observation. 

\proclaim {Lemma 1.3} Let $\ca F$ be a fusion system on $S$, let $\ca F'$ be a fusion system 
on $S'$, and let $\a:S\to S'$ be a fusion-preserving homomorphism. Then the mapping $P\maps P\a$ 
from objects of $\ca F$ to objects of $\ca F'$, together with the set of mappings 
$$ 
\a_{P,Q}:Hom_{\ca F}(P,Q)\to Hom_{\ca F'}(P\a,Q\a) \ \ (P,Q\leq S) 
$$ 
defines a functor $\a^*:\ca F\to\ca F'$. 
\qed 
\endproclaim

In view of the preceding result, a fusion-preserving homomorphism $\a$ will also be referred to as a 
homomorphism of fusion systems. Notice that the inverse of a fusion-preserving isomorphism is 
fusion-preserving, and is therefore an isomorphism of fusion systems.

\definition {1.4 (Subsystems and generation)} In the special case where $S\leq S'$ and the inclusion map 
$S\to S'$ is fusion-preserving, we say that $\ca F$ is a {\it subsystem} of $\ca F'$. Here are a few 
examples, by way of introducing some further terminology and notation.  
\roster 

\item "{(1)}" Every fusion system $\ca F$ on $S$ is a subsystem of the total fusion system 
$\bar{\ca F}_S$ on $S$. 

\item "{(2)}" $\ca F_S(S)$ is a subsystem of $\ca F$ for each fusion system $\ca F$ on $S$, by 1.1(2). 
We shall refer to $\ca F_S(S)$ as the {\it trivial fusion system} on $S$.  

\item "{(3)}" Let $\ca F$ be a fusion system on $S$, and let $\frak E$ be a non-empty set of subsystems  
of $\ca F$, such that each $\ca E\in\frak E$ is again a fusion system on $S$. For $\ca E\in\bold E$, write 
$Hom(\ca E)$ for the set of all $\ca E$-homomorphisms. Then $\bigcap\{Hom(\ca E)\mid \ca E\in\bold E\}$ 
is the set of homomorphisms of a subsystem of $\ca F$, to be denoted $\bigcap\bold E$. 

\item "{(4)}" Let $\ca F$ be a fusion system on $S$, let $T\leq S$ be a subgroup of $S$and let $\Phi$ 
be a set of $\ca F$-isomorphisms $\phi:X\to Y$ such that the domain $X$ of $\phi$ and the 
image $Y$ of $\phi$ are subgroups of $T$ for all $\phi\in\Phi$. There is then a subsystem $\<\Phi\>_T$ of 
$\bar{\ca F}_T$, defined as the intersection (as defined in (3)) of the subsystems $\ca E$ 
of $\bar{\ca F}_T$ on $T$, such that $\Phi$ is a subset of $Hom(\ca E)$. Then $\<\Phi\>_T$ is a 
subsystem also of $\ca F$, to be referred to as the subsystem {\it generated by} $\Phi$. Thus, 
$\<\Phi\>_T$ is the fusion system on $T$ whose isomorphisms are identity maps and 
compositions of restrictions of members of $\Phi$. 

\endroster 
\enddefinition

For the remainder of this section, $\ca F$ will be a fixed fusion system on $S$. For any subgroup $X\leq S$ 
write $X^{\ca F}$ for the set of subgroups of $S$ that are $\ca F$-isomorphic to $X$. The elements 
of $X^{\ca F}$ are the {\it $\ca F$-conjugates} of $X$.

\definition {Definition 1.5} Let $V$ be a subgroup of $S$ and let $\Phi$ be the set of 
$\ca F$-isomorphisms $\phi:X\to Y$ such that $X$ and $Y$ are subgroups of $N_S(V)$, and such that 
$\phi$ extends to an $\ca F$-isomorphism $\phi':XV\to YV$ which restricts to an automorphism of $V$. 
Let $\Psi$ be the subset of $\Phi$ consisting of those $\ca F$-isomorphisms $\psi:X\to Y$ such that $X$ and 
$Y$ are subgroups of $C_S(T)$, and such that $\psi$ extends to an $\ca F$-isomorphism $\psi':XV\to YV$ whose 
restriction to $V$ is the identity map. Write $N_{\ca F}(V)$ for $\<\Phi\>_{N_S(V)}$, and 
$C_{\ca F}(V)$ for $\<\Psi\>_{C_S(V)}$. 
\enddefinition

The following result is immediate from the definitions.

\proclaim {Lemma 1.6} Let $U,V\leq S$ be subgroups of $S$, and suppose that there exists an 
$\ca F$-isomorphism $\a:N_S(U)\to N_S(V)$ such that $U\a=V$. Then $\a$ is an isomorphism 
$N_{\ca F}(U)\to N_{\ca F}(V)$. Also, if $\b:C_S(U)U\to C_S(V)V$ is an $\ca F$-isomorphism 
such that $U\b=V$, then the restriction of $\b$ to $C_S(U)$ is an isomorphism 
$C_{\ca F}(U)\to C_{\ca F}(V)$. 
\qed  
\endproclaim

\definition {Definition 1.7} A set $\G$ of subgroups of $S$ is {\it $\ca F$-invariant} if 
$\G$ is closed with respect to $\ca F$-conjugacy ($X\in\G\implies X^{\ca F}\sub\G$). We say that $\G$ is 
{\it $\ca F$-closed} if $\G$ is non-empty, $\ca F$-invariant, and closed with respect to overgroups in 
$S$ ($X\in\G$, $X\leq Y\leq S$ $\implies$ $Y\in\G$). 
\enddefinition 

\definition {Definition 1.8} Let $T$ be a subgroup $S$. Then: 
\roster 

\item "{$\cdot$}" $T$ is {\it normal} in $\ca F$ if $\ca F=N_{\ca F}(T)$. 

\item "{$\cdot$}" $T$ is {\it strongly closed} in $\ca F$ if the set $\G$ of all subgroups of $T$ is 
$\ca F$-invariant. 

\item "{$\cdot$}" $T$ is {\it weakly closed} in $\ca F$ if $\{T\}$ is $\ca F$-invariant. 

\endroster 
\enddefinition 

\proclaim {Lemma 1.9} Let $\ca U$ be a set of subgroups of $S$ such that each $U\in\ca U$ is normal in 
$\ca F$. Then the subgroup $\<\ca U\>$ of $S$ generated by $\ca U$ is normal in $\ca F$. 
Thus, there is a largest subgroup of $S$ which is normal in $\ca F$. 
\endproclaim 

\demo {Proof} Zorn's lemma yields the existence of a maximal $U\in\ca U$. Let also $V\in\ca U$, and let 
$\phi:X\to Y$ be an $\ca F$-isomorphism. Then $\phi$ extends to an $\ca F$-isomorphism 
$\phi_1:XU\to YU$ with $U\phi_1=U$, and then $\phi_1$ extends to an $\ca F$-isomorphism 
$\phi_2:XUV\to YUV$ with $V\phi_2=V$. Thus $UV\norm\ca F$, and the maximality of $U$ then yields 
$V\leq U$. 
\qed 
\enddemo 

We end this section by introducing a category of localities, and a functor  to the 
category $\frak F$ of fusion systems. Recall from I.3.4 that to any locality $(\ca L,\D,S)$ there is 
associated a fusion system $\ca F:=\ca F_S(\ca L)$, generated by the conjugation maps $c_g:S_g\to S_{g\i}$ 
for $g\in\ca L$. Also from I.3.4 we have the conjugation maps $c_w:S_w\to S_{w\i}$, for $w\in\bold W(\ca L)$. 
Thus each $\ca F$-isomorphism is the restriction of some $c_w$ to a subgroup of $S_w$. 

\proclaim {Lemma 1.10} Let $(\ca L,\D,S)$ and $(\ca L',\D',S')$ be localities, with fusion systems 
$\ca F=\ca F_S(\ca L)$ and $\ca F'=\ca F_{S'}(\ca L')$. Let $\l:\ca L\to\ca L'$ be a homomorphism of 
partial groups such that $S\l\leq S'$. Then the restriction $\l\mid_S$ of $\l$ to $S$ is fusion-preserving 
with respect to $\ca F$ and $\ca F'$. That is, $\l\mid_S$ is a morphism $\ca F\to\ca F'$ in the category 
of fusion systems. Moreover, if $\l$ is a projection (cf. I.5.4) then $\l\mid_S$ maps $Ob(\ca F)$ onto 
$Ob(\ca F')$, and maps $Mor(\ca F)$ onto $Mor(\ca F')$. 
\endproclaim 

\demo {Proof} For any subgroup $X\leq S$ and any $g\in\ca L$ write $\bar X$ for $X\l$ and $\bar g$ for 
$g\l$. Also, for $X\leq S$ write $\l_X$ for $\l\mid X$. Thus, $\l_S$ is fusion-preserving if for each 
$\ca F$-homomorphism $\a:X\to Y$ there exists an $\ca F'$-homomorphism $\b:\bar X\to\bar Y$ such that 
$$ 
\l_X\circ\b=\a\circ\l_Y. \tag*
$$ 
(The uniqueness of $\b$ then follows from the observation that $\bar x\b=(x\a)\l_Y$ for all $x\in X$.) 
In the case that $\a=c_g$ for some $g\in\ca L$ the formula (*) is immediate from the definition of 
homomorphism of partial groups, and $\b=c_{\bar g}$ in that case. Since $\ca F$ is generated by such 
homomorphisms $c_g$, (*) then follows from the further observation that if also $h\in\ca L$ and 
$Y^h\leq Z\leq S$, then 
$$ 
c_g\circ c_h\circ\l=\l\mid_X\circ c_{\bar g}\circ c_{\bar h}.  
$$ 
Thus $\l\mid_S$ is fusion-preserving.  

Suppose that $\l$ is a projection. Then $\l$ maps $S$ onto $S'$ by I.5.3(c), and so 
$\l\mid_S$ maps $Ob(\ca F)$ onto $Ob(\ca F')$. Let $\bar w\in\bold W(\ca L')$. By I.4.13(e) there exists 
$w\in\bold W(\ca L)$ such that that the $\l\mid_S$-preimage of $S'_{\bar w}$ in $S$ is a subgroup of 
$S_w$. The mapping $\l^*:\bold W(\ca L)\to\bold W(\ca L')$ of free monoids, induced by $\l$, sends 
$w$ to $\bar w$, and then $\l\mid_S$ sends $c_w$ to $c_{\bar w}$. Thus, $\l\mid_S$ maps 
$Mor(\ca F)$ onto $Mor(\ca F')$. 
\qed 
\enddemo

\proclaim {Corollary 1.11} Let $\frak L$ be the category whose objects are pairs $(\ca L,S)$ where 
$\ca L$ is a locality and where $S$ is a designated Sylow $p$-subgroup of $\ca L$, and 
whose morphisms 
$$ 
(\ca L,S)@>\l>>(\ca L',S')
$$ 
are homomorphisms $\ca L@>\l>>\ca L'$ of partial groups such that $\l$ maps $S$ into $S'$. 
Let $\frak F$ be the category of fusion systems and fusion-preserving homomorphisms. There is then 
a functor $\t:\frak L\to\frak F$ which sends $(\ca L,S)$ to $\ca F_S(\ca L)$ and which sends $\l$ to 
the restriction of $\l$ to $S$. 
\qed 
\endproclaim

\vskip .2in 
\noindent 
{\bf Section 2: Stratified fusion systems} 
\vskip .1in 

We begin with a few definitions regarding partially ordered sets. 
\vskip .1in 

\definition {Definition 2.1} A poset $\S$ is {\it finite-dimensional} if there is an upper bound 
(denoted $dim(\S)$) to the lengths of monotone chains in $\S$. If $\S$ is finite-dimensional and 
$X\in\S$ then we write $dim_{\S}(X)$ for the supremum of the lengths of monotone increasing chains that 
terminate in $X$. 
\enddefinition 

A {\it sub-poset} of a poset $\S$ is by definition a subset of $\S$ with the partial ordering induced 
from $\S$. 

\definition {Definition 2.2} Let $\S$ be a poset. A {\it retraction} of $\S$ consists of a pair 
$(\Omega,\star)$, where $\Omega$ is a sub-poset of $\S$ and where $\star$ is a surjective 
morphism $\S\to\Omega$ of posets which restricts to the identity map on $\Omega$. The retraction 
$(\Omega,\star)$ is {\it finite-dimensional} if $\Omega$ is finite-dimensional. 
\enddefinition 

For the remainder of this section, $\ca F$ will be a fusion system on a group $S$. 
Denote by $Sub(S)$ the poset (partially ordered by inclusion) of subgroups of $S$. A sub-poset 
$\Omega$ of $Sub(S)$ is {\it $\ca F$-invariant} if $X^{\ca F}\sub\Omega$ for all $X\in\Omega$. Write 
$\1$ for the identity subgroup of $S$. 

\definition {Definition 2.3} Let $\ca F$ be a fusion system on $S$. A {\it stratification} on $\ca F$ 
consists of a finite-dimensional retraction $(\Omega,\star)$ of $Sub(S)$ satisfying the following 
conditions. 
\roster 

\item "{(St1)}" $\Omega$ is $\ca F$-invariant. 

\item "{(St2)}" $X\leq X^\star$ for all $X\in Sub(S)$. 

\item "{(St3)}" Each $\ca F$-isomorphism $\phi:X\to Y$ extends to an $\ca F$-isomorphism 
$\psi:X^\star\to Y^\star$. 

\endroster 
When the above conditions hold, we say also that $(\ca F,\Omega,\star)$ is a {\it stratified fusion system} 
on $S$. As $*$ is surjective by 2.2 one may also say, more simply, that $*$ is a stratification on 
$\ca F$, and that $(\ca F,\star)$ is a stratified fusion system. 
\enddefinition 

\vskip .1in 
Given a stratified fusion system $(\ca F,\Omega,\star)$ on $S$ and a subgroup $X\leq S$, write 
$dim_{\Omega}(X)$ for $dim_{\Omega}(X^\star)$. 

\proclaim {Lemma 2.4} Let $(\ca F,\Omega,\star)$ be a stratified fusion system on $S$ and let $X$ and 
$Y$ be be subgroups of $S$.  
\roster 

\item "{(a)}" $dim_\Omega(X)=dim_\Omega(Y)$ if $Y\in X^{\ca F}$. 

\item "{(b)}" If $X\leq Y\leq X^\star$ then $X^\star=Y^\star$. 

\item "{(c)}" Let $\phi:X\to Y$ be an $\ca F$-homomorphism, and let $X_0\leq X$ be a subgroup of $X$ 
such that $(X_0)^\star\leq X$. Then $(X_0)^\star\phi=(X_0\phi)^\star$. 

\endroster 
\endproclaim 

\demo {Proof} By (St3) the proof of (a) reduces to the case $X\in\Omega$, and then (a) follows from 
(St1). Assume $X\leq Y\leq X^\star$ as in (b). From $X\leq Y$ we have $X^\star\leq Y^\star$ as $\star$ 
is a morphism of posets. From $Y\leq X^\star$ we have $Y^\star\leq (X^\star)^\star=X^\star$, and thus 
(b) holds. Finally, assume the setup of (c). Then (a) yields 
$dim_\Omega((X_0)^\star\phi)=dim_\Omega((X_0\phi)^\star)$, while  
$(X_0)^\star\phi\in\Omega$ by (St1). Since $X_0\phi\leq (X_0)^\star\phi$, (c) follows. 
\qed 
\enddemo 

Recall from 1.9 that there is a largest subgroup $X\leq S$ such that $X$ is normal in $\ca F$. 
We refer to $X$ as the {\it socle} of $\ca F$, denoted $Soc(\ca F)$. Condition (St3) then yields: 

\proclaim {Lemma 2.5} $\1^\star$ is normal in $\ca F$. In particular, we have $\1^\star\leq Soc(\ca F)$.  
\qed 
\endproclaim

Henceforth, we assume that $\ca F$ comes equipped with a fixed stratification $(\Omega,\star)$, and we  
write $dim(X)$ for $dim_\Omega(X)$.

\proclaim {Lemma 2.6} Let $V\leq S$ be a subgroup of $S$, set $R=N_S(V)$, and set $D=C_S(V)$.  
For each subgroup $X\leq R$ and each subgroup $A\leq D$ set 
$$ 
X^{\bull}=(VX)^\star\cap R\ \ \text{and}\ \  A^{\ccirc}=(VX)^\star\cap D.  
$$ 
Further, set 
$$ 
N_{\Omega}(V)=\{X^{\bull}\mid X\leq R\}\ \ \text{and}\ \ C_{\Omega}(A)=\{A^{\ccirc}\mid Y\leq D\}.   
$$ 
Then:
\roster 

\item "{(a)}" $(N_{\Omega}(V),\bull)$ is a stratification on $N_{\ca F}(V)$ $($to be denoted also 
$N_{(\Omega,\star)}(V))$, and $(C_{\Omega}(V),\ccirc)$ is a stratification on $C_{\ca F}(V)$ $($to be 
denoted also $C_{(\Omega,\star)}(V))$. 

\item "{(b)}" The map $X\maps X^\star$ is an injective homomorphism of 
posets $N_{\Omega}(V)\to\Omega$, and the map $A\maps A^{\star}$ is an injective homomorphism of 
posets $C_{\Omega}(V)\to\Omega$. 

\item "{(c)}" If $dim(N_{\Omega}(V))=dim(\Omega)$ then $V\leq 1^\star$ and $R^\star=S$, and if 
$dim(N_{\Omega}(V))=dim(\Omega)$ then $Z(V)\leq 1^\star$ and $(VD)^\star=S$. 

\endroster   
\endproclaim 

\demo {Proof} Set $\Psi=N_{\Omega}(V)$ and $\ca E=N_{\ca F}(V)$. Let $X$ and $Y$ be subgroups of $R$ with 
$X\leq Y$. Then $(VX)^\star\leq (VY)^\star$, and intersecting with $R$ then yields $X^\bull\leq Y^\bull$. 
Thus $\bull$ is a morphism $Sub(R)\to\Psi$ of posets. If $(VX)^\star=(VY)\star$ then $X^{\bull}=Y^{\bull}$, 
so the map $\Psi\to\Omega$ given by $X\maps(VX)^\star$ is an injective morphism of posets, 
and thus $\Psi$ is finite-dimensional. 

We have 
$$ 
VX\leq VX^{\bull}\leq VX^\star\leq (VX)^\star, 
$$ 
and so $(VX)^\star=(VX^{\bull})^\star$ by 2.4(b). Then 
$$ 
(X^{\bull})^{\bull}=(VX^{\bull})^\star\cap R=(VX^\star)\cap R=X^{\bull}, 
$$ 
and thus $(\Psi,\bull)$ is a retraction. The condition (St2) for $(\Psi,\bull)$ is immediate from the 
definition of $\bull$. Now let $X$ and $Y$ be any subgroups of $R$ (not necessarily with $X\leq Y$), and let 
$\l:X\to Y$ be an $\ca E$-isomorphism. By the definition of $\ca E$, $\l$ extends to an $\ca F$-isomorphism 
$\w{\l}:VX\to VY$ which leaves $V$ invariant, and then $\w{\l}$ extends to an $\ca F$-isomorphism 
$\mu:(VX)^\star\to (VY)^\star$ with $V\mu=V$. The restriction of $\mu$ to $X^{\bull}$ is then an 
$\ca E$-isomorphism $X^{\bull}\to Y^{\bull}$, and so $(\Psi,\bull)$ satisfies (St3). Moreover, we have 
$$ 
\align
Y^{\bull}&=(X\l)^{\bull}=R\cap (VX\w{\l})^\star=R\cap((VX)^\star)\mu \\
         &=(R\cap (VX)^\star)\w{\l}=X^{\bull}\mu, 
\endalign 
$$ 
which shows that $\Psi$ is $\ca E$-invariant. That is, $(\Psi,\bull)$ satisfies (St3), and thus 
(a) holds for $(N_{\Omega}(V),\bull)$. If instead we take $\Psi=C_{\Omega}(V)$ and replace 
$\bull$ with $\ccirc$ then nearly the same argument shows that (a) holds for $(C_{\Omega}(V),\ccirc)$. 
Indeed, the only change in wording that is required here is that $\w{\l}$ extends to an $\ca F$-isomorphism 
$\mu:(VX)^\star\to (VY)^\star$ such that $\mu$ restricts to the identity map on $V$.

Let $(V^{\bull}=X_0<\cdots<X_d=R)$ be a monotone chain of maximal length in $N_{\Omega}(V)$. If 
either $V\nleq\1^\star$ or $R^\star\neq S$ then the chain 
$$ 
\1^\star\leq V^\star<\cdots<R^\star\leq S
$$ 
contains a monotone chain of length $d+1$ in $\Omega$, and hence $dim(N_{\Omega}(V))<dim(\Omega)$. 
This completes the proof of (b) and (c) for $N_{\Omega}(V)$. 
The proof of (b) and (c) for $C_{\Omega}(V)$ is again nearly identical, and is omitted. 
\qed 
\enddemo 

\definition {Definition 2.7} Let $(\ca F,\Omega,\star)$ be a stratified fusion system on $S$
and let $V\leq S$ be a subgroup of $S$. Then $V$ is {\it fully normalized} 
in $(\ca F,\Omega,\star)$ if $dim(N_S(U))\leq dim(N_S(V))$ for all $U\in V^{\ca F}$; and   
$V$ is {\it fully centralized} in $(\ca F,\Omega,\star)$ if $dim(C_S(U)U)\leq dim(C_S(V)V)$ for all 
$U\in V^{\ca F}$. 
\enddefinition

Recall (1.7) that a set $\G$ of subgroups of $S$ is  $\ca F$-closed if $\G$ is non-empty, 
$\ca F$-invariant, and closed with respect to overgroups in $S$. 

\definition {Definition 2.8} Let $\ca F$ be a fusion system on $S$ and let $Y\leq S$ be a subgroup of 
$S$. Then $Y$ is {\it normalizer-inductive} in $\ca F$ if: 
\roster 

\item "{($\cdot$)}" For each $X\in Y^{\ca F}$ there exists an $\ca F$=homomorphism $\phi:N_S(X)\to N_S(Y)$ 
with $X\phi=Y$. 

\endroster 
$Y$ is {\it centralizer-inductive} in $\ca F$ if: 
\roster 

\item "{($\cdot$)}" For each $X\in Y^{\ca F}$ there exists an $\ca F$=homomorphism 
$\psi:C_S(X)X\to C_S(Y)Y$ with $X\psi=Y$. 

\endroster 
Let $\G$ be an $\ca F$-closed set of subgroups of $S$. Then $\ca F$ is {\it $\G$-inductive} if  
for each $X\in\G$ there exists $Y\in X^{\ca F}$ such that $Y$ is normalizer-inductive in $\ca F$. 
$\ca F$ is {\it inductive} if $\ca F$ is $Sub(S)$-inductive, where $Sub(S)$ is the set of all 
subgroups of $S$.  
\enddefinition 

\proclaim {Lemma 2.9} Let $(\ca F,\Omega,\star)$ be a stratified fusion system on $S$, let $\G$ be an 
$\ca F$-closed set of subgroups of $S$, and assume that $\ca F$ is $\G$-inductive. Let $Y\in\G$. 
\roster 

\item "{(a)}" $Y$ is normalizer-inductive in $\ca F$ if and only if $Y$ is fully normalized in 
$(\ca F,\Omega,\star)$. 

\item "{(b)}" $Y$ is centralizer-inductive in $\ca F$ if and only if $Y$ is fully centralized in 
$(\ca F,\Omega,\star)$. 

\item "{(c)}" If $Y$ is fully normalized in $(\ca F,\Omega,\star)$ then $Y$ is fully centralized in 
$(\ca F,\Omega,\star)$. 

\endroster 
\endproclaim 

\demo {Proof} If $Y$ has the inductive property in $\ca F$ then clearly 
$dim_{\Omega}(N_S(X))\leq dim_{\Omega}(N_S(Y))$ for each $X\in Y^{\ca F}$, and thus $Y$ is fully 
normalized in $(\ca F,\Omega,\star)$. On the other hand, let $X\in\G$ and suppose that $X$ is fully 
normalized in $(\ca F,\Omega,\star)$. As $\ca F$ is $\G$-inductive, there exists 
an $\ca F$-conjugate $Y$ of $X$ such that $Y$ has the inductive property in $\ca F$. Thus there 
exists an $\ca F$-homomorphism $\phi:N_S(X)\to N_S(Y)$ with $X\phi=Y$. As $X$ is fully normalized we 
get $dim_{\Omega}(N_S(X))=dim_{\Omega}(N_S(Y))$. Then $\phi$ extends to an $\ca F$-isomorphism 
$\phi^*:N_S(X)^*\to N_S(Y)^*$, and $\phi^*$ then restricts to an isomorphism 
$$ 
N_{N_S(X)^*}(X)\to N_{N_S(Y)^*}(Y). 
$$ 
But $N_{N_S(X)^*}(X)=N_S(X)$ and $N_{N_S(Y)^*}(Y)-N_S(Y)$, and so $\phi$ is an isomorphism 
$N_S(X)\to N_S(Y)$. It follows that $X$ has the inductive property in $\ca F$, and this completes the 
proof of (a). The proof of (b) is similar in every respect, since 
$$ 
C_{(C_S(X)X)^*}(X)=C_S(X)
$$ 
for any subgroup $X\leq S$; and the details are left to the reader. Point (c) is immediate from (a) and (b). 
\qed 
\enddemo

In the interest of simplification of language and notation, we shall assume from now on that the fusion 
system $\ca F$ is endowed with a fixed stratification $(\Omega,\star)$. In view of 2.9, there will be no 
ambiguity in saying ``$V$ is fully normalized in $\ca F$" or ``$V$ is fully centralized in $\ca F$" if 
$V\in\G$ and $\ca F$ is $\G$-inductive. Similarly, if $X$ is a subgroup of $N_S(V)$ then the phrase 
``$X$ is fully normalized in $N_{\ca F}(V)$" should be understood relative to the stratification 
$N_{\Omega}(V)$, and we write $dim_V(X)$ for the dimension of $X$ relative to 
$N_{\Omega}(V)$. Similarly, we may speak of a subgroup of $C_S(V)$ being fully normalized in 
$C_{\ca F}(V)$.

\proclaim {Lemma 2.10} Let $\ca F$ be $\G$-inductive,   
let $T$ be strongly closed in $\ca F$, let $Q$ be a subgroup of $S$, and set $V=Q\cap T$. 
Assume $V$ is in $\G$, $V$ is fully normalized in $\ca F$, and $Q$ is fully normalized in 
$N_{\ca F}(V)$. Then $Q$ is fully normalized in $\ca F$.  
\endproclaim 

\demo {Proof} Let $P\in Q^{\ca F}$ such that $P$ is fully normalized in $\ca F$, let $\phi:N_S(Q)\to N_S(P)$ 
be an $\ca F$-homomorphism with $Q\phi=P$, and set $U=V\phi$. Then $U=P\cap T$, and $N_S(P)\leq N_S(U)$. 
As $V$ is fully normalized in $\ca F$ there exists an $\ca F$-homomorphism $\psi:N_S(U)\to N_S(V)$ 
with $U\psi=V$. Set $Q'=P\psi$. Then $Q'=(Q\phi)\psi$ is an $N_{\ca F}(V)$-conjugate of $Q$. 
Since $N_S(Q)=N_{N_S(V)}(Q)$ (and similarly for $Q'$), and since $Q$ 
is fully normalized in $N_{\ca F}(V)$, we obtain 
$$ 
dim_V(N_S(Q'))\leq dim_V(N_S(Q)).  
$$ 
Then 2.6(b) yields 
$$ 
dim(N_S(Q'))\leq dim(N_S(Q)).  
$$ 
The sequence 
$$ 
N_S(Q)@>\phi>>N_S(P)@>\psi>>N_S(Q')  
$$ 
of injective homomorphisms then shows that $dim(N_S(P))=dim(N_S(Q))$, and so $Q$ is fully normalized in 
$\ca F$. 
\qed 
\enddemo

\proclaim {Lemma 2.11} Let $\ca F$ be $\G$-inductive and let $X\in\G$. Suppose that $X^{\star}$ is fully 
normalized in $\ca F$. Then there exists an $N_{\ca F}(X^\star)$-conjugate $X'$ of $X$ such that 
$X'$ is fully normalized in $\ca F$. Moreover, for any such $X'$ we have $X^\star=(X')^\star$, and 
$X'$ is fully normalized in $N_{\ca F}(X^*)$. 
\endproclaim 

\demo {Proof} Let $Y\in X^{\ca F}$ be fully normalized in $\ca F$ and let $\a:X\to Y$ be an 
$\ca F$-isomorphism. Then $\a$ extends to an $\ca F$-isomorphism $\b:X^*\to Y^*$. As $X^*$ 
is fully normalized in $\ca F$ there exists an $\ca F$-homomorphism $\phi:N_S(Y^*)\to N_S(X^*)$ with 
$Y^*\phi=X^*$. Set $X'=Y\phi$. Then $X'$ is an $N_{\ca F}(X^*)$-conjugate of $X$, and since  
$N_S(Y)\leq N_S(Y^*)$ we have $N_S(Y)\phi\leq N_S(X')$. Then $dim(N_S(X'))=dim(N_S(Y))$ and $X'$ is 
fully normalized in $\ca F$. We have $X^\star=(X_1)^\star$ for any $N_{\ca F}(X^*)$-conjugate $X_1$ of $X$, 
by 2.4(b). Any $\ca F$-homomorphism $\eta:N_S(X_1)\to N_S(X')$ with $X_1\eta=X'$ 
is then an $N_{\ca F}(X^*)$-homomorphism, and so $X'$ is fully normalized in $N_{\ca F}(X^*)$. 
\qed
\enddemo 

The hypothesis in 2.11, that $X$ is in $\G$, can be weakened if $N_{\ca F}(X^*)$ is inductive. 

\proclaim {Lemma 2.12} Let $\ca F$ be $\G$-inductive, and let $X\leq S$ be a subgroup of $S$ such that 
$X^\star\in\G$. Assume that $X^\star$ is fully normalized in $\ca F$ and that $X$ is fully normalized in 
$N_{\ca F}(X^\star)$. Assume further that $N_{\ca F}(X^\star)$ is inductive. Then $X$ is fully 
normalized and normalizer-inductive in $\ca F$. 
\endproclaim 

\demo {Proof} Let $Y\in X^{\ca F}$. Then $Y^\star\in(X^\star)^{\ca F}$ by (St3). As $X^\star\in\G$  
and $\ca F$ is $\G$-inductive, there exists an $\ca F$-homomorphism $\phi:N_S(Y^\star)\to N_S(X^\star)$ 
with $Y^\star\phi=X^\star$. Set $X_1=Y\phi$. Then $X_1$ is an $N_{\ca F}(X^\star)$-conjugate of $X$. 
By hypothesis $N_{\ca F}(X^\star)$ is inductive and $X$ is fully normalized in $N_{\ca F}(X^\star)$, 
so there exists an 
$N_{\ca F}(X^\star)$-homomorphism $\psi:N_S(X_1)\to N_S(X)$ with $X_1\psi=X$. Let $\phi_1$ be the 
restriction of $\phi$ to $N_S(Y)$. Then $\phi_1\circ\psi$ maps $N_S(Y)$ into $N_S(X)$, and sends 
$Y$ to $X$. As the dimension function is $\ca F$-equivariant by 2.4(a), we conclude that 
$dim(N_S(Y))\leq dim(N_S(X))$, and thus $X$ is fully normalized in $\ca F$. 
\qed 
\enddemo 

\proclaim {Lemma 2.13} Let $\ca F$ be $\G$-inductive, and let $V\in\G$. 
Then there exists an $\ca F$-conjugate $U$ of $V$ such that both $U$ and $Soc(N_{\ca F}(U))$ are 
fully normalized in $\ca F$. 
\endproclaim 

\demo {Proof} We may assume to begin with that $V$ is fully normalized in $\ca F$/ 
Set $\w V=Soc(N_{\ca F}(V))$. Then $V\leq\w V$, so $\w V\in\G$. Let $\w U\in\w V^{\ca F}$ 
be fully normalized in $\ca F$, let $\phi:N_S(\w V)\to N_S(\w U)$ be an $\ca F$-homomorphism which 
maps $\w V$ to $\w U$, and set $U=V\phi$. Set $R=N_S(V)$ and $Q=N_S(U)$. As $V\in\G$ and $V$ is fully 
normalized in $\ca F$ there exists an $\ca F$-homomorphism $\psi:Q\to R$ with $U\psi=V$. Let $\l:R\to R$ 
be the composition of $\phi\mid_R$ followed by $\psi\mid_{R\phi}$. Then $\l$ is an 
$N_{\ca F}(V)$-endomorphism of $R$. Here $R\in N_{\Omega}(V)$ (cf. 2.6). Since  
$$ 
dim_{N_{\Omega}(V)}(R)=dim_{N_{\Omega}(V)}(R\l) 
$$ 
it follows that $\l$ is surjective, and hence an automorphism of $R$. Then $\phi$ restricts to an 
isomorphism $R\to Q$, and 1.3 implies that $\phi$ induces an isomorphism $N_{\ca F}(V)\to N_{\ca F}(U)$ 
of fusion systems. Then $\w U=Soc(N_{\ca F}(U))$, and thus $Soc(N_{\ca F}(U))$ is fully normalized in 
$\ca F$. Since $dim(Q)=dim(R)$, also $U$ is fully normalized in $\ca F$. 
\qed 
\enddemo

\proclaim {Lemma 2.14} Let $\ca F$ be inductive, and let $\ca E$ be a fusion subsystem of 
$\ca F$ on $S$. Suppose that $Aut_{\ca E}(V)=Aut_{\ca F}(V)$ for all $V$ such that $V$ is fully normalized 
in $\ca F$. Then $\ca E=\ca F$. 
\endproclaim 

\demo {Proof} Suppose $\ca E\neq\ca F$. Then there exists an $\ca E$-isomorphism $\phi:X\to Y$ such 
that $\phi$ is not an $\ca F$-isomorphism. Among all such, choose $\phi$ so that $dim(X)$ is as 
large as possible. Then $X\in\Omega$, and $X\neq S$ since $Aut_{\ca E}(S)=Aut_{\ca F}(X)$. As $\ca F$ 
is inductive there exist $\ca F$-homomorphisms $\l:N_S(X)\to N_S(V)$ and $\r:N_S(Y)\to N_S(V)$ such 
that $X\l=V=Y\r$. As $dim(N_S(X))>dim(X)$ the maximality of $dim(X)$ implies that 
$\l:N_S(X)\to Im(\l)$ is an $\ca E$-isomorphism, and similarly for $\r:N_S(Y)\to Im(\r)$. Let 
$\l_0$ and $\r_0$ be the restrictions of $\l$ and $\r$ to $X$ and $Y$, respectively, and set  
$\a=\l_0\i\circ\phi\circ\r_0$. Then $\a$ is an $\ca F$-automorphism of $V$, hence an $\ca E$-automorphism. 
Then $\phi=\l_0\circ\a\circ\r_0\i$ is an $\ca E$-isomorphism. 
\qed 
\enddemo

\vskip .2in 
\noindent 
{\bf Section 3: The fusion system of a centric locality} 
\vskip .1in

We begin by recalling a few notions and results from section 2 of Part I. Let $(\ca L,\D,S)$ be a 
locality. For each $w\in\bold W(\ca L)$ there is then (I.2.9) a conjugation isomorphism  
$$ 
c_w:S_w\to S_{w\i}.  
$$ 
If $w$ is the empty word we take $S_w=S$ and $c_w$ to be the identity map on $S$. Otherwise $S_w$ is the 
set of all $g\in S$ such that $g$ is conjugated successively into $S$ by the entries of $w$. The fusion system 
$\ca F_S(\ca L)$ is by definition the fusion system on $S$ generated by the set of all such mappings $c_w$. 

Definition I.2.9 also introduced the mapping 
$$ 
X\maps X^\star=\bigcap\{S_w\mid X\leq S_w,\ w\in\bold W(\ca L)\} 
$$ 
on the set $Sub(S)$ of subgroups of $S$, and introduced the poset 
$$ 
\Omega=\Omega_S(\ca L)=\{X^\star\mid X\leq S\}. 
$$ 
Then 
$$ 
\Omega=\{S_w\mid w\in\bold W(\ca L)\}=\{X\leq S\mid X=X^\star\},  
$$
by I.2.13 and I.2.11(a). 

As defined following I.3.4, a subgroup $X$ of $S$ is {\it fully normalized} in $\ca F$ if 
$dim_\Omega(N_S(X))\geq dim_\Omega(N_S(X'))$ for every $\ca F$-conjugate $X'$ of $X$. The following lemma 
shows, among other things, that $(\Omega_S(\ca L),\star)$ is a stratification on $\ca F$, and hence that 
the definition following I.3.4 is consistent with that of 2.7.  

\proclaim {Lemma 3.1} Let $(\ca L,\D,S)$ be a locality, and set $\ca F=\ca F_S(\ca L)$. Then 
$(\Omega_S(\ca L),\star)$ is a stratification on $\ca F$, and $\ca F$ is $\D$-inductive. 
Indeed, for any $Q\in\D$ such that $Q$ is fully normalized in $\ca F$, and any $\ca F$-conjugate $P$ of 
$Q$, there exists $g\in\ca L$ with $N_S(P)\leq S_g$ and with $P^g=Q$. 
\endproclaim 

\demo {Proof} Set $\Omega=\Omega_S(\ca L)$. Then $(\Omega,\star)$ is a retraction of $Sub(S)$ by I.2.11(a). 
It is  immediate from the definition of $X^\star$ that $X\leq X^\star$, while I.2.12 shows that $\Omega$ 
is $\ca F$-invariant and that every $\ca F$-isomorphism $X\to Y$ extends to an $\ca F$-isomorphism 
$X^\star\to Y^\star$. Thus $(\Omega,\star)$ satisfies the conditions (St1) through (St3) in the 
definition 2.3 of stratification, as required. 

Let $Q\in\D$ such that $Q$ is fully normalized in $(\ca F,\Omega,\star)$, and let $P\in Q^{\ca F}$ be an 
$\ca F$-conjugate of $P$. Then $Q=P^g$ for some $g\in\ca L$ by I.2.3(c), and conjugation by $g$ is 
an isomorphism $N_{\ca L}(P)\to N_{\ca L}(Q)$ by I.2.3(b). Also, as $Q$ is fully normalized we have 
$N_S(Q)\in Syl_p(N_{\ca L}(Q))$ by I.3.10, and so there exists $x\in N_{\ca L}(Q)$ such that 
$(N_S(P)^g)^x\leq N_S(Q)$. Here $(g,x)\in\bold D$ via $P$, so we may replace $g$ with $gx$ and 
obtain $N_S(P)^g\leq N_S(Q)$. 
\qed 
\enddemo

The stratification $(\Omega_S(\ca L),\star)$ will be said to be {\it induced from} $\ca L$. 
\vskip .1in

\proclaim {Lemma 3.2} Let $(\ca L,\D,S)$ be a locality, set $\ca F=\ca F_S(\ca L)$, and let 
$P\in\D$ with $P$ fully normalized in $\ca F$. Then $N_{\ca F}(P)=\ca F_{N_S(P)}(N_{\ca L}(P))$,  
and $C_{\ca F}(P)=\ca F_{C_S(P)}(C_{\ca L}(P))$.  
\endproclaim 

\demo {Proof} As $P$ is fully normalized in $\ca F$ we have $N_S(P)\in Syl_p(N_{\ca L}(P))$ by I.3.10. 
Then $C_S(P)\in Syl_p(C_{\ca L}(P))$ by I.3.13(a), and it follows that $P$ is fully centralized in $\ca F$.  
Let $\phi:X\to Y$ be an $N_{\ca F}(P)$-isomorphism between two subgroups $X$ and $Y$ of $N_S(P)$ 
containing $P$. As $P\in\D$ we get $\phi=c_g$ for some $g\in N_{\ca L}(P)$ by I.2.3(c), and this shows that 
$N_{\ca F}(P)=\ca F_{N_S(P)}(N_{\ca L}(P))$. 

Next, let $\psi:U\to V$ be a $C_{\ca F}(P)$-isomorphism 
between two subgroups $U$ and $V$ of $C_S(P)$. By the definition of $C_{\ca F}(P)$, $\psi$ extends to an 
$\ca F$-isomorphism $\w\psi:UP\to VP$ such that $\w\psi$ restricts to the identity map on $P$. Then 
$\w\psi=c_g$ for some $g\in C_{\ca L}(P)$, and thus $C_{\ca F}(P)=\ca F_{C_S(P)}(C_{\ca L}(P))$. 
\qed 
\enddemo

\definition {Definition 3.3} Let $\ca F$ be a fusion system system on $S$. A subgroup $X\leq S$ is 
{\it centric} in $\ca F$ (or is {\it $\ca F$-centric}) if $C_S(Y)\leq Y$ for all $Y\in X^{\ca F}$. 
The set of all $\ca F$-centric subgroups of $S$ is denoted $\ca F^c$. 
\enddefinition 

\definition {Definition 3.4} Let $(\ca L,\D,S)$ be a locality on $\ca F$. Then $\ca L$ is {\it centric} if 
$$ 
\bold D(\ca L)=\{w\in\bold W(\ca L)\mid S_w\in\ca F^c\}. 
$$  
\enddefinition 

Notice that if $(\ca L,\D,S)$ is a centric locality on $\ca F$ then also $(\ca L,\ca F^c,S)$ is a centric 
locality on $\ca F$. Thus, we may always take the set of objects of a centric locality on $\ca F$ 
to be $\ca F^c$. 
\vskip .1in

For the remainder of this section we shall be working with the following setup.

\definition {Hypothesis 3.5} $(\ca L,\D,S)$ is a centric locality on $\ca F$, and 
$(\ca F,\Omega,\star)$ is the stratification induced from $\ca L$. Further, we are given an 
$\ca F$-closed set $\G$ of subgroups of $S$ such that $\D\sub\G$ and such that $(\ca F,\Omega,\star)$ is 
$\G$-inductive. We adopt the linguistic and notational abbreviations introduced following 2.9.  
\enddefinition

\proclaim {Lemma 3.6} Assume 3.5 and let $V\in\G$. 
\roster 

\item "{(a)}" If $V$ is either fully centralized or fully normalized in $\ca F$, and $C_S(V)\leq V$, then 
$V\in\ca F^c$. 

\item "{(b)}" If $V$ is fully normalized in $\ca F$, and $X$ is a subgroup of $N_S(V)$ containing $V$, 
then $X$ is centric in $\ca F$ if and only if $X$ is centric in $N_{\ca F}(V)$. 

\item "{(c)}" If $V$ is fully centralized in $\ca F$, and $X$ is a subgroup of $C_S(V)$ containing $Z(V)$, 
then $XV$ is centric in $\ca F$ if and only if $X$ is centric in $C_{\ca F}(V)$. 

\endroster 
\endproclaim 

\demo {Proof} Point (a) is immediate from 2.9. Now let $V$ and $X$ be given as in (b). If $X$ is centric 
in $\ca F$ then $X$ is centric in any subsystem of $\ca F$ on any subgroup of $S$ containing $X$, so 
$X$ is centric in $N_{\ca F}(V)$. On the other hand, suppose that $X$ is centric in $N_{\ca F}(V)$. In order 
to show that $X$ is centric in $\ca F$ we are free to replace $X$ with any $N_{\ca F}(V)$-conjugate of $X$, 
and so we may assume that $X$ is fully centralized in $N_{\ca F}(V)$. Let $Y\in X^{\ca F}$ with $Y$ 
fully normalized in $\ca F$. Then $Y\in\G$ and there exists an $\ca F$-homomorphism 
$\phi:C_S(X)X\to C_S(Y)Y$ with $X\phi=Y$. Here $C_S(Y)Y\leq N_S(V\phi)$, so there exists an 
$\ca F$-homomorphism $\eta:C_S(Y)Y\to N_S(V)$ with $V(\phi\circ\eta)=V$. Set $X'=X(\phi\circ\eta)$. 
Then $X'$ is an $N_{\ca F}(V)$-conjugate of $X$, and since $X$ is centric in $N_{\ca F}(V)$ we then 
have $C_S(X')\leq X'$. Then $C_S(Y)\leq Y$, and so $Y\in\ca F^c$ by (a). This completes the 
proof of (b). The proof of (c) is similar, and is 
left to the reader. 
\qed 
\enddemo 

\proclaim {Lemma 3.7} Assume Hypothesis 3.5 and let $V\in\G$, with $V$ fully normalized in $\ca F$. Set 
$$ 
\D_V=\{P\in\D\mid V\norm P\},  
$$ 
and set 
$$ 
\bold D_V=\{w\in \bold W(N_{\ca L}(V))\mid N_{S_w}(V)\in\D_V\}.   
$$ 
Let $\ca L_V$ be the set of words of length $1$ in $\bold D_V$, and let $\ca F_V$ be the fusion system on 
$N_S(V)$ whose morphisms are the conjugation maps $c_w:X\to Y$ where $X$ and $Y$ are subgroups of 
$N_S(V)$ and where $w\in\bold W(\ca L_V)$ with $X\leq S_w$ and with $X^w\leq Y$. Then the following hold. 
\roster 

\item "{(a)}" $\ca L_V$ is a partial group via the restriction of the product in $\ca L$ to 
$\bold D_V$, and via the restriction of the inversion map associated with $\ca L$; and 
in this way the inclusion map $\ca L_V\to\ca L$ is a homomorphism of partial groups. 

\item "{(b)}" $(\ca L_V,\D_V,N_S(V))$ is a locality on $\ca F_V$. 

\item "{(c)}"Let $(\Omega_V,\sharp)$ be the stratification on $\ca F_V$ induced from $\ca L_V$. Then 
$X=X^\star\cap N_S(V)$ for all $X\in\Omega_V$, and the map $X\maps X^\star$ is an injective homomorphism of 
posets $\Omega_V\to\Omega$. Moreover, if $V\nleq \1^\star$ or if $N_S(V)^\star\neq S$, then 
$dim(\Omega_V)<dim(\Omega)$. 

\item "{(d)}" If $\ca F_V=N_{\ca F}(V)$ then $(\ca L_V,\D_V,N_S(V))$ is a centric locality. 

\endroster
\endproclaim 

\demo {Proof} Point (a) is a straightforward exercise with definition I.1.1, and the proof is omitted.  
Observe that $\D_V$ is an $N_{\ca F}(V)$-closed set of subgroups of $N_S(V)$. Since $\ca F_V$ is a fusion 
subsystem of $N_{\ca F}(V)$ it is then immediate from the definition of $\bold D_V$ that $(\ca L_V,\D_V)$ 
is objective, as defined in I.2.1. 

As $V$ is fully normalized in $\ca F$, I.3.9 implies that $N_S(V)$ is a maximal $p$-subgroup 
of $N_{\ca L}(V)$, and hence also a maximal $p$-subgroup of $\ca L_V$. 
By I.2.15 $\Omega_V$ is finite-dimensional, so we have (b). 

Set $R=N_S(V)$ and let $X$ be a subgroup of $R$. Then $X^\sharp=R\cap S_w$ for some $w\in\bold W(\ca L_V)$ 
by I.2.13, and so $R\cap X^\star\leq X^\sharp$. As $X\leq X^\star\cap R$ we then have $X=R\cap X^\star$ in 
the case that $X=X^\sharp$. Now let $X\leq Y\leq R$ with both $X=X^\sharp$ and $Y=Y\sharp$. If 
$X^\star=Y^\star$ then $X=R\cap X^\star=R\cap Y^\star=Y$, and this shows that the map $X\maps X^\star$ is 
an injective homomorphism $\Omega_V\to\Omega$. 

By I.3.7 there is a largest normal $p$-subgroup $O_p(\ca L)$ of $\ca L$, and $O_p(\ca L)=\1^\star$. 
Let $O_p(\ca L_V)<\cdots<Y$ be a monotone chain in $\Omega_V$. If $V\nleq O_p(\ca L)$ then there is a 
longer monotone chain 
$$ 
O_p(\ca L)< O_p(\ca L_V)^\star<\cdots Y^\star
$$ 
in $\Omega$. Similarly, if $N_S(V)^\star\neq S$ then 
$$ 
O_p(\ca L_V)^\star<\cdots Y^\star<S 
$$ 
is a longer monotone chain in $\Omega$, and so $dim(\Omega_V)<dim(\Omega)$. This yields (c). 

Point (d) follows from the observation, given by 3.6(b), that if $\ca F_V=N_{\ca F}(V)$ then 
$\bold D(\ca L_V)=\{w\in\bold W(\ca L_V)\mid N_{S_w}(V)\in N_{\ca F}(V)^c\}$. 
\qed 
\enddemo

The proof of the following result is essentially the same as the proof of points (a), (b), and (d) of 
the preceding lemma (with 3.6(c) invoked in place of 3.6(b)). The details are left to the reader.

\proclaim {Lemma 3.8} Assume Hypothesis 3.5 and let $V\in\G$, with $V$ fully centralized in $\ca F$. Set 
$$ 
\S_V=\{Q\leq C_S(V)\mid Z(V)\leq Q\ \text{and}\ QV\in\D\}. 
$$ 
and set 
$$ 
\bold E_V=\{w\in \bold W(C_{\ca L}(V))\mid C_{S_w}(V)V\in\S_V\}.   
$$ 
Let $\ca C_V$ be the set of words of length 1 in $\bold E_V$. Let $\ca E_V$ be the fusion system on 
$C_S(V)$ whose morphisms are the conjugation maps $c_w:X\to Y$ where $X$ and $Y$ are subgroups of 
$C_S(V)$ and where $w\in\bold W(\ca C_V)$ with $X\leq S_w$ and with $X^w\leq Y$. Then the following hold. 
\roster 

\item "{(a)}" $\ca C_V$ is a partial group via the restriction of the product in $\ca L$ to 
$\bold E_V$, and via the restriction of the inversion map associated with $\ca L$; and 
in this way the inclusion map $\ca C_V\to\ca L$ is a homomorphism of partial groups. 

\item "{(b)}" $(\ca C_V,\S_V,C_S(V))$ is a locality on $\ca E_V$.   

\item "{(c)}" Assume that $\ca E_V=C_{\ca F}(V)$. Then 
$(\ca C_V,\S_V,C_S(V))$ is a centric locality. 

\endroster 
\qed 
\endproclaim

The remainder of this section will be devoted to the proof of the following result.

\proclaim {Theorem 3.9} Assume Hypothesis 3.5. 
\roster 

\item "{(a)}" $\ca F$ is inductive. 

\item "{(b)}" If $V$ is fully normalized in $\ca F$ then the locality $\ca L_V$ given by 3.7 is centric, 
on $N_{\ca F}(V)$. 

\item "{(c)}" If $V$ is fully centralized in $\ca F$ then the locality $\ca C_V$ given by 3.8 is centric, 
on $C_{\ca F}(V)$. 

\endroster 
\endproclaim

We assume Hypothesis 3.5 throughout the remainder of this section, and we begin the proof of 3.9 by making 
a few initial reductions. 

\proclaim {Lemma 3.10} Suppose that 3.9(b) holds relative to $\G$. That is, suppose that for all 
$V\in\G$ with $V$ fully normalized in $\ca F$ the locality $\ca L_V$ is centric on $N_{\ca F}(V)$. 
Then also 3.9(c) holds relative to $\G$. 
\endproclaim 

\demo {Proof} Let $V\in\G$ be given with $V$ fully centralized in $\ca F$, and let $V'\in V^{\ca F}$ 
with $V'$ fully normalized in $\ca F$. As $\ca F$ is $\G$-inductive so there exists 
an $\ca F$-homomorphism $\phi:N_S(V)\to N_S(V')$ with $V\phi=V'$. Here $C_S(V)V\in\ca F^c$ by 3.6(c), 
so $C_S(V)V\in\D$ and $\phi=c_g$ for some $g\in\ca L$ with $C_S(V)V\leq S_g$. Inductivity implies that 
$V'$ is fully centralized, and then conjugation by $g$ induces an isomorphism $\ca C_V\to\ca C_{V'}$. 
From this it follows that if $\ca C_{V'}$ is a centric locality on $C_{\ca F}(V')$ then $\ca C_V$ is 
a centric locality on $C_{\ca F}(V)$, and thus we may assume from the outset that $V$ is fully 
normalized in $\ca F$. Then 3.9(b) yields the centric locality $\ca L_V$ on $N_{\ca F}(V)$. 

Notice that the definition of $\S_V$ in 3.8 involves $\ca L$, but remains unchanged when $\ca L$ is 
replaced by $\ca L_V$. For that reason the definition of $\ca C_V$ remains unchanged as well.   
Since $C_{\ca F}(V)=C_{N_{\ca F}(V)}(V)$ we may now assume that $\ca L=\ca L_V$. That is, we may 
take $V\norm\ca L$. 

Set $T=C_S(V)V$ and set $\ca N=C_{\ca L}(V)V$. Then $\ca N\norm\ca L$ and $T=S\cap\ca N$. Set 
$H=N_{\ca L}(T)$. Then $H$ is a subgroup of $\ca L$ since $T\in\ca F^c$. 
Let $u=(f_1,\cdots,f_n)\in\bold D(\ca C_V)$, let $h\in H$, and for each $i$ with $1\leq i\leq n$ 
set $u_i=(h\i,f_i,h)$. Set $u'=u_1\circ\cdots\circ u_n$ and set $P=C_{S_u}(V)V$. Then 
$u'\in\bold D(\ca L)$ via $(C_{S_g}(V)V)^h$, and thus $H$ acts on $\ca C_V$ as a group of automorphisms. 
Set $\ca E=\ca F_T(\ca C_V)$. Then also $H$ acts as a group of automorphisms of $\ca E$. 

For $A\leq T$ set $A^{\bull}=A^\star\cap T$, and let $\Omega_0$ be the set, partially ordered by 
inclusion, of all $A\leq T$ such that $A=A^{\bull}$. There is then an injective mapping 
$\Omega_0\to\Omega$ given by $A\maps A^\star$, and thus $\Omega_0$ is finite-dimensional. Further, 
by I.3.2 we have $A<N_T(A)$ for $T\neq A\in\Omega_0$. 

Assume by way of contradiction to 3.9(c) that $\ca E\neq C_{\ca F}(V)$. Then there exists a 
$C_{\ca F}(V)$-isomorphism $\a:X\to Y$ such that $\a$ is not an $\ca E$-isomorphism. Here $\a$ 
extends to an $\ca F$-isomorphism $\w{\a}:VX\to VY$ such that the restriction of $\w{\a}$ to $V$ is the 
identity map, and we may take $Z(V)\leq X\cap Y$. Among all such $\a$ with $\a$ not an $\ca E$-isomorphism, 
we assume that $\a$ has been chosen so that $dim(VX)$ is as large as possible. As $\w{\a}$ extends to 
an $\ca F$-isomorphism $(VX)^\star\to(VY)^\star$ we may assume $VX\in\Omega_0$. As  
$\a$ is not an $\ca E$-isomorphism we have $X\neq C_S(V)$, so $T\neq VX$. Then $VX<N_T(VX)$ by the 
result of the preceding paragraph. As $(VX)^\star\cap T=VX$ and $N_T(VX)^\star\cap T\geq N_T(VX)$ 
we obtain $(VX)^\star < N_T(VX)^\star$, and so $dim(VX)<dim(N_T(VX))$. 

Let $D\in (VX)^{\ca F}$ be fully normalized in $\ca F$, and let $\r:N_T(VX)\to N_T(D)$ be an 
$\ca F$-homomorphism with $(VX)\r=D$. Then $\r=c_w$ for some $w\bold W(\ca L)$. As $\ca L=C_{\ca L}(V)H$ 
$\r$ is then a composition of $\ca F$-isomorphisms of the form $c_g$ and $c_h$ for $g\in C_{\ca L}(V)$ 
and $h\in H$. The maximality of $dim(XV)$ in the choice of $\a$ implies that each such $c_g$ restricts 
to an $\ca E$-isomorphism. As $\ca E$ is $H$-invariant it follows that $w$ be chosen to be of the 
form $w_0\circ(h)$ where $w_0\in\bold W(C_{\ca L}(V))$ and $h\in H$. Similarly there exists an 
$\ca F$-homomorphism $\s:N_T(VY)\to D$ with $(TY)\s=D$, and such that $\s=c_{\bar w_0}\circ c_{\bar h}$ 
with $\bar{w_0}\in\bold W(C_{\ca L}(V))$ and $\bar h\in H$. 

Set $X_0=X^{w_0}$ and $Y_0=Y^{\bar w_0}$. Then $VX_0=D^{h\i}$ and $VY_0=D{\bar h\i}$. As 
$C_S(VX_0)=C_T(X_0)$, $VX_0$ is fully centralized in $\ca F$; and similarly for $VY_0$. 
Let $\b$ be the $C_{\ca F}(V)$-isomorphism $c_{(w_0)\i}\circ\a\circ c_{\bar w_0}:X_0\to Y_0$. Then $\b$ 
is not an $\ca E$-isomorphism, as otherwise $\a$ is an $\ca E$-isomorphism. Thus we may assume that 
$w_0$ and $\bar w_0$ are empty words, and that $VX$ and $VY$ are fully centralized in $\ca F$. 

We have the centric locality $\ca L_D$ on $N_{\ca F}(D)$, and we have the partial normal subgroup 
$\ca K_D:=C_{\ca L_D}(D)D\norm\ca L_D$, where $S\cap\ca K_D=C_S(D)D=C_T(D)D$. Set $R=C_T(D)D$ and set 
$G=N_{\ca L_D}(R)$. Then $R\in\ca F^c$, so $R\in\D$, and $G$ is a subgroup of $\ca L_D$. Set 
$\ca F_D=\ca F_{N_S(D)}(\ca L_D)$. The Frattini Lemma (I.4.10) yields $\ca L_D=\ca K_DG$, so every 
$\ca F_D$-automorphism of $D$ extends to an $\ca F_D$-automorphism of $R$. We have 
$\ca F_D=N_{\ca F}(D)$ by 3.9(b), so the automorphism $\g=c_{h\i}\circ\w{\a}\circ c_{\bar h}$ of $D$ 
extends to an $\ca F_D$-automorphism $\w{\g}$ of $R$. Then $\w{\g}=c_g$ for some $g\in G$, and so 
$\w{\a}=c_h\circ c_g\circ c_{{\bar h}\i}$. Here $(h,g,{\bar h}\i)\in\bold D(\ca L)$ via $R^{h\i}$, 
and we set $a=hg\bar h\i$. Thus $\w{\a}=c_a$, so $a\in C_V$, and we have shown that $\a\in\ca E_V$.  
This contradiction completes the proof. 
\qed 
\enddemo

\proclaim {Lemma 3.11} If $\ca L$ is a group then $\ca L$ is not a counterexample to 3.9. 
\endproclaim 

\demo {Proof} Assume that $\ca L$ is a group, and let $\S$ be the set of all subgroups of $S$. Then 
$(\ca L,\S,S)$ is a locality on $\ca F$. By 3.1 $\ca F$ is inductive, and we have 3.9(a).  

By 3.2, if $V$ is fully normalized in $\ca F$ then the group $N_{\ca L}(V)$ is a locality on $N_{\ca F}(V)$. 
By definition, $\bold D_V$ is the set of all $w\in\bold D$ such that $N_{S_w}(V)\in\ca F^c$.  
By 3.6(b) yields $\bold D_V=\{w\in\bold D\mid N_{S_w}(V)\in N_{\ca F}(V)^c\}$, and thus 
the locality $(N_{\ca L}(V),\D_V,N_S(V))$ is centric. Thus 3.9(b) holds, and then 3.9 holds in its 
entirety by 3.10. 
\qed 
\enddemo 

\proclaim {Lemma 3.12} Set $X=O_p(\ca F)$, and set $\ca L_X=N_{\ca L}(X)$. Then $(\ca L_X,\D,S)$ is a 
centric locality on $\ca F$, and  $dim(\Omega_S(\ca L_X))\leq dim(\Omega)$.  Moreover, if $\ca L$ is a 
counterexample to 3.9 then so is $\ca L_X$. 
\endproclaim 

\demo {Proof} We have $X\norm S$ so I.3.14 implies that $(\ca L_X,\D,S)$ is a locality. Every 
$\ca F$-isomorphism $\phi:P\to Q$ extends to an $\ca F$-isomorphism $\psi:XP\to XQ$ such that $X\psi=X$, so 
$\ca F=\ca F_S(\ca L_X)$. Then $\bold D(\ca L_X)=\{w\in\bold W(\ca L_X)\mid S_w\in\ca F^c\}$, 
and thus $(\ca L_X,\D,S)$ is a centric locality on $\ca F$. An application of 2.6 with $X$ in the role 
of $V$ yields $dim(\Omega_S(\ca L_X))\leq dim(\Omega)$.  

Assume now that $\ca L_X$ is not a counterexample to 3.9. Then $\ca F$ is inductive, and for $V$ fully 
normalized in $\ca F$ the locality $(\ca L_X)_V$ is centric on $N_{\ca F}(V)$. Notice that then 
$N_{\ca F}(V)=\ca F_{N_S(V)}((\ca L_X)_V)$ is a subsystem of $\ca F_{N_S(V)}(\ca L_V)$, and that 
$\ca F_{N_S(V)}(\ca L_V)$ is a subsystem of $N_{\ca F}(V)$. Thus $\ca L_V$ is a locality on $N_{\ca F}(V)$,   
and then 3.6(b) implies that $\ca L_V$ is a centric locality. Thus $\ca L$ is not a counterexample to 3.9 
if $\ca L_X$ is not a counterexample. 
\qed 
\enddemo

\vskip .1in 
Let $\frak G$ be the set of all $\ca F$-closed subsets $\G$ of $Sub(S)$ containing $\D$ such that 3.9   
holds relative to $\G$. That is (and in view of 3.10) for $\G\in\frak G$: 
\roster 

\item "{(1)}" $(\ca F,\Omega,\star)$ is $\G$-inductive, and   

\item "{(2)}" if $V\in\G$ is fully normalized in $\ca F$ then the locality $\ca L_V$ given by 3.7 is centric 
on $N_{\ca F}(V)$. 

\endroster 
We have $\D\in\frak G$ by 3.1 (which yields (1)) and by 3.2 and 3.6(b) (which together yield (2)). 
Thus $\frak G$ is non-empty. We regard $\frak G$ as a poset via inclusion. By Zorn's Lemma there exists a 
maximal $\G\in\frak G$, and we fix such a maximal $\G$ in what follows.  

\vskip .1in 
Recall (I.3.7) that $\1^\star$ is the largest normal $p$-subgroup $O_p(\ca L)$ of $\ca L$, and 
notice that $O_p(\ca L)\leq O_p(\ca F)$ and that $\Omega_S(N_{\ca L}(O_p(\ca F)))$ is contained in $\Omega$ 
as a sub-poset. By 3.12 we may then assume: 

\definition {3.13} Among all counterexamples to Theorem 3.9, $\ca L$ has been chosen so that 
$dim(\Omega)$ is as small as possible, and so that $O_p(\ca L)=O_p(\ca F)$. 
\enddefinition

Recall that a subgroup $V$ of $S$ is said to be normalizer-inductive in $\ca F$ if 
there exists an $\ca F$-homomorphism $N_S(U)\to N_S(V)$ which maps $U$ to $V$.

\proclaim {Lemma 3.14} Let $X\leq S$ be a subgroup of $S$. Suppose that $X^*\in\G$, $X^*$ is fully normalized 
in $\ca F$, and $X$ is fully normalized in $N_{\ca F}(X^*)$. Then $X$ is fully 
normalized in $\ca F$, and $X$ is normalizer-inductive in $\ca F$.  
\endproclaim 

\demo {Proof} As $X^*\in\G$ and $X^*$ is fully normalized in $\ca F$, we know that $X^*$ is 
normalizer-inductive in $\ca F$ and that $\ca L_{X^*}$ is a centric locality on $N_{\ca F}(X^*)$. 
Let $Y$ be an $\ca F$-conjugate of $X$. Then $Y^*$ is an $\ca F$-conjugate of $X^*$, and there is an 
$\ca F$-homomorphism $\phi:N_S(Y^*)\to N_S(X^*)$ with $Y^*\phi=X^*$. Set $X'=Y\phi$. Then $X'$ is an 
$N_{\ca F}(X^*)$-conjugate of $X$. 

Suppose now that $X\nleq O_p(\ca L)$. Then $X\nleq O_p(\ca F)$, so $dim(\Omega_{X^*})<dim(\Omega)$ by 2.6.  
Then $\ca L_{X^*}$ is not a counterexample to 3.9, and then $X$ is normalizer-inductive in $N_{\ca F}(X)$. 
As $(X')^*=X^*$ there then exists an $N_{\ca F}(X^*)$-homomorphism $\psi: N_S(X')\to N_S(X)$ with 
$X'\psi=X$. Then $\phi\circ\psi$ is an $\ca F$-homomorphism $N_S(Y)\to N_S(X)$ with $(Y\phi)\psi= X$. 
Thus $X$ is normalizer-inductive in $\ca F$ in this case. 

On the other hand, suppose that $X\leq O_p(\ca L)$. Then $X^*=O_p(\ca L)$ by I.3.7. Thus 
$N_{\ca F}(X^*)=\ca F$, and so $X$ is fully normalized in $\ca F$. Set $T=O_p(\ca L)C_S(O_p(\ca L))$ 
and set $H=N_{\ca L}(T)$. Then $T\in\ca F^c$, so $T\in\D$ and $H$ is a subgroup of $\ca L$. Set 
$\ca N=C_{\ca L}(O_p(\ca L))O_p(\ca L)$. Then $\ca N\norm\ca L$ and $T=S\cap\ca N$. 
The Frattini Lemma (I.4.-10) yields $\ca L=\ca NH$, and thus 
$$ 
X^{\ca F}=X^{\ca F_S(H)}\tag*. 
$$ 
By 3.11 there exists an $\ca F_S(H)$-conjugate $Z$ of $X$ such that $Z$ has the $\ca F_S(H)$-inductive 
property. Thus there exists $h\in H$ with $N_S(X)^h\leq N_S(Z)$ and with $X^h=Z$. As $X$ is fully 
normalized in $\ca F$ we conclude that $dim(N_S(X))=dim(N_S(Z))$. Then 
$(N_S(X)^\star)^h=N_S(Z)^\star$, and so $N_S(X)^h=N_S(Z)$. 
This shows that $X$ has the $\ca F_S(H)$-inductive property, and then (*) implies that $X$ 
is normalizer-inductive in $\ca F$. In particular, we have shown that $X$ is fully normalized in $\ca F$. 
\qed 
\enddemo 

\proclaim {Corollary 3.15} If $\Omega\sub\G$ then $\ca F$ is inductive. 
\endproclaim 

\demo {Proof} Assume that $\ca F$ is not inductive and let $X\leq S$ be a subgroup of $S$ such that 
that no $\ca F$-conjugate of $X$ is normalizer-inductive in $\ca F$. 
Assuming that $\Omega\sub\G$, then $X^*\in\G$. 
Replacing $X$ by a suitable $\ca F$-conjugate, we may take $X^*$ to be fully normalized in $\ca F$, 
and then $X$ to be fully normalized in $N_{\ca F}(X^*)$. A contradiction is then given by 3.14. 
\qed 
\enddemo 

\proclaim {Lemma 3.16} The counterexample $\ca L$ may be chosen so that 3.13 holds and so that  
$\Omega\nsub\G$. 
\endproclaim 

\demo {Proof} Assume $\Omega\sub\G$. Then $\ca F$ is inductive by 3.15. Let $X\leq S$ with $X$ fully 
normalized in $\ca F$. By 3.7, applied with 
$Sub(S)$ in the role of $\G$, we have the locality $(\ca L_X,\D_X,N_S(X))$ on the fusion system 
$\ca F_X:=\ca F_{N_S(X)}(\ca L_X)$. As $\ca L$ is a counterexample to 3.9, and in view of 3.10, 
it then suffices to show that $\ca L_X$ is a centric locality on $N_{\ca F}(X)$. 
Then by 3.7(d) it suffices to show that $\ca F_X=N_{\ca F}(X)$. 

Set $R=N_S(X)$. Then $R\in\ca F^c$ by 3.6(b), and so $R\in\D$. Now 
let $X'\in X^{\ca F}$ with also $X'$ fully normalized in $\ca F$, and set $R'=N_S(X')$. There then 
exists an $\ca F$-homomorphism $\phi:R\to R'$ with $X\phi=X'$, and since $R\in\D$ we obtain $\phi=c_g$ 
for some $g\in\ca L$ with $R^g\leq R'$ and with $X^g=X'$. Similarly there exists $h\in\ca L$ with 
$(R')^h\leq R$ and with $(X')^h=X$. Set $f=gh$ (product defined via $R$). Then $f\in\ca L_X$ since 
$f\in N_{\ca L}(X)$ and since $R\leq S_f$. Set $\Omega_X=\Omega_R(\ca L_X)$ and let 
$\s=(D_0<\cdots<D_k)$ be a chain of maximal length in $\Omega_X$. Each $(D_i)^f$ is in $\Omega_X$,  
and the image of $\s$ under conjugation by $f$ is again a chain of maximal length in $\Omega_X$. Then 
$D_k=R$ and $R^f=R$. Thus $R^g=R'$, and this shows that conjugation by $g$ is an isomorphism 
$\ca L_X\to\ca L_{X'}$, and  inducing an isomorphism $\ca F_X\to\ca F_{X'}$. But also, conjugation by 
the $\ca F$-isomorphism $c_g$ induces an isomorphism $N_{\ca F}(X)\to N_{\ca F}(X')$ by 1.6, from 
which it  follows that if $\ca F_{X'}=N_{\ca F}(X')$ then 
$\ca F_X=N_{\ca F}(X)$. We are therefore free to replace $X$ by any fully normalized $\ca F$-conjugate 
of $X$, and so by 3.14 we may assume that $X$ has been chosen so that $X^*$ is fully normalized in $\ca F$. 

Suppose that $X\nleq O_p(\ca L)$. Then, as in the proof of 3.14, $\ca L_{X^*}$ is not a counterexample 
to 3.9, and so $(\ca L_{X^*})_X$ is a centric locality on $N_{N_{\ca F}(X^*)}(X)$. For any 
$w\in\bold W(\ca L)$ with $X^w=X$ we have $(X^*)^w=X^*$, so $N_{N_{\ca F}(X^*)}(X)=N_{\ca F}(X)$. 
By definition, $\bold D(\ca L_X)$ is the set of all $w\in\bold W(N_{\ca L}(X))$ such that 
$N_{S_w}(X)\in\D$. For any such $w$ we have $w\in\bold W(N_{\ca L}(X^*))$ and $N_{S_w}(Y^*)\in\D$, 
and thus $\bold D(\ca L_X)\sub\bold D(\ca L_{X^*})$. This shows that $(\ca L_{X^*})_X=\ca L_X$, and 
so $\ca F_X=N_{\ca F}(X)$, as required. We have shown: 
\roster 

\item "{(*)}" Given $X$ fully normalized in $\ca F$ with $X\notin\G$ we have $X\leq O_p(\ca L)$, and 
thus every subgroup of $S$ which is not contained in $O_p(\ca L)$ is a member of $\G$. 

\endroster 
We next show: 
\roster 

\item "{(**)}" $\ca F_X$ is inductive. 

\endroster 
Suppose false. Then $\ca L_X$ is a counterexample to 3.9. Set $\ca K=N_{\ca L_X}(O_p(\ca F_X))$ and set 
$\L=\Omega_R(\ca K)$. Then 
$$ 
dim(\L)\leq dim(\Omega_R(\ca L_X))\leq dim(\Omega), 
$$ 
by two applications of 2.6. Notice that $\ca F_X=\ca F_R(\ca K)$, and let $\G_X$ be 
the largest $\ca F_X$-closed set of subgroups of $R$ such that $\ca K$ satisfies the conclusion of 3.9 
relative to $\G_X$. (That is, $\G_X$ is the union of all the $\ca F_X$-closed sets $\S$ 
such that $\ca K$ is $\S$-inductive, and such that $(\ca K)_V$ is centric on $N_{\ca F_X}(V)$ 
for all $V\in\S$ with $V$ fully normalized in $\ca F_X$.) If $\L\sub\G_X$ then (**) is given by 
3.15. So we may assume $\L\nsub\G_X$; and thus the counterexample $\ca K$ to 3.9 satisfies  
the conditions of 3.13 and the condition required for the fulfillment of the lemma being proved. That is, 
we have the lemma in this case, with $\ca K$ in place of $\ca L$. As the lemma is assumed to be false, we 
obtain (**). 

As $\ca F_X\neq N_{\ca F}(X)$ there exists a pair $A,A'$ of subgroups of $R$ containing $X$, and an 
$N_{\ca F}(X)$-isomorphism $\a:A\to A'$ such that $\a$ is not an $\ca F_X$-isomorphism. Let  
$(\Omega_X,\bull)=N_{(\Omega,\star)}(X)$ be the stratification on $N_{\ca F}(V)$ given by 2.6, and write 
$dim_X$ for the associated dimension function.  
Then $\a$ extends to an $N_{\ca F}(X)$-isomorphism $A^\bull\to(A')^\bull$, and we may therefore 
assume that $A$ and $A'$ are in $\Omega_X$. Further, we may assume that $\a$ has been chosen so that 
$dim_X(A)$ is as large as possible.  

Let $B\in A^{\ca F_X}$ and $B'\in(A')^{\ca F_X}$ be fully normalized in $\ca F_X$. By (**)  
there are $\ca F_X$-homomorphisms $\l:N_R(A)\to N_R(B)$ and $\l':N_R(A')\to N_R(B')$ such that $A\l=B$ and 
$A'\l'=B'$. Set $\b=(\l\mid_A)\i\circ\a\circ(\l'\mid_{A'})$. Then $\b:B\to B'$ is an 
$N_{\ca F}(X)$-isomorphism, and $\b$ is not an $\ca F_X$-isomorphism. Thus we may assume 
that $A$ and $A'$ are fully normalized in $\ca F_X$. We note that since $\a$ is not an 
$\ca F_X$-isomorphism we have $A\notin\D_X$. In particular, $A\neq R$. 

Suppose first that $A\nleq O_p(\ca L)$. Then $A\in\G$ by (*).  
Let $D\in A^{\ca F}$ be fully normalized in $\ca F$. There is then an $\ca F$-homomorphism 
$\r:N_S(A)\to N_S(D)$ with $A\r=D$. Set $Y=X\r$. As $X$ is fully normalized in $\ca F$ there is an 
$\ca F$-homomorphism $\psi:N_S(Y)\to R$ with $Y\psi=X$. Set $A_1=D\psi$, and let $\xi$ be the composition of 
$\r$ restricted to $N_R(A)$ followed by $\psi$ restricted to $N_R(A)\r$. Then 
$\psi$ and $\xi\i$ induce isomorphisms 
$$ 
C_S(D)D@>\psi>>C_S(A_1)A_1@>\xi\i>>C_S(A)A, 
$$ 
and so $A$ is fully centralized in $\ca F$. Then $C_S(A)A\in\ca F^c$ by 3.6(a), and thus 
$C_S(A)A\in\D_X$. Similarly, $C_S(A')A'\in\D_X$. Thus, in order to obtain a contradiction to the assumption 
that $\a$ is not in $\ca F_X$ it suffices to show that $\a$ extends to an $N_{\ca F}(X)$-isomorphism 
$\w\a:C_S(A)A\to C_S(A')A'$. 

Let $\r_0:C_S(A)A\to C_S(D)D$ be the isomorphism given by restriction of $\r$. Similarly, there is an 
isomorphism $\s_0:C_S(A')A'\to C_S(D)D$ sending $A'$ to $D$. Let $\r_1:A\to D$ and 
$\s_1:A'\to D$ be the obvious restrictions, and set $\b=(\r_1)\i\circ\a\circ\s_1$. Then 
$\b\in Aut_{\ca F}(D)$. As $A\in\G$ we have also $D\in\G$, and so $N_{\ca F}(D)=\ca F_D$. Set 
$K=N_{\ca L_D}(C_S(D)D)$. Then $\ca L_D=C_{\ca L_D}(D)K$ by the Frattini Lemma, and so $\b$ is given 
by conjugation by an element of $K$. Thus $\b$ extends to an $\ca F$-automorphism $\w\b$ of 
$C_S(D)D$, and then $\r_0\circ\w\b\circ(\s_0)\i$ is an extension of $\a$ to an isomorphism 
$C_S(A)A\to C_S(A')A'$. As pointed out at the end of the preceding paragraph, this shows that 
$\a$ is in fact an $\ca F_X$-isomorphism. 

Finally, we consider the case where $A\leq O_p(\ca L)$. Set $\ca N=C_{\ca L}(O_p(\ca L))$, 
$E=C_S(O_p(\ca L))$, $T=O_p(\ca L)E$, and $H=N_{\ca L}(T)$. Then $T\in\ca F^c$, so $T\in\D$, and so 
$H$ is a subgroup of $\ca L$. By I.7.7 we have $O_p(\ca L)\ca N\norm\ca L$ and $S\cap O_p(\ca L)\ca N=T$, 
and the Frattini Lemma (I.4.10) then yields 
$$ 
\ca L= O_p(\ca L)\ca N H=\ca N H. 
$$ 
Thus $\a$ is given by conjugation by an element $h\in H$. Then $\a$ extends to an 
$N_{\ca F}(X)$-isomorphism of $AE\to A'E$, and we may therefore take $E$ to be contained in $A$. 
Then $E\leq O_p(\ca L)$, so $O_p(\ca L)\in\ca F^c$, and $H=\ca L$. Thus $\ca L$ is a group, and then 
3.11 shows that $\ca L$ is not a counterexample to 3.9. This contradiction completes the proof.  
\qed 
\enddemo

In view of 3.14 and 3.16 we may assume henceforth that we are given $V\in\Omega$ with $V\notin\G$, 
and such that every overgroup of $V$ in $S$ is in $\G$. As $\G$ is $\ca F$-invariant, we may 
assume that $V$ is fully normalized in $\ca F$. Set $\G^+=\G\bigcup V^{\ca F}$. Thus, $\G^+$ is 
$\ca F$-closed.

\proclaim {Lemma 3.17} $\ca F$ is $\G^+$-inductive. 
\endproclaim 

\demo {Proof} Since $\G^+=\G\cap V^{\ca F}$, it suffices to show that $V$ is normalizer-inductive in 
$\ca F$. That is, we must show:  
\roster 

\item "{(*)}" For all $U\in V^{\ca F}$ there exists $w\in\bold W(\ca L)$ such that $U^w=V$ 
and such that $N_S(U)\leq S_w$. 

\endroster 
In fact, in view of a later application (8.2, below), we shall need to establish (*) under conditions 
other than those arising strictly from Hypothesis 3.5. To that end, we make the following list of 
assumptions (where we do {\it not} assume that $\D=\ca F^c$ as in 3.5). 
\roster 

\item "{(A)}" $(\ca L,\D,S)$ is a locality on $\ca F$, with stratification $(\Omega,\star)$ 
induced from $\ca L$, and with $\D\sub\ca F^c$.  

\item "{(B)}" We are given an $\ca F$-closed set $\G$ containing $\D$, such that $\ca F$ is 
$\G$-inductive. 

\item "{(C)}" For each $X\in\G$ such that $X$ is fully normalized in $\ca F$ we are given a 
centric locality $(\ca L_X,\D_X,N_S(X))$ on $N_{\ca F}(X)$, where $\D_X$ is the set of all 
$P\in\D$ such that $X\norm P$, and where $\ca L_X$ is the set of all $g\in N_{\ca L}(X)$ such 
that $N_{S_g}(X)\in\D_X$. Moreover, if $dim(\Omega_X)< dim(\Omega)$ then $\ca L_X$ satisfies the 
conclusion of 3.9. 

\item "{(D)}" We are given $V\in\Omega$ such that $V$ is fully normalized in $\ca F$ and such that 
the set $\G^+:=\G\cup V^{\ca F}$ is $\ca F$-closed. 

\endroster 
We now begin the proof of (*). For any $U\in V^{\ca F}$ set 
$$ 
\bold W_{U,V}=\{w\in\bold W(\ca L)\mid U\leq S_w,\ U^w=V\}. 
$$ 
We assume by way of contradiction that there exists $U\in V^{\ca F}$ such that there exists no 
$w\in\bold W_{U,V}$ with $N_S(U)\leq S_w$. Among all such counterexamples, choose $U$ so that first the 
integer 
$$ 
m:=max\{dim(N_{S_w}(U))\mid w\in\bold W_{U,V}\}
$$
is as large as possible, and then so that the integer 
$$ 
n:=min\{\ell(w)\mid w\in\bold W_{U,V},\ dim(N_{S_w}(U))=m\} 
$$ 
is as small as possible, where $\ell(w)$ is the length of $w$ as a word in the alphabet $\ca L$. 
Now choose $w\in \bold W_{U,V}$ so as to attain first the maximum value $m$ and then, subject 
to that condition, the minimum value $n$. Write $w=(g_1,\cdots,g_n)$. Thus, there is a sequence 
$$ 
U@>g_1>>V_1@>g_2>>\cdots @>g_{n-1}>>V_{n-1}@>g_n>>V 
$$ 
of conjugation isomorphisms, where each $V_i$ is a subgroup of $S$. 

Set $w_0=(g_2,\cdots,g_n)$ and set $D=N_{S_w}(U)$. Then $D^{g_1}\leq N_{S_{w_0}}(V_1)$, and so 
$dim(D)\leq dim(N_{S_{w_0}}(V_1)$. If this last inequality is strict then $V_1$ is not a counter-example 
(by the maximality of $m$), and otherwise $V_1$ is not a counter-example by the minimality of $n$. 
Thus there exists $w_1\in\bold W_{V_1,V}$ with $N_S(V_1)\leq S_{w_1}$. Set $E=N_{S_{g_1}}(U)$ and 
set $w'=(g_1)\circ w_1$. Then $w'\in\bold W_{U,V}$, and $E\leq S_{w'}$ since $E^{g_1}\leq N_S(V_1)$. 
The maximality of $m$ in the choice of $w$ then yields $dim(D)=dim(E)$. 
As $D\leq E$ we obtain $D^*=E^*$, and then $N_{D^*}(U)=N_{E^*}(U)$. Thus 
$D=N_{S_{g_1}}(U)$. As $U\notin\G$, while $S_{g_1}\in\D\sub\G$, 
$U$ is a proper subgroup of $S_{g_1}$. Since $U\in\Omega$, $U$ is then a proper subgroup 
of $N_{S_{g_1}}(U)$ by I.3.2, hence a proper subgroup of $D$, and so $D\in\G$. 

Let $D_1\in D^{\ca F}$ with $D_1$ fully normalized in $\ca F$, and let $\ca L_{D_1}$ be the centric locality 
on $\ca F_{D_1}$ given by (C). If $D_1\leq O_p(\ca L)$ then $V\leq O_p(\ca L)$, and then $V=O_p(\ca L)$ 
since $V\in\Omega$. In that case we have $U=V$, the identity map on $S$ is an $\ca F$-homomorphism 
$N_S(U)\to N_S(V)$, and $V$ is not a counterexample. Thus $D_1\nleq O_p(\ca L)$, so I.2.15 yields 
$dim(\Omega_{D_1})< dim(\Omega)$, and thus (C) implies that $\ca L_{D_1}$ satisfies the conclusion of 3.9. 

As $\ca F$ is $\G$-inductive there exists $w_1\in\bold W(\ca L)$ with 
$D^{w_1}=D_1$ and with $N_S(D)\leq S_{w_1}$. Set $U_1=U^{w_1}$, and let $U_1'$ be an 
$\ca F_{D_1}$-conjugate of $U_1$ such that $U_1'$ is fully normalized in $\ca F_{D_1}$. Thus there exists  
$v\in\bold W(\ca L_{D_1})$ with $N_{N_S(D_1)}(U_1)\leq S_v$ and with $(U_1)^v=U_1'$. 
Then $N_{N_S(U)}(D)\leq S_{w_1\circ v}$. Upon replacing $w_1$ with $w_1\circ v$, we obtain: 
\roster 

\item "{(1)}" There exists $w_1\in\bold W(\ca L)$ such that $N_{N_S(U)}(D)\leq S_{w_1}$, and having the 
property that $D_1=D^{w_1}$ is fully normalized in $\ca F$ and that $U_1=U^{w_1}$ is fully normalized in 
$\ca F_{D_1}$. 

\endroster 

We have $(D^\star)^{w_1}=(D_1)^\star$ as $(\Omega,\star)$ is a stratification on $\ca F$, and 
$$ 
N_{N_S(U)}(D^\star)D^\star@>c_{w_1}>> N_{N_S(U_1)}((D_1)^\star)D_1^\star. \tag2
$$ 
Suppose that $U_1$ is not a counterexample to (*). That is, assume that there exists $w_2\in\bold W(\ca L)$ 
such that $(U_1)^{w_2}=V$ and with $N_S(U_1)\leq S_{w_2}$. Set $E=(D_1)^{w_2}$. Then 
$(D_1^\star)^{w_2}=E^\star$, and 
$$ 
N_{N_S(U_1)}(D_1^\star)D_1^\star@>c_{w_2}>>N_{N_S(V)}(E^\star)E^\star. \tag3
$$
Composing (2) and (3) yields 
$$ 
N_{N_S(U)}(D^\star)D^\star@>c_{w_1\circ w_2}>>N_{N_S(V)}(E^\star)E^\star, 
$$ 
where $U^{w_1\circ w_2}=V$. 

We now claim that $dim(N_{S_{w_1\circ w_2}}(U))>m$. Indeed, since $D$ is a proper subgroup of $N_S(U)$, it 
follows from I.3.2 that $D$ is a proper subgroup of $N_{N_S(U)}(D^\star cap N_S(U))$. Since 
$N_{D^\star}(U)=D$ we then have $N_{N_S(U)}(D)\nleq D^*$, and so $dim(N_{N_S(U)}(D)D^\star)>dim(D^\star)=m$. 
As 
$$ 
N_{N_S(U)}(D)D^\star\leq (N_{N_S(U)}(D)D)^\star=(N_{N_S(U)}(D))^\star, 
$$ 
we obtain $dim(N_{N_S(U)}(D))>m$. As $N_{N_S(U)}(D)\leq S_{w_1\circ w_2}$, the claim is thereby 
proved. As this result violates the maximality of $m$ we conclude that $U_1$ is a counter-example to 
the lemma. Thus, we may assume to begin with that $U=U_1$ and $D=D_1$. That is: 
\roster 

\item "{(4)}" $D$ is fully normalized in $\ca F$, and $U$ is fully normalized in $\ca F_D$.  

\endroster 

Set $E=D^w$ and let $w'\in\bold W$ such that $N_S(E)\leq S_{w'}$ and such that $E^{w'}=D$. Set $U_2=V^{w'}$. 
Then $U_2=U^{w\circ w'}$ while $D^{w\circ w'}=D$, and thus conjugation by $w\circ w'$ induces an 
$N_{\ca F}(D)$-isomorphism $U\to U_2$. As $D\in\G$, and $D$ is fully normalized in $\ca F$, we have 
$N_{\ca F}(D)=\ca F_D$. As $U$ is fully normalized in $\ca F_D$, and since we have already seen that 
$\ca L_D$ satisfies the conclusion of 3.9, there exists $v'\in\bold W(\ca L_D)$ such that 
$N_{N_S(D)}(U_2)\leq S_v$ and such that $(U_2)^v=U$. There is then a sequence of two conjugation maps: 
$$
N_{N_S(V)}(E)E^\star@>c_{w'}>>N_{N_S(U_2)}(D)D^\star@>c_{v'}>>N_{N_S(U)}(D)D^\star.  
$$
Set $X=N_{N_S(V)}(E^\star)^{w'\circ v'}$. Then $dim(X)=dim(N_{N_S(V)}(E)E^\star)$ and then, as in the proof 
of (4), $dim(X)>m$. Set $w''=(w'\circ v)\i$. Then $U^{w''}=V$, and $X\leq S_{w''}$. This contradicts the 
maximality of $m$, and thereby completes the proof. 
\qed 
\enddemo

\proclaim {Lemma 3.18} We have $V=O_p(N_{\ca F}(V))$.  
\endproclaim 

\demo {Proof} Suppose false, and set $\w V=O_p(N_{\ca F}(V))$. Then $\w V\in\G$. If $\w V\leq O_p(\ca L)$ 
then $V=O_p(\ca L)$ since $V\in\Omega$. Thus the assumption that $V<\w V$ implies that $\w V\nleq O_p(\ca L)$. 
By 2.13 we may assume that $\w V$ is fully normalized in $\ca F$, and then 3.7 yields the centric 
locality $\ca L_{\w V}$ which, by 3.7(c) and the minimality of $dim(\Omega)$, is not 
a counterexample to 3.9. Since $(\ca L_{\w V})_V=\ca L_V$ and $N_{N_{\ca F}(\w V)}(V)=N_{\ca F}(V)$, we 
conclude that $\ca F_V=N_{\ca F}(V)$, and $\ca L$ itself is not a counterexample. 
\qed 
\enddemo 

\proclaim {Lemma 3.19} We have $V\nleq O_p(\ca L)$. 
\endproclaim 

\demo {Proof} Suppose $V\leq O_p(\ca L)$. Then $V=O_p(\ca L)$ as $V\in\Omega$. Then $\ca L_V=\ca L$ and
$N_{\ca F}(V)=\ca F$, and thus $\ca L_V$ is a centric locality on $N_{\ca F}(V)$. With 3.17 we then 
have points (a) and (b) of 3.9 with respect to $\G^+$, and then also 3.9(c) holds with respect to $\G^+$ 
by 3.10. This contradicts the maximality in the choice of $\G$.  
\qed 
\enddemo 

Set $R=N_S(V)$, and recall once more from 2.6 that the stratification $(N_\Omega(V),\bull)$ on 
$N_{\ca F}(V)$. Then 3.18 and the definition in 2.6 yields 
$$ 
N_\Omega(V)=\{X\leq N_S(V)\mid V\leq X=R\cap X^\star\}, 
$$ 
and $X^\bullet=R\cap(VX)^\star$ for $X\leq N_S(V)$. It is with respect to this stratification that we may 
say that a subgroup $Y$ of $R$ is (or is not) fully normalized in $N_{\ca F}(V)$. 

\vskip .1in 
The following notions, which have been extracted from the proof of [Theorem 2.2 in 5A], are essential  
to the remainder of the analysis. Let $\ca P$ be the set of all pairs $(X,U)$ of subgroups 
of $S$, such that $U\in V^{\ca F}$, $X\leq N_S(U)$, and $U$ is a proper subgroup of $X$. For any 
$(X,U)\in\ca P$ define $(X,U)^{\ca F}$ to be the set of all $(X\phi,U\phi)$ with 
$\phi\in Hom_{\ca F}(X,S)$, and write $(X,U)\phi$ for $(X\phi,U\phi)$.  Set 
$$
N_S(X,U)=N_S(X)\cap N_S(U) 
$$ 
for $(X,U)\in\ca P$ (and thereby forego the possible interpretation of $N_S(X,U)$ as the set 
- necessarily empty since $dim(U)<dim(X)$ - of all $g\in S$ such that $X^g\leq U$). A pair $(X,U)\in\ca P$  
is defined to be {\it fully normalized} 
if $dim(N_S(X,U))\geq dim(N_S(X',U'))$ for all $(X',U')\in(X,U)^{\ca F}$. For any $(X,U)\in\ca P$ let 
$N_{\ca F}(X,U)$ be the fusion system on $N_S(X,U)$ generated by the set of all $\ca F$-homomorphisms 
$\phi:P\to N_S(X,U)$ such that $X\norm P$ and such that $(X,U)\phi=(X,U)$.

\proclaim {Lemma 3.20} Let $Y$ be fully normalized in $N_{\ca F}(V)$ with $V<Y$. Then the following hold. 
\roster 

\item "{(a)}" $(Y,V)$ is a fully normalized pair.  

\item "{(b)}" There exists a fully normalized pair $(X,U)\in (Y,V)^{\ca F}$ such that $X$ is fully 
normalized in $\ca F$ and $U$ is fully normalized in $N_{\ca F}(X)$. Moreover, for any such $(X,U)$ 
there exists an $\ca F$-isomorphism $\psi:N_S(X,U)\to N_S(Y,V)$ with $(X,U)\psi=(Y,V)$. 

\item "{(c)}" Let $Y'\in Y^{N_{\ca F}(V)}$. Then there exists an $N_{\ca F}(V)$-homomorphism 
$\phi:N_S(Y',V)\to N_S(Y,V)$ with $Y'\phi=Y$. 

\endroster 
\endproclaim  

\demo {Proof} Evidently $(Y,V)\in\ca P$. Let $(X,U)\in\ca P^{\ca F}$. 
As $V$ is fully normalized in $\ca F$, 3.17 yields the existence of an $\ca F$-homomorphism 
$\phi:N_S(U)\to N_S(V)$ with $V=U\phi$. Set $Y'=X\phi$. Then $N_S(X,U)\phi\leq N_S(Y',V)$, and thus 
$dim(N_S(X,U))\leq dim(N_S(Y',V))$. As $Y'$ is an $N_{\ca F}(V)$-conjugate of $Y$, where $Y$ is fully 
normalized in $N_{\ca F}(V)$, we have also $dim_V(N_S(Y',V))\leq dim_V(N_S(Y,V))$. The injection  
$N_{\Omega}(V)\to\Omega$ given by 2.6 is a homomorphism of posets, so $dim(N_S(Y',V))\leq dim(N_S(Y,V))$, 
and thus $dim(N_S(X,U))\leq dim(N_S(Y,V))$. As $(X,U)$ is fully normalized we conclude that 
$dim(N_S(X,U))= dim(N_S(Y,V))$. Thus $(Y,V)$ is a fully normalized pair, and (a) holds. 

Next, let $X\in Y^{\ca F}$ be fully normalized in $\ca F$ and let $\psi:N_S(Y,U)\to N_S(X)$ be an 
$\ca F$-homomorphism with $Y\psi=X$. Set $U'=U\psi$. As $V<Y$ we have $Y\in\G$, so also $X\in\G$. As $X$ 
is fully normalized in $\ca F$, 3.7 yields the centric locality $\ca L_X$, and then 
$N_{\ca F}(X)=\ca F_X$ as $X\in\G$. Let $U\in(U')^{\ca F_X}$ be fully normalized in $\ca F_X$. 
As $dim(V)<dim(X)$ we have $X\nleq O_p(\ca L)$, so $dim(\Omega_{N_S(X)}(\ca L_X))<dim(\Omega)$. 
Then $\ca F_X$ is inductive, and so  
there exists an $\ca F_X$-homomorphism $\l:N_S(X,U')\to N_S(X,U)$ with $U'\l=U$. Upon replacing $\psi$ 
with $\psi\circ\l$ we obtain $\psi:N_S(Y,V)\to N_S(X,U)$, and then (a) yields $dim(N_S(X,U))=dim(N_S(Y,V))$. 
Thus $(X,U)$ is a fully normalized pair. Moreover, 
$$ 
N_S(Y,V)=N_{N_S(Y,V)^\star}(Y,V)\psi=N_{N_S(X,U)^\star}(X,U)
$$ 
so that $\psi:N_S(Y,V)\to N_S(X,U)$ is an isomorphism. This completes the proof of (b). 

Finally, let $(Y',V)\in(Y,V)^{N_{\ca F}(V)}$, and let $(X,U)$ and $\psi:N_S(Y,V)\to N_S(X,U)$ be given 
as in (b). As $Y'\in X^{\ca F}$ and $X$ is fully normalized in $\ca F$ there exists an 
$\ca F$-homomorphism $\l:N_S(Y',V)\to N_S(X)$ with $Y'\l=X$. Set $U'=V\l$. Then $U'=U(\psi\i\circ\l)$ 
is an $N_{\ca F}(X)$-conjugate of $U$. We have already seen that $N_{\ca F}(X)=\ca F_X$ is inductive, 
so there is an $\ca F_X$-homomorphism $\r:N_S(X,U')\to N_S(X,U)$ 
with $U'\r=U$. Then $\l\circ\r\circ\psi\i$ is a homomorphism $N_S(Y',V)\to N_S(Y,V)$ which maps 
$(Y',V)$ to $(Y,V)$. Thus (c) holds. 
\qed 
\enddemo

\proclaim {Lemma 3.21} Let $\S$ be the set of all $R\in N_{\ca F}(V)^c$ such that $V\leq R$ and 
such that $R$ is fully normalized in $N_{\ca F}(V)$. Then $\S\sub\G$, and $N_{\ca F}(V)$ 
is generated by the union of the groups $Aut_{N_{\ca F}(V)}(R)$ for $R\in\S$. 
\endproclaim 

\demo {Proof} 
By definition, $N_{\ca F}(V)$ is generated by the set $\Phi$ of $N_{\ca F}(V)$-isomorphisms $\phi:P\to Q$ 
with $V\leq P\cap Q$. Given $\phi\in\Phi$, we shall say that a sequence $(\a_1,\cdots,\a_k)$ is a 
{\it $\S$-decomposition} of $\phi$ if there exists a sequence $\s=(R_1,\cdots,R_k)$ of members of $\S$ 
and a sequence $\pi=(P_0,\cdots,P_k)$ of $N_{\ca F}(V)$-conjugates of $P$ such that 
\roster 

\item "{(i)}" $P_0=P$ and $P_k=Q$, 

\item "{(ii)}" $\a_i\in Aut_{N_{\ca F}(V)}(R_i)$ for all $i$, 

\item "{(iii)}" each $\a_i$ restricts to an isomorphism $\b_i:P_{i-1}\to P_i$ (where $P_{i-1}$ and 
$P_i$ are subgroups of $R_i$), and 

\item "{(iv)}" $\phi=\b_1\circ\cdots\circ\b_k$. 

\endroster 
It suffices to show that each $\phi\in\Phi$ has a $\S$-decomposition, in order to prove the lemma.  

Among all $(\phi:P\to Q)\in\Phi$ such that $\phi$ has no $\S$-decomposition, choose $\phi$ so that 
$dim_{N_{\Omega}(V)}(P)$ is as large as possible. Suppose first that $P\neq V$, so that $(P,V)\in\ca P$. Let 
$Y$ be an $N_{\ca F}(V)$-conjugate of $P$ (hence also of $Q$) such that $Y$ is fully normalized in 
$N_{\ca F}(V)$. By 3.20(c) there exist $N_{\ca F}(V)$-homomorphisms $\l$ and $\r$ 
$$ 
N_S(P,V)@>\l>>N_S(Y,V)@<\r<<N_S(Q,V) 
$$ 
with $P\l=Y=Q\r$. If $P=N_S(V)$ then $\phi\in Aut_{N_{\ca F}(V)}(N_S(V))$, and the sequence $(\phi)$ 
of length 1 is itself a $\S$-decomposition of $\phi$. Thus $P$ is a proper subgroup of $N_S(V)$, and since 
$P\in N_\Omega(V)$ it follows from I.3.2 that $P<N_S(P,V)$. Thus 
$$ 
dim_{N_\Omega(\ca F)}(P)<dim_{N_\Omega(\ca F)}(N_S(P,V)). 
$$ 
The maximality of $dim_{N_{\Omega}(V)}(P)$ now implies that $\l$ has a 
$\S$-decomposition, and the same is then true of $\r$. The same is also true of the restriction 
$\l_0$ of $\l$ to $P$ and the restriction $\r_0$ of $\r$ to $Q$. Set $\a=\l_0\i\circ\phi\circ\r_0$. 
Then $\a\in Aut_{N_{\ca F}(V)}(Y)$. If $\a$ has a $\S$-decomposition then so does $\phi$, since 
$\phi=\l_0\circ\a\circ(\r_0)\i$. Thus $\a$ has no $\S$-decomposition, and in this way we have reduced 
the problem to the case where $P=Q=Y$ and where $\phi=\a$. 

By 3.20(b) there exists $(X,U)\in(Y,V)^{\ca F}$ with $X$ fully normalized in $\ca F$ and with $U$ fully 
normalized in $N_{\ca F}(X)$. If $dim(\ca L_X)=dim(\ca L)$ then $X\leq O_p(\ca L)$ by 3.7(c), and then 
$V=O_p(\ca L)=X$ as $V\in\Omega$. This contradicts 3.19, so in fact $dim(\ca L_X)<dim(\ca L)$, and so 
$\ca L_X$ satisfies the conclusion of 3.9.
Then $N_{\ca F}(X,U)$ is the fusion system of the centric locality 
$(\ca L_X)_U$. Let $\psi:Y\to X$ be an $\ca F$-homomorphism which maps $V$ to $U$, and set 
$\b=\psi\i\circ\a\circ\psi$. Then $\b\in Aut_{(\ca L_X)_U}(X)$. By definition, the fusion system 
$(\ca F_X)_U$ is generated by conjugation isomorphisms $P'\to Q'$ with $X\leq P',Q'$ and where $P'$ and 
$Q'$ are centric in $(\ca F_X)_U$. Thus $\b$ has a $\S'$-decomposition, where $\S'$ is the set of all 
$N_{\ca F}(X,U)$-centric subgroups of $N_S(X,U)$ containing $X$. Since $N_{\ca F}(X,U)$ is isomorphic to 
$N_{\ca F}(Y,V)$ by 3.20(b), and since $N_{\ca F}(Y,V)$-centric subgroups of $N_S(Y,V)$ which contain $Y$ 
are centric in $N_{\ca F}(Y)$, it follows that $\a$ has a $\S$-decomposition. This 
contradicts the choice of $\phi$, yielding the following result. 
\roster 

\item "{(*)}" Let $\phi:P\to Q$ be an $N_{\ca F}(V)$-isomorphism, where $P$ properly contains $V$. 
Then $\phi$ has a $\S$-decomposition. 

\endroster 
It remains to show that each $\g\in Aut_{\ca F}(V)$ has a $\S$-decomposition. By definition, $\g$ 
can be written as a composition of $\ca F$-isomorphisms: 
$$ 
V=V_0@>c_{g_1}>>V_1@>c_{g_2}>>\cdots@>c_{g_{n-1}}>>V_{n-1}@>c_{g_n}>>V_n=V,  
$$ 
for some $w=(g_1,\cdots,g_n)\in N_{\bold W}(V)$. 
For each $i$ there exists an $\ca F$-homomorphism $\eta_i:N_S(V_i)\to N_S(V)$ with $V_i\eta_i=V$, and 
where we may take $\eta_0=\eta_n$ to be the identity map on $N_S(V)$. For each $i$ with $1\leq i\leq n$ set 
$$ 
P_i=N_{S_{g_i}}(V_{i-1})\eta_{i-1}\ \ \text{and}\ \ Q_i=N_{S_{g_i\i}}(V_i)\eta_i. 
$$
Then $P_i$ and $Q_i$ are subgroups of $N_S(V)$, and we have the composition 
$$ 
P_i@>\eta_{i-1}\i>>N_{S_{g_i}}(V_{i-1})@>c_{g_i}>>N_{S_{g_i\i}}(V_i)@>\eta_i>>Q_i, 
$$ 
which is an $N_{\ca F}(V)$-isomorphism $\psi_i:P_i\to Q_i$. Let $\phi_i$ be the restriction of 
$\psi_i$ to $V$. Then $\phi_i\in Aut_{\ca F}(V)$, and $\g=\phi_1\circ\cdots\circ\phi_n$. If each 
$\phi_i$ has a $\S$-decomposition then so does $\g$. As $\g$ is a counter-example there is then an 
index $j$ such that $\phi_j$ has no $\S$-decomposition. Then $\psi_j$ has no $\S$-decomposition, and then 
(*) implies that $P_i=V=Q_i$ for all $i$. Thus $N_{S_{g_1}}(V)=V$, whence $S_{g_1}=V$ and $V\in\D$. As 
$\D\sub\G$, and $V\notin\G$, we have obtained the contradiction that proves the lemma. 
\qed 
\enddemo

\demo{\bf Proof of Theorem 3.9} We have $V\in\Omega$ with $V\notin\G$, $V$ is fully normalized 
in $\ca F$, and every strict overgroup of $V$ in $S$ is in $\G$. We have found in 3.17  
that $\ca F$ is $\G^+$-inductive, where $\G^+=\G\cup V^{\ca F}$. Let $\ca L_V$ be the locality given 
by 3.7 (with $\G^+$ in the role of $\G$), and set $R-N_S(V)$. Then $R$ is a maximal $p$-subgroup of 
$\ca L_V$, by $\G$-inductivity. Set $\ca F_V=\ca F_R(\ca L_V)$, and let $P\in(\ca F_V)^c$ with $V\leq P$. 
Then $P\in\ca F^c$ by 3.6(b), and hence $\ca L_V$ is a centric locality on $\ca F_V$.  Thus, in 
order to show that the conclusion of 3.9 holds relative to $\G^+$, it suffices to show (in view of 3.10) 
that $N_{\ca F}(V)$ is equal to its subsystem $\ca F_V$. By 3.21 it then suffices to show:  
$$ 
Aut_{N_{\ca F}(V)}(Y)=Aut_{\ca F_V}(Y) \tag* 
$$ 
for all $Y\in\S_V$, where $\S_V$ is the set of all $Y\in N_{\ca F}(V)^c$ such that $V\leq Y$ and 
such that $Y$ is fully normalized in $N_{\ca F}(V)$.  

Let $Y\in\S_V$. Then $Y\in\ca F^c$ by 3.6(b), and hence $Y\in\D$. Let $\b\in Aut_{N_{\ca F}(V)}(Y)$. 
Then $\b\in N_{Aut_{\ca F}(Y)}(V)$, and $\b$ is then given by conjugation by an element 
$g\in N_{\ca L}(Y)\cap\ca L_V$. That is, $\b\in Aut_{\ca F_V}(Y)$, and thus 
$Aut_{N_{\ca F}(V)}(Y)\sub Aut_{\ca F_V}(Y)$. The reverse inclusion holds since $\ca F_V$ is a fusion 
subsystem of $N_{\ca F}(V)$. Thus (*) holds, and the proof is complete. 
\qed 
\enddemo

\vskip .2in 
\noindent 
{\bf  Section 4: Expansions} 
\vskip .1in 

Recall from I.3.14 that if $(\ca L,\D,S)$ is a locality on $\ca F$, and $\D'$ is an $\ca F$-closed 
subset of $\D$, then there is a locality $(\ca L',\D',S)$, called the restriction of 
$\ca L$ to $\D'$, and defined by   
$$ 
\bold D(\ca L')=\{w\in\bold D(\ca L)\mid S_w\in\D'\}. 
$$ 
We may write $\ca L\mid_{\D'}$ for the restriction of $\ca L$ to $\D'$. 

\definition {Definition 4.1} Let $(\ca L,\D,S)$ and $(\ca L',\D',S)$ be localities. Then $\ca L'$ is an 
{\it expansion} of $\ca L$ if $\ca L$ is the restriction $\ca L'\mid_\D$ of $\ca L$ to $\D$. Such an 
expansion is {\it elementary} if there exists $R\in\D'$ such that $\D'=\D\bigcup R^{\ca F}$, where 
$\ca F=\ca F_S(\ca L)$.  

More generally, let $(\ca L,\D,S)$ and $(\ca L',\D',S')$ be localities, and let $\eta:\ca L\to\ca L'$ 
be an injective homomorphism of partial groups. Then $\ca L'$ is an {\it expansion of $\ca L$ via $\eta$}  
(or simply {\it $\eta$ is an expansion of $\ca L$}) if: 
\roster 

\item "{(E1)}" $\eta$ restricts to an isomorphism $S\to S'$, 

\item "{(E2)}" $\D\eta\sub\D'$.   

\item "{(E3)}" $Im(\eta)$ is equal to the restriction of $\ca L'$ to $\D$. 

\endroster 
Such an expansion is then {\it elementary} if there exists $R\in\D'$ such that 
$\D'=\D\eta\bigcup R^{\ca F'}$, where $\ca F'=\ca F_{S'}(Im(\eta))$.  
\enddefinition

\definition {Hypothesis 4.2} We are given a locality $(\ca L,\D,S)$ on $\ca F$, and we are given 
$R\in\Omega_S(\ca L)$ such that the following hold. 
\roster 

\item "{(1)}" $\D\bigcup R^{\ca F}$ is $\ca F$-closed. That is, each subgroup $P\leq S$ such that $P$ 
properly contains an $\ca F$-conjugate of $R$ is in $\D$. 

\item "{(2)}" $N_{\ca L}(R)$ is a subgroup of $\ca L$, and $N_S(R)\in Syl_p(N_{\ca L}(R))$.  

\endroster 
\enddefinition

Hypothesis 4.2 will remain in effect for the remainder of this section. Notice that 4.2 holds vacuously 
if $R\in\D$.

\proclaim {Theorem 4.3} Assume Hypothesis 4.2, and set $\D^+=\D\bigcup R^{\ca F}$. Then there exists an  
elementary expansion $(\ca L^+,\D^+,S)$ of $(\ca L,\D,S)$. More precisely:  
\roster 

\item "{(a)}" $\ca L$ is the restriction of $\ca L^+$ to $\D$, $N_{\ca L^+}(R)=N_{\ca L}(R)$,  
$\ca F_S(\ca L^+)=\ca F$, and the stratification on $\ca F$ induced from $\ca L^+$ is identical with 
$(\Omega,\star)$. 

\item "{(b)}" Let $(\ca L',\D',S')$ be an expansion of $\ca L$ via $\eta:\ca L\to\ca L'$, and set 
$\ca F'=\ca F_{S'}(\ca L')$. Assume that $\D'=\D\eta\bigcup(R\eta)^{\ca F'}$, and assume that $\eta$ 
restricts to an isomorphism $N_{\ca L}(R)\to N_{\ca L'}(R\eta)$ of groups. Then there is a unique 
homomorphism $\b:\ca L^+\to\ca L'$ such that $\b$ restricts to $\eta$ on $\ca L$; and $\b$ is then 
an isomorphism. 

\item "{(c)}" Let $\D^+_0$ be the set of all $P\in\D^+$ such that $P$ contains an $\ca F$-conjugate of $R$, 
and set $\D_0=\ca D_0^+\cap\D$. Then $\D_0$ and $\D_0^+$ are $\ca F$-closed, and we thereby have the 
restriction $\ca L_0^+$ of $\ca L^+$ to $\D_0^+$ and the restriction $\ca L_0$ of $\ca L$ to $\D_0$. 
Moreover, $\ca L^+$ is the pushout (in the category of partial groups) of the diagram 
$$ 
\ca L_0^+@<<<\ca L_0 @>>>\ca L 
$$ 
of inclusion maps. 

\endroster 
\endproclaim

The remainder of this section is devoted to the proof of Theorem 4.3. Along the way, 
especially in 4.14 through 4.17, we explicitly determine the partial group structure of 
$\ca L^+$, and these results will play an important computational role in section 5.

\vskip .1in 
Since $R\in\Omega$, we have $R<N_P(R)$ for every 
$P\in\D$ such that $R<P$ by I.3.3. Since $\Omega$ and $\D$ are $\ca F$-invariant, the following 
result is then immediate from 4.2(1).

\proclaim {Lemma 4.4} Let $U\in R^{\ca F}$ and let $P\in\D$ with $U\leq P$. Then $N_P(U)\in\D$. 
\qed 
\endproclaim

Let $\bY$ be the set of elements $y\in\ca L$ such that $R\leq S_y$ and 
$N_S(R^y)\leq S_{y\i}$. For $V\in R^{\ca F}$ set $\bY_V=\{y\in\bY\mid V=R^y\}$.

\proclaim {Lemma 4.5} Let $V\in R^{\ca F}$. Then $\bY_V\neq\nset$, and $N_S(V)\leq N_S(R)^y$ for each 
$y\in\bY_V$.  
\endproclaim 

\demo {Proof} We have $N_S(V)^{y\i}\leq N_S(R)$ for each $y\in\bY_V$ by the definition of $\bY_V$, and 
then $N_S(V)\leq N_S(R)^y$ for such $y$. Thus it suffices to show that $\bY_V\neq\nset$. 

By definition, each $\ca F$-isomorphism $\phi:V\to R$ can be factored as a composition of conjugation 
maps 
$$ 
V @>c_{x_1}>> V_1@>>> \cdots @>>>V_{k-1} @>c_{x_k}>>R \tag* 
$$
for some word $w=(x_1,\cdots,x_k)\in\bold W(\ca L)$ with $V\leq S_w$. Write $V^w=R$ to indicate this. Among 
all $V\in R^{\ca F}$ for which $\bY_V=\nset$, choose $V$ so that the minimum length $k$, taken over all words 
$w$ for which $V^w=R$, is as small as possible. Then, subject to this condition, choose $w$ so that 
$dim(N_{S_w}(U))$ is as large as possible. 

Suppose first that $k=1$. Thus $w=(x)$ for some $x\in\ca L$ with $V^x=R$. Set $P=N_{S_x}(V)$ and 
set $\w P=N_{N_S(V)}(P)$. Then $P\in\D$ by 4.4. Conjugation by $x$ then induces an isomorphism 
$N_{\ca L}(P)\to N_{\ca L}(P^x)$ by I.2.3(b). Thus $\w P^x$ is a $p$-subgroup of $N_{\ca L}(R)$. As 
$N_S(R)\in Syl_p(N_{\ca L}(R))$ by 4.2(2), there exists $z\in N_{\ca L}(R)$ with 
$(\w P^x)^z\leq N_S(R)$, and then $(x,z)\in\bold D$ via $P$. Thus $\w P^{xz}$ is defined and is a 
subgroup of $N_S(R)$. Since $dim(\w P)>dim(P)$ by I.3.4, and we contradict the 
maximality of $dim(N_{S_w}(V))$ in the choice of $V$ and $w$. This shows that $k>1$. 

The minimality of $k$ now yields $\bY_{V_1}\neq\nset$, where $V_1$ is defined by (*). Thus $V^{(c,d)}=R$, 
where $d$ may be chosen in $\bY_{V_1}$. Set $A=N_{S_c}(V)$. Then $A\in\D$ by 4.4, 
and $A^{c}\leq N_S(V_1)$. As $b\in\bY_{V_1}$ we then have $w\in\bold D$ via $A$. Thus $k=1$, in 
contradiction to what has already been shown, and completing the proof.  
\qed 
\enddemo 

Set $\bX=\{x\i\mid x\in\bY\}$ and set 
$$ 
\Phi=\bX\times N_{\ca L}(R)\times\bY. 
$$ 
For any $\phi=(x\i,h,y)\in\Phi$ set $U_\phi=R^x$ and $V_\phi=R^y$. Thus: 
$$ 
U_\phi@>x\i>>R@>h>>R@>y>>V_\phi, \quad\text{and}\quad N_S(U_\phi)@>x\i>>N_S(R)@<y\i<<N_S(V_\phi)
$$ 
are diagrams of conjugation maps, labelled by the conjugating elements. (In the first of these diagrams 
the conjugation maps are isomorphisms, and in the second they are maps into $N_S(R)$.)  
\vskip .1in 

We remind the reader of the fundamental fact (I.2.1) that $\bold D=\bold D_\D$. That is, a word 
$w\in\bold W(\ca L)$ is in $\bold D$ if and only if $S_w\in\D$.

\definition {(4.6)} Define a relation $\sim$ on $\Phi$ in the following way. For $\phi=(x\i,h,y)$ and 
$\bar\phi=(\bar x\i,\bar h,\bar y)$ in $\Phi$, write $\phi\sim\bar\phi$ if 
\roster 

\item "{(i)}" $U_\phi=U_{\bar\phi}$, $V_\phi=V_{\bar\phi}$, and 

\item "{(ii)}" $(\bar xx\i)h=\bar h(\bar y y\i)$. 

\endroster 
\enddefinition 

The products in 4.6(ii) are well-defined. Namely, by 4.6(i), $(\bar x,x\i)\in\bold D$ via 
$N_S(U)^{\bar x\i}$, $\bar xx\i\in N_{\ca L}(R)$, and then $(\bar xx\i,h)\in\bold D$ since 
$N_{\ca L}(R)$ is a group (4.2(2)). The same considerations apply to $(\bar y,y)$ and 
$(\bar h,\bar y y\i)$. 

One may depict the relation $\sim$ by means of a commutative diagram, as follows. 
$$ 
\CD 
U   @>x\i>> R@>g>>R@>y>>V \\ 
@|  @A \bar x x\i AA  @AA\bar y y\i A  @| \\ 
U   @>>\bar x\i> R@>>\bar g> R@>>\bar y> V 
\endCD 
$$

\proclaim {Lemma 4.7} $\sim$ is an equivalence relation on $\Phi$.
\endproclaim

\demo {Proof} Evidently $\sim$ is reflexive and symmetric. Let $\phi_i=(x_i\i,g_i,y_i)\in\Phi$ 
($1\leq i\leq 3$) with $\phi_1\sim\phi_2\sim\phi_3$. Then $R^{x_1}=R^{x_3}$ and 
$R^{y_1}=R^{y_3}$. Notice that 
$$ 
\align 
(x_3,x_2\i,x_2,x_1\i)&\in\bold D\quad\text{via $N_S(U)^{x_3\i}$ and,} \\ 
(y_3,y_2\i,y_2,y_1\i)&\in\bold D\quad\text{via $N_S(V)^{y_3\i}$}. 
\endalign 
$$ 
Computation in the group $N_{\ca L}(R)$ then yields 
$$ 
\align 
(x_3 x_1\i)g_1&=(x_3 x_2\i)(x_2 x_1\i)g_1=(x_3 x_2\i)g_2(y_2 y_1\i) \\ 
&=g_3(y_3 y_2\i)(y_2 y_1\i)=g_3(y_3 y_1\i), 
\endalign 
$$
which completes the proof of transitivity. 
\qed 
\enddemo

\proclaim {Lemma 4.8} Let $\psi\in\Phi$ and set $U=U_\psi$ and $V=V_\phi$. Let $x\in\bY_U$ and 
$y\in\bY_V$. Then there exists a unique $h\in N_{\ca L}(V)$ such that $\psi\sim(x\i,h,y)$. 
\endproclaim 

\demo {Proof} Write $\psi=(\bar x\i,\bar h,\bar y)$. Then $(\bar x,x\i)\in\bold D$ via $N_S(U)^{\bar x\i}$, 
and $\bar xx\i\in N_{\ca L}(R)$. Similarly $(\bar y,y\i)\in\bold D$ and $\bar yy\i\in N_{\ca L}(R)$. 
As $N_{\ca L}(R)$ is a subgroup of $\ca L$ we may form the product $h:=(x\bar x\i)\bar h(\bar yy\i)$ 
and obtain $(\bar xx\i)h=\bar h(\bar yy\i)$. Setting $\phi=(x\i,h,y)$, we thus have $\phi\sim\psi$. 
If $h'\in N_{\ca L}(R)$ with also $(x\i,h',y)\sim\psi$ then $h'=(x\bar x\i)\bar h(\bar yy\i)=h$. 
\qed 
\enddemo 

For future reference we record the following observation, even though it is simply part of 
the definition of the relation $\sim$.

\proclaim {Lemma 4.9} Let $C$ be a $\sim$-class of $\Phi$, let $\phi=(x\i,g,y)\in C$, and 
set $U=R^x$ and $V=R^y$. Then the pair $(U,V)$ depends only on $C$, and not on the choice of  
representative $\phi$. 
\qed 
\endproclaim

\proclaim {Lemma 4.10} Let $\ca L_0$ be the set of all $g\in\ca L$ such that $S_g$ contains an 
$\ca F$-conjugate $U$ of $R$. For each $g\in\ca L_0$ set 
$$ 
\Phi_g=\{\phi\in\Phi\cap\bold D\mid \Pi(\phi)=g\},   
$$ 
and set 
$$ 
\ca U_g=\{U\in R^{\ca F}\mid U\leq S_g\}. 
$$ 
\roster 

\item "{(a)}" $\ca U_g$ is the set of all $U_\phi$ such that $\phi\in\Phi_g$. 

\item "{(b)}" $\Phi_g$ is a union of $\sim$-classes.  

\item "{(c)}" Let $U\in\ca U_g$, set $V=U^g$, and let $x\in\bX_U$ and $y\in\bX_V$. Then 
$(x,g,y\i)\in\bold D$, $h:=xgy\i\in N_{\ca L}(U)$, and $\phi:=(x\i,h,y)\in\Phi_g$. If also  
$(x\i,h',y)\in\Phi_g$ then $h=h'$. 

\item "{(d)}" Let $\phi$ and $\psi$ in $\Phi_g$. Then  
$$ 
\phi\sim\psi\iff U_\phi=U_\psi\iff V_\phi=V_\psi. 
$$ 

\endroster 
\endproclaim 

\demo {Proof} Let $g\in\ca L_0$, let $U\in\ca U_g$, and set $V=U^g$. Let $(x\i,y)\in\bX_U\times\bY_V$ and 
set $w=(x,g,y\i)$. Then $w\in\bold D$ via 
$$
(N_{S_g}(U)^{x\i},N_{S_g}(U),N_{S_{g\i}}(V),N_{S_{g\i}}(V)^{y\i}), 
$$ 
by 4.5. Set $h=\Pi(w)$. Then $h\in N_{\ca L}(R)$ since $R^x=U$ and $R^y=V$. Set $w'=(x\i,x,g,y\i,y)$. Then 
$w'\in\bold D$ via $N_{S_w}(U)$, and $g=\Pi(w')=\Pi(x\i,h,y)$ by $\bold D$-associativity (I.1.4(b)). If 
also $(x\i,h',y)\in\bold D$ with $\Pi(x\i,h',y)=g$ then $h=h'$ by the cancellation rule (I.1.4(e)). This 
establishes (c), and shows that $\ca U_g\sub\{U_\phi\mid\phi\in\Phi_g\}$. The opposite inclusion is 
immediate (cf. I.2.3(c)), so also (a) is established. 

Let $\phi\in\Phi_g$ and let $\bar\phi\in\Phi$ with $\phi\sim\bar\phi$. Write $\phi=(x\i,h,y)$ and 
$\bar\phi=(\bar x\i,\bar h,\bar y)$, and set 
$$ 
w'=(\bar x\i,\bar x,x\i,h,y,\bar y\i,\bar y). \tag*
$$ 
As $\phi\in\bold D$ we have $N_{S_\phi}(U)\in\D$ by 4.4, and then $w'\in\bold D$ via $N_{S_\phi}(U)$. 
Now $\Pi(w')=\Pi(\phi)=g$, while also 
$$ 
\Pi(w')=\Pi(\bar x\i,\Pi((\bar xx\i),h,(y\bar y\i)),\bar y)=\Pi(\bar x\i,\bar h,\bar y)=\Pi(\bar\phi), 
\tag**
$$ 
and thus $\bar\phi\in\Phi_g$. This proves (b). 

It remains to prove (d). So, let $\phi,\psi\in\Phi_g$. If $\phi\sim\psi$ then $U_\phi=U_\psi$ and 
$V_\phi=V_\psi$ by 4.9. On the other hand, assume that $U_\phi=U_\psi$ or that $V_\phi=V_\psi$. Then 
both equalities obtain, since $V_{\phi}=(U_{\phi})^g$ and $V_{\psi}=(U_{\psi})^g$. Write 
$\phi=(x\i,h,y)$ and $\psi=(\bar x\i,\bar h,\bar y)$ in the usual way, and define $w'$ as in (*). 
Then $w'\in\bold D$ via $N_{S_\phi}(V_\phi)$, and $\Pi(w')=\Pi(\phi)=g$. Set 
$$ 
h'=\Pi((\bar xx\i)h(y\bar y\i)). 
$$ 
Then $g=\Pi(w')=\Pi(\bar x\i,h',\bar y)$, so $(\bar x\i,h',\bar y)\in\Phi_g$. Then 4.8 yields 
$h'=\bar h$, and thus $(\bar x\i,h',\bar y)=\bar\phi$. This shows that $\phi\sim\bar\phi$, completing 
the proof of (d). 
\qed 
\enddemo

\definition {(4.11)} We now have a partition of the disjoint union $\ca L\bigsqcup\Phi$ (and a 
corresponding equivalence relation $\approx$ on $\ca L\bigsqcup\Phi$) by means of three types of 
$\approx$-classes, as follows.   
\roster 

\item "{$\cdot$}" Singletons $\{f\}$, where $f\in\ca L$ and where $S_f$ contains no 
$\ca F$-conjugate of $R$ (classes whose intersection with $\Phi$ is empty). 

\item "{$\cdot$}" $\sim$-classes $[\phi]$ such that $[\phi]\cap\bold D=\nset$  
(classes whose intersection with $\ca L$ is empty).  

\item "{$\cdot$}" Classes $\Phi_g\cup\{g\}$, where $g\in\ca L$ and where $S_g$ contains an 
$\ca F$-conjugate of $R$ (classes having a non-empty intersection with both $\ca L$ and $\Phi$). 

\endroster 
For any element $E\in\ca L\cup\Phi$, write $[E]$ for the $\approx$-class of $E$. (Thus $[E]$ is also 
a $\sim$-class if and only if $[E]\cap\bold D=\nset$.) Let $\ca L^+$ be the set of all $\approx$-classes. 
Let $\ca L_0^+$ be the subset of $\ca L^+$, consisting of those $\approx$-classes whose  
intersection with $\Phi$ is non-empty. That is, the members of $\ca L_0^+$ are the $\approx$-classes 
of the form $[\phi]$ for some $\phi\in\Phi$. 
\enddefinition

By I.1.1 there is an ``inversion map" $w\maps w\i$ on $\bold W(\ca L)$, given by 
$(g_1,\cdots,g_n)\i=(g_n\i,\cdots,g_1\i)$. The following result is then a straightforward 
consequence of the definitions of $\Phi$, $\sim$, and $\approx$.

\proclaim {Lemma 4.12} The inversion map on $\bold W(\ca L)$ preserves $\Phi$. Further,  
for each $E\in\ca L\cup\Phi$ we have $E\i\in\ca L\bigcup\Phi$, and $[E\i]$ is the 
set $[E]\i$ of inverses of members of $[E]$. 
\qed 
\endproclaim

\definition {Definition 4.13} For any $\g=(\phi_1,\cdots,\phi_n)\in\bold W(\Phi)$ let $w_\g$ be the 
word $\phi_1\circ\cdots\circ\phi_n$ in $\bold W(\ca L)$. Let $\G$ be the set of all $\g\in\bold W(\Phi)$ 
such that  $S_{w_\g}$ contains an $\ca F$-conjugate of $R$. Let $\bold D_0^+$ be the set of all sequences 
$w=([\phi_1],\cdots,[\phi_n])\in\bold W(\ca L_0^+)$ for which there exists a sequence 
$\g$ of representatives for $w$ with $\g\in\G$. We shall say that $\g$ is a {\it $\G$-form} of $w$. 
\enddefinition

The following lemma shows how to define a product $\Pi_0^+:\bold D_0^+\to\ca L_0^+$.

\proclaim {Lemma 4.14} Let $w=([\phi_1],\cdots,[\phi_n])\in\bold D_0^+$, and let 
$\g=(\phi_1,\cdots,\phi_n)$ be a $\G$-form of $w$. Write $\phi_i=(x_i\i,h_i,y_i)$. 
\roster 

\item "{(a)}" $(y_i,x_{i+1}\i)\in\bold D$ and $y_ix_{i+1}\i\in N_{\ca L}(R)$ for each $i$ with $1\leq i<n$. 

\item "{(b)}" Set 
$$ 
w_0=(h_1,y_1x_2\i,\cdots,y_{n-1}x_n\i,h_n). 
$$ 
Then $w_0\in\bold W(N_{\ca L}(R))$ and $(x_1\i,\Pi(w_0),y_n)\in\Phi$. Moreover: 

\item "{(c)}" The $\approx$-class $[x_1\i,\Pi(w_0),y_n]$ depends only on $w$, and not on 
the choice of $\G$-form of $w$. 

\endroster 
\endproclaim 

\demo {Proof} Let $U\in R^{\ca F}$ and let $x,y\in\bY_U$. Then $(y,x\i)\in\bold D$ via $N_S(U)^{y\i}$, 
and then $yx\i\in N_{\ca L}(R)$. This proves (a), and shows that $w_0\in\bold W(N_{\ca L}(R))$. As 
$N_{\ca L}(R))$ is a subgroup of $\ca L$, (b) follows. 

Let $\bar\g=(\bar\phi_1,\cdots,\bar\phi_n)$ be any $\G$-form of $w$, write 
$\bar\phi_i=(\bar x_i\i,\bar h_i,\bar y_i)$, and define $\bar w_0$ in analogy with $w_0$. Set 
$U_0=U_{\phi_1}$ and $\bar U_0=U_{\bar\phi_1}$. For each $i$ with $1\leq i\leq n$ set $U_i=V_{\phi_i}$ 
and $\bar V_i=U_{\bar\phi_i}$. Thus:  
$$ 
U_{i-1}@>x_i\i>>R@>h_i>>R@>y_i>>U_i\quad\text{and}\quad 
\bar U_{i-1}@>\bar x_i\i>>R@>\bar h_i>>R@>\bar y_i>>\bar U_i. 
$$ 

Suppose that there exists an index $j$ with $U_j\neq\bar U_j$. As $\phi_j\approx\bar\phi_j$ it follows 
from 4.9 and 4.10 that $\phi_j$ and $\bar\phi_j$ are in $\bold D$, and that there is an element 
$g_j\in\ca L_0$ such that $\phi_j$ and $\bar\phi_j$ are in $\Phi_g$. If $j<n$ then  
$$ 
U_{j+1}=(U_j)^g\neq(\bar U_j)^g=\bar U_{j+1},  
$$ 
and if $j>0$ one obtains $U_{j-1}\neq\bar U_{j-1}$ in similar fashion, by consideration of $w\i$ 
via 4.13. Thus, for each index $i$ 
we have $U_i\neq\bar U_i$. Then $\phi_i$ and $\bar\phi$ lie in distinct $\sim$-classes by 4.9. As  
$\phi_i\approx\bar\phi_i$ it follows that $\phi_i$ and $\bar\phi_i$ are members of $\Phi\cap\bold D$, 
and that there is a word $v=(g_1,\cdots,g_n)\in\bold W(\ca L)$ with $g_i=\Pi(\phi_i)=\Pi(\bar\phi_i)$. 
As $\<U_0,\bar U_0\>\leq S_v$, we have $v\in\bold D$. 

Set $P=S_v$ and set $P_0=N_P(U_0)$. Then $P_0\in\D$ by 4.4. Set $P_i=(P_{i-1})^{g_i}$ for $1\leq i\leq n$.
Then 
$$
P_{i-1}^{x_i\i}\leq N_S(R)\ \ \text{and}\ \ (P_i)^{y_i\i}\leq N_S(R).  
$$ 
Thus conjugation by $g_i$ maps $P_{i-1}^{x_i\i}$ to $(P_i)^{y_i\i}$, and this shows: 
\roster 

\item "{(*)}" Set $w_\g=\g_1\circ\cdots\circ\g_n$. Then $w_\g\in\bold D$ via $P_0$. 

\endroster 
Similarly, one obtains $w_{\bar\g}\in\bold D$ via $N_P(\bar U_0)$, where 
Now $\bold D$-associativity (I.1.4(b)) yields: 
$$ 
\Pi(x_1\i,\Pi(w_0),y_n)=\Pi(w_\g)=\Pi(v)=\Pi(w_{\bar \g})=\Pi(\bar x_1\i,\Pi(\bar w_0),\bar y_n).  
$$ 
Thus $\Pi(x_1\i,\Pi(w_0),y_n)$ and $(\bar x_1\i,\Pi(\bar w_0),\bar y_n)$ lie in the same fiber of 
$\Pi:\Phi\cap\bold D\to\ca L_0$, and thus 
$(x_1\i,\Pi(w_0),y_n)\approx(\bar x_1\i,\Pi(\bar w_0),\bar y_n)$. This reduces (c) to the following claim. 
\roster 

\item "{(**)}" Suppose that $U_i=\bar U_i$ for all $i$. Then 
$(x_1\i,\Pi(w_0),y_n)\sim(\bar x_1\i,\Pi(\bar w_0),\bar y_n)$. 

\endroster 
Among all counter-examples to (**), let $w$ be chosen so that $n$ is as small as possible. Then 
$n\neq 0$ (i.e. $w$ and $\bar w$ are non-empty words), as there is otherwise 
nothing to verify. If $n=1$ (so that $w=(\phi_1)$ and $\bar w=(\bar\phi_1)$) 
then $w_0=(h)$, $\bar w_0=(\bar h)$, and (**) follows since $\phi\sim\bar\phi$. Thus, $n\geq 2$. 

Set $\psi=(x_1\i,h_1(y_1x_2\i)h_2,y_2)$, and similarly define $\bar\psi$. Then $\psi$ and $\bar\psi$ 
are in $\Phi$. As $\phi_i\sim\bar\phi_i$ for $i=1,2$, one has the commutative diagram 
$$ 
\CD 
U_0@>x_1\i>>R@>h_1>>R@>y_1>>U_1@>x_2\i>>R@>h_2>>R@>y_2>>U_2  \\  
@|  @A{\bar x_1x_1\i}AA @AA{\bar y_1y_1\i}AA @| @A{\bar x_2x_2\i}AA @AA{\bar y_2y_2\i}A @|   \\  
U_0@>\bar x_1\i>>R@>\bar h_1>>R@>\bar y_1>>U_1@>\bar x_2\i>>R@>\bar h_2>>R@>\bar y_2>>U_2
\endCD 
$$ 
of conjugation maps. This diagram then collapses to the commutative diagram 
$$ 
\CD 
U_0 @>x_1\i>> R @>h_1(y_1x_2\i)h_2>> R @>y_2>>U_2  \\  
@|  @A{\bar x_1x_1\i}AA  @AA{\bar y_2y_2\i}A  @|   \\
U_0 @>\bar x_1\i>> R @>\bar h_1(\bar y_1\bar x_2\i)\bar h_2>> R @>\bar y_2>>U_2 
\endCD 
$$
which shows that $\psi\sim\bar\psi$. Here $(\psi,\phi_3,\cdots,\phi_n)$ and 
$(\bar\psi,\bar\phi_3,\cdots,\bar\phi_n)$ are $\G$-forms of the word 
$u=([\psi],[\phi_3],\cdots,[\phi_n])$. Set 
$$ 
u_0=(h_1(y_1x_2\i)h_2,y_2x_3\i,\cdots,y_{n-1}x_n\i,h_n), 
$$ 
and similarly define $\bar u_0$. Then $\Pi(u_0)=\Pi(w_0)$ and $\Pi(\bar u_0)=\Pi(\bar w_0)$. 
The minimality of $n$ then yields 
$$ 
(x\i,\Pi(w_0),y_n)\approx(\bar x_1\i,\Pi(\bar w_0),\bar y_n). 
$$ 
This proves (**), and thereby completes the proof of (c). 
\qed 
\enddemo

\proclaim {Proposition 4.15} There is a mapping $\Pi^+:\bold D_0^+\to\ca L_0^+$, given by 
$$ 
\Pi_0^+(\nset)=[\1,\1,\1],    
$$ 
and by 
$$
\Pi^+_0(w)=[x_1\i,\Pi(w_0),y_n], \tag* 
$$ 
on non-empty words $w\in\bold D_0^+$, where $w_0$ is given by a $\G$-form of $w$ as in 4.14(b).  
Further, there is an involutory bijection on $\ca L_0^+$ given by 
$$ 
[x\i,g,y]\i=[y\i,g\i,x]. \tag** 
$$
With these structures, $\ca L_0^+$ is a partial group. 
\endproclaim 

\demo {Proof} The reader may refer to I.1.1 for the conditions (1) through (4) defining the notion of 
partial group. Condition (1) requires that $\bold D_0^+$ contain all words of length 1 in the 
alphabet $\ca L_0^+$, and that $\bold D_0^+$ be closed with respect to decomposition. (That is, 
if $u$ and $v$ are two words in the free monoid $\bold W(\ca L_0^+)$, and the 
concatenation $u\circ v$ is in $\bold D_0^+$, then $u$ and $v$ are in $\bold D_0^+$.) Both of these   
conditions are immediate consequences of the definition of $\bold D_0^+$. 

That $\Pi_0^+$ is a well-defined mapping is given by 4.14. 
The proof that $\Pi_0^+$ satisfies the conditions I.1.1(2) 
($\Pi_0^+$ restricts to the identity map on words of length 1) and I.1.1(3): 
$$ 
\bu\circ\bv\circ\w\in\bold D_0^+\implies\Pi_0^+(\bu\circ\bv\circ\bv)=\Pi_0^+(\bu\circ\Pi_0^+(\bv)\circ\bw)
$$ 
are then straightforward, and may safely be omitted. 

The inversion map $[x\i,h,y]\maps [y\i,h\i,x]$ is well-defined by 4.12. Evidently this mapping is an 
involutory bijection, and it extends to an involutory bijection 
$$ 
(C_1,\cdots,C_n)\i=(C_n\i,\cdots,C_1\i) 
$$ 
on $\bold W(\ca L_0^+)$. It thus remains to verify I.1.1(4). That is, we must check that  
$$ 
w\in\bold D_0^+\implies w\i\circ w\in\bold D_0^+\ \ \text{and}\ \ 
\Pi_0^+(w\i\circ w)=[\1,\1,\1]. 
$$ 
In detail: take $w=(C_1,\cdots,C_n)$ and let $\g=(\phi_1,\cdots,\phi_n)$ be a $\G$-form of $w$, where  
$\phi_i$ is written as $[x_i\i,h_i,y_i)$. 
One easily verifies that $\g\i\circ\g\in\G$, and hence $w\i\circ w\in\bold D_0^+$. Now 
$$ 
\Pi_0+(w\i\circ w)=[y_n\i,\Pi(u_0),y_n], 
$$ 
where 
$$  
u_0=(g_n\i,x_ny_{n-1}\i,\cdots,x_2y_1\i,g_1\i,x_1x_1\i,g_1,y_1x_2\i),\cdots,y_{n-1}x_n\i,g_n). 
$$ 
One observes that $\Pi(u_0)=\1$, and so $\Pi_0+(w\i\circ w)=[y_n\i,\1,y_n]$. Now observe that 
$(y_n\i,\1,y_n)\equiv(\1,\1,\1)$, since $(y_n\i,\1,y_n)$ and $(\1,\1,\1)$ are in $\bold D$ and  
since $\Pi(y_n\i,\1,y_n)=\1=\Pi(\1,\1,\1)$. Thus $\Pi_0+(w\i\circ w)=[\1,\1,\1]=\Pi_0^+(\nset)$. Thus 
I.1.1(4) holds, and the proof is complete. 
\qed 
\enddemo

\proclaim {Lemma 4.16} Let $\D_0^+$ be the set of all $P\in\D^+$ such that $P$ contains an $\ca F$-conjugate 
of $R$, and set $\D_0=\D_0^+\cap\D$. Then $\D_0$ and $\D_0^+$ are $\ca F$-closed. Let $\ca L_0$ be 
the restriction of $\ca L$ to $\D_0$ and let $\l:\ca L_0\to\ca L_0^+$ be the 
mapping $g\maps\Phi_g\cup\{g\}$. Then $\l$ is an injective homomorphisms of partial groups. 
\endproclaim 

\demo {Proof} Evidently $\D_0$ and $\D_0^+$ are $\ca F$-closed. Set $\bold D_0=\bold D(\ca L_0)$. Now  
let $v=(g_1,\cdots,g_n)\in\bold D_0$, set $w=([g_1],\cdots,[g_n])$, and let $U\in R^{\ca F}$ be 
chosen so that $U\leq S_v$. By 4.6 there exists a word $\g=(\phi_1,\cdots,\phi_n)\in\bold W(\Phi)$ 
such that $\phi_i\in[g_i]$ and such that $U\leq S_{w_\g}$, where $w_\g$ is the word 
$\phi_1\circ\cdots\circ\phi_n\in\bold W(\ca L)$. Then $\g$ is a $\G$-form of $w$, and so $w\in\bold D_0^+$. 
Set $P=S_w$. The proof of the intermediary result (*) in the proof of 4.14 shows that $w_\g\in\bold D$ via 
$N_P(U)$. Write $\phi_i=(x_i\i,h_i,y_i)$, and form the word $w_0\in\bold W(N_{\ca L}(R))$ as in 4.14(b). 
Set $g=\Pi(v)$. Then $g=\Pi(w_\g)=\Pi(x_1\i,w_0,y_n)$ by $\bold D$-associativity, and thus  
$(x_1\i,w_0,y_n)\in\Phi_g$. That is, we have $[g]=[x_1\i,w_0,y_n]$. This shows that $\l$ is a homomorphism 
of partial groups. Since $\Phi_g\cap\ca L=\{g\}$ by 4.10(c), $\l$ is injective.  
\qed
\enddemo

In the Appendix B to Part I  it was shown that the category of partial groups 
has certain colimits, and that in some cases these commute with the forgetful functor to the 
category of pointed sets. In particular, by Theorem I.B.5 the pushout of 
$(\ca L_0^+@<\l<<\ca L_0@>\iota>>\ca L)$ exists, and its underlying set is formed by taking the disjoint 
union of $\ca L_0^+$ with $\ca L$, and by then identifying $\ca L_0$ with $Im(\l)$. This yields the 
following result.

\proclaim {Proposition 4.17} There is a partial group $\ca L^+$ which is a pushout of the diagram 
$$ 
(\ca L_0^+@<\l<<\ca L_0@>\iota>>\ca L) 
$$
of injective homomorphisms. Specifically, $\bold D(\ca L^+)$ is the pushout of the diagram 
$$ 
(\bold D(\ca L_0^+)@<\w\l<<\bold D(\ca L_0)@>\w\iota>>\bold D(\ca L)) 
$$ 
of maps of sets; where $\w\iota$ is the obvious inclusion map, and where $\w\l$ is the restriction to 
$\bold D(\ca L_0)$ of the map $\l^*:\bold W(\ca L_0)\to\bold W(\ca L_0^*)$ of free monoids. The product 
$\Pi^+:\bold D(\ca L^+)\to\ca L^+$ is the defined by $\Pi_0^+$ on $\bold D(\ca L^+_0)$, and by $\Pi$ 
on $\bold D(\ca L)$. 
\endproclaim

Let $\D^+$ be the union of $\D$ with the set of all subgroups $P\leq S$ such that $P$ contains an 
$\ca F$-conjugate of $R$. The following lemma prepares the way for showing that the partial group 
$(\ca L^+,\D^+)$ is objective.

\proclaim {Lemma 4.18} Let $[\phi]\in\ca L_0^+$, and let $\l:\ca L\to\ca L^+$ be the homomorphism of 
partial groups given by $g\maps[g]$. For any subset $X$ of $S$ write $[X]$ for $X\l$. 
\roster 

\item "{(a)}" Suppose that $[\phi]\cap\ca L=\nset$, let $\phi\in[\phi]$, and set $U=U_{[\phi]}$. 
Let $S_{[\phi]}$ be the set of all $a\in S$ such that $[a]^{[\phi]}$ is defined in $\ca L^+$, 
and such that $[a]^{[\phi]}\in[S]$. Then $U\norm S_{[\phi]}=S_\phi$.  

\item "{(b)}" Let $A\in R^{\ca F}$ be an $\ca F$-conjugate of $R$ such that $[A]^{[\phi]}$ is defined 
in $\ca L^+$, with $[A]^{[\phi]}\sub[S]$. Then $[A]^{[\phi]}=[B]$ for some $B\in R^{\ca F}$, and 
there exists $\phi\in[\phi]$ such that $A=U_\phi$ and $B=V_\phi$. 

\endroster 
\endproclaim 

\demo {Proof} Let $A,B$, and $[\phi]$ be given as in (b), and suppose first that $[\phi]\cap\ca L\neq\nset$. 
Then $[\phi]=\Phi_g\cup\{g\}$ for some $g\in\ca L_0$. As $\l$ is a homomorphism of partial groups, we obtain 
$[A]^{[\phi]}=[A^g]\leq [S]$, and then $A^g\leq S$ since $\l$ is injective. Now 4.6(b) shows that there exists 
$\phi\in[\phi]$ with $A=U_\phi$ and $B=V_\phi$. Thus, (b) holds in this case. On the other hand, suppose that 
$[\phi]\cap\ca L=\nset$, let $\phi\in[\phi]$, and suppose that (a) holds. Then $A\leq S_\phi$. If 
$A\neq U_\phi$ then $S_\phi\in\D$ by 4.2(1), and then $\phi\in\bold D$ and $\Pi(\phi)\in [\phi]\cap\ca L$. 
As $[\phi]\cap\ca L=\nset$ we conclude that in fact $A=U_\phi$. Since $[\phi]\i\cap\ca L=\nset$ and 
$[B]^{[\phi]\i}=[A]$, we obtain also $V=V_\phi$. Thus (b) holds provided that (a) holds. 

Assume now that $[\phi]\cap\ca L=\nset$, let $a\in S$ be given so that $[a]^{[\phi]}\in[S]$, and set 
$[b]=[a]^{[\phi]}$. Here $[\phi]$ is a $\sim$-class, and 4.9 shows that the pair $(U,V):=(U_\phi,V_\phi)$ 
is constant over all $\phi\in[\phi]$. Set $w=([\phi]\i,[a],[\phi])$. Then $w\in\bold D_0^+$ by 
hypothesis, so there is a $\G$-form $\g=(\phi\i,\psi,\bar\phi)$ of $w$. This means that, upon 
setting $w\g=\phi\i\circ\psi\circ\bar\phi$, we have $V\leq S_{w_\g}$. The uniqueness of $(U,V)$ for 
$[\phi]$ (and of $(V,U)$ for $[\phi]\i$) yields $U=U_\psi=V_\psi$, and then $a\in N_S(U)$ 
since $a=\Pi(\psi)$. Since $[b]^{[\phi]\i}=[a]$ we similarly obtain $b\in N_S(V)$. (Notice also 
that we may take $\bar\phi$ to be $\phi$.) As $x\in\bY_U$ and $y\in\bY_V$ we obtain:  
\roster 

\item "{(*)}" $a^{x\i},b^{y\i}\in N_S(R)$. 

\endroster 

Let $\phi\in[\phi]$ and write $\phi=(x\i,h,y)$. Set $\theta=(a\i x\i,\1,x a^2)$. 
As $a\in N_S(U)$ we get $xa^2\in\bold Y_U$ and $a\i x\i\in\bold X_U$, and thus 
$\theta\in\Phi$. Moreover, we have $\theta\in\bold D$ via $N_S(U)$, and $\Pi(\theta)=a$. Thus 
$\theta\in\Phi_a$, and $(\phi\i,\theta,\phi)$ is a $\G$-form of $w$. Then 
$$
\Pi^+(w)=[y\i,h\i(x(a\i x\i))((xa^2)x\i)h,y] 
$$
by the definition (4.14(b)) of $\Pi^+$. Observing now that 
$$ 
\text{$(x,a\i,x\i,x,a^2,x\i)\in\bold D$}, 
$$ 
 via $N_S(U)^{x\i}$, we obtain 
$$
[b]=\Pi^+(w)=[y\i,h\i(xax\i)h,y]=[y\i,(a^{x\i})^h,y]. 
$$ 

We now claim that $(y\i,b^{y\i},y)\in [b]$. 
In order to see this, one observes first of all that $(y\i,b^{y\i},y)\in\Phi$. Also, since 
$b^{y\i}$ normalizes $N_S(V)^{y\i}$ we have $(y\i,y,b,y\i,y)\in\bold D$ via $N_S(V)$. Then 
$\Pi(y\i,b^{y\i},y)=b$, and the claim is proved. Thus: 
$$ 
[y\i,(a^{x\i})^h,y]=[y\i,b^{y\i},y]. \tag**
$$ 
Application of $\Pi$ to both sides of (**) then yields 
$$ 
((a^{x\i})^h)^y=(b^{y\i})^y,  
$$ 
and then $(a^{x\i})^h=b^{y\i}$ by the cancellation rule in $\ca L$. As $a^{x\i},b^{y\i}\in S$ by (*), we 
conclude that $a\in S_\phi$. This completes the proof of (a), and thereby completes the proof of the lemma. 
\qed 
\enddemo

\proclaim {Proposition 4.19} $(\ca L^+,\D^+,S)$ is a locality. Moreover: 
\roster 

\item "{(a)}" $\ca L$ is the restriction $\ca L^+\mid_\D$ of $\ca L^+$ to $\D$. 

\item "{(b)}" $N_{\ca L^+}(R)=N_{\ca L}(R)$. 

\item "{(c)}" $\ca F_S(\ca L^+)=\ca F$, and the stratification $(\Omega,\star)$ on $\ca F$ induced from 
$\ca L$ is the same as that induced from $\ca L^+$. 

\endroster 
\endproclaim 

\demo {Proof} The conjugation maps $c_{[\phi]}:U\to V$ for $U\in R^{\ca F}$ and $\phi\in\Phi$ are 
$\ca F$-homomorphisms of the form $c_{x\i}\circ c_h\circ c_y$, so $\ca F_S(\ca L^+)=\ca F$. 
In order to show that $(\ca L^+,\D^+)$ satisfies the conditions (O1) and (O2) in the definition I.2.1 
of objectivity, note that condition (O2) is then the requirement that $\D^+$ be $\ca F$-closed. 
Since $\D$ is $\ca F$-closed, and $\D^+$ is given by attaching to $\D$ an $\ca F$-conjugacy class 
$R^{\ca F}$ and all overgroups in $S$ of members of $R^{\ca F}$, (O2) then holds for $(\ca L^+,\D^+)$. 
Condition (O1) requires that $\bold D^+$ be equal to $\bold D_{\D^+}$. This means: 
\roster 

\item "{(*)}" The word $w=(C_1,\cdots,C_n)\in\bold W(\ca L^+)$ is in $\bold D^+$ $\iff$ there exists 
a sequence $(X_0,\cdots,X_n)\in\bold W(\D^+)$ such that $X_{i-1}^{C_i}=X_i$ for all $i$, $(1\leq i\leq n)$. 

\endroster 
Here we need only be concerned with the case $w\in\bold W(\ca L_0^+)$ since 
$\bold D^+=\bold D_0^+\cup[\bold D]$, and since $\bold D=\bold D_\D$. The implication $\implies$ in (*) is 
then given by the definition of $\bold D_0^+$, with $X_0\in R^{\ca F}$. The reverse implication is given 
by 4.18. Thus $(\ca L^+,\D^+)$ is objective. 

Let $\w S$ be a $p$-subgroup of $\ca L^+$ containing $S$, and let $a\in N_{\w S}(S)$. As $S\in\D$ we get 
$a\notin\ca L_1^+$, so $a\in N_{\ca L}(S)$, and then $a\in S$ since $S$ is a maximal $p$-subgroup of $\ca L$. 
Thus $\w S=S$, $S$ is a maximal $p$-subgroup of $\ca L^+$, and $(\ca L^+,\D^+,S)$ is a pre-locality. The 
restriction of $\ca L^+$ to $\D$ is by definition the partial group whose product is the restriction 
of $\Pi^+$ to $\bold D_\D$, whose underlying set is the image 
of $\Pi^+\mid\bold D_\D$, and whose inversion map is inherited from $\ca L^+$. That is, (a) holds. 

Let $[\phi]\in\ca L_1^+$, let $\phi=(x\i,h,y)\in[\phi]$, and set $U=U\phi$. Then $U\norm S_{[\phi]}=S_\phi$ 
by 4.18(a), and the conjugation map $c_{[\phi]}:S_\phi\to S$ is then the composite 
$c_x\i\circ c_h\circ c_y$ applied to $S_\phi$. Thus $c_{[\phi]}$ is an $\ca F$-homomorphism, and this 
yields (c). Suppose now that $[\phi]\in N_{\ca L^+}(R)$. Then $x,h$, and $y$ are in $N_{\ca L}(R)$, 
and then $\phi\in\bold D$ as $N_{\ca L}(R)$ is a subgroup of $\ca L$. Then $[\phi]\notin\ca L_1^+$, 
and so (b) holds. 

We have seen that $\ca F_S(\ca L^+)=\ca F$, and (c) is an immediate consequence of that observation.  
In particular, $\ca L^+$ is finite-dimensional. 
In order to show that $\ca L^+$ is a locality, it remains to show that every subgroup $H$ of $\ca L^+$ is 
countable and locally finite. Here $H$ is a subgroup of $N_{\ca L^+}(P)$ for some $P\in\D^+$ by I.2.16. 
If $H\leq\ca L$ then there is nothing to show, so assume that $H\nleq\ca L$. Then 
$P\notin\D$ and $P$ contains an $\ca F$-conjugate $U$ of $R$. By 4.2(1) we have $dim(U)=dim(P)$, 
so $U=P$ (again because $U\in\Omega$). Then $H$ is conjugate in $\ca L^+$ to a subgroup of $N_{\ca L}(R)$. 
Since $N_{\ca L}(R)$ is a subgroup of $\ca L$ by 4.2(2), $H$ is then countable and locally finite. 
\qed 
\enddemo  

\proclaim {Lemma 4.20} Let $(\ca L',\D',S')$ be an expansion of $\ca L$ via $\eta:\ca L\to\ca L'$, and 
set $\ca F'=\ca F_{S'}(\ca L')$. Assume that $\D'=\D\eta\bigcup(R\eta)^{\ca F'}$, and assume that $\eta$ 
restricts to an isomorphism $N_{\ca L}(R)\to N_{\ca L'}(R\eta)$ of groups. Then there is a unique 
homomorphism $\b:\ca L^+\to\ca L'$ such that $\b$ restricts to $\eta$ on $\ca L$; and $\b$ is then 
an isomorphism. 
\endproclaim 

\demo {Proof} Set $\bold D'=\bold D(\ca L')$ and write $\Pi'$ for the product in $\ca L'$. For simplicity of 
notation, we may regard $\eta$ as an inclusion map. Thus $(\ca L',\D^+,S)$ is a locality. Let 
$\phi=(x\i,h,y)$ and $\bar\phi=(\bar x\i,\bar h,\bar y)$ be members of $\Phi$ such that $\phi\sim\bar\phi$. 
Set $U=U_\phi$. Then $U=U_{\bar\phi}$ by 4.9, and $(\bar x,x\i)\in\bold D$ via $N_S(U)^{\bar x\i}$. Then 
$\w\Pi(\bar x,x\i)=\Pi(\bar x,x\i)$ as $\ca L$ is the restriction of $\ca L'$ to $\D$. Similarly, 
we obtain $\Pi'(\bar y,y\i)=\Pi(\bar y,y\i)$. Then: 
$$ 
\Pi'(\bar x,x\i,h)=\Pi'(\bar xx\i,h)=\Pi(\bar xx\i,h)=\Pi(\bar h,\bar yy\i)=\Pi'(\bar h,\bar y,y\i).
\tag1  
$$ 
Observe that both $(\bar x\i,\bar x,x\i,h,y)$ and $(\bar x\i,\bar h,\bar y,y\i,y)$ are in 
$\bold D^+\cap\bold D'$ (via the obvious conjugates of $U$). Then (1) yields  
$$ 
\Pi'(\phi)=\Pi'(\bar x\i,\bar x,x\i,h,y)=\Pi'(\bar x\i,\bar h,\bar y,y\i,y)=\Pi'(\bar\phi).  
$$ 
If $\phi\in\bold D$ then $\Pi'(\phi)=\Pi(\phi)$, so we have shown that there is a well-defined 
mapping 
$$ 
\b:\ca L^+\to\ca L' 
$$ 
such that $\b$ restricts to $\eta$ (the identity map) on $\ca L$, and sends the element $[\phi]\in\ca L_0^+$ 
to $\Pi'(\phi)$.  

We now show that $\b$ is a homomorphism of partial groups. As $\ca L^+$ is a pushout, in the manner 
described in 4.17(c), and since the inclusion $\ca L\to\ca L'$ is a homomorphism, it suffices to 
show that the restriction $\b_0$ of $\b$ to $\ca L_0^+$ is a homomorphism. 

Let $w=([\phi_1],\cdots,[\phi_n])\in\bold D^+$, let $\g=(\phi_1,\cdots,\phi_n)$ be 
a $\G$-form of $w$, and set $w_\g=(\phi_1,\cdots,\phi_n)$. Write $\phi_i=(x_i\i,h_i,y_i)$. Then 
$\Pi^+(w)=[x_1\i,\Pi(w_0),y_n]$, where $w_0\in\bold W(N_{\ca L}(R))$ is the word 
$$ 
w_0=(h_1,y_1x_2\i,\cdots,y_{n-1}x_n\i,h_n). 
$$
given by 4.14(b). Let $\b*$ be the induced mapping $\bold W(\ca L^+)\to\bold W(\ca L')$ of free monoids. 
Then 
$$ 
\Pi'(w\b^*)=\Pi'([\phi_1]\b,\cdots,[\phi_n]\b)=\Pi'(\Pi'(\phi_1),\cdots,\Pi'(\phi_n))=\Pi'(w_\g) 
$$ 
by $\bold D'$-associativity, and  
$$ 
(\Pi^+(w))\b=[x_1\i,\Pi(w_0),y_n]\b=\Pi'(x_1\i,\Pi(w_0),y_n).       
$$ 
Set $w_0'=(h_1,y_1,x_2\i,\cdots,y_{n-1},x_n\i,h_n)$. Then $w_0'\in\bold D'$ via $U_{\phi_1}$, and 
$\Pi'(w_0')=\Pi'(w_0)=\Pi(w_0)$. Then  
$$
(\Pi^+(w))\b=(x_1\i)\circ\Pi'(w_0')\circ(y_n))=\Pi'(w_\g)=\Pi'(w\b^*),  
$$ 
and so $\b_0$ is a homomorphism. As already mentioned, this result implies that $\b$ is a homomorphism. 

Let $f\in\ca L'$ with $f\notin\ca L$. Then $S_f\notin\D$ (as $\ca L=\ca L'\mid_\D$), and 
so $S_f$ contains a unique $U\in R^{\ca F}$. Set $V=U^f$. We note that the inversion map on $\ca L$ 
is induced from the inversion map on $\ca L'$, by I.1.13. If $V\in\D$ then $f\i\in\ca L$, and then 
$f\in\ca L$. Thus $V\notin\D$, so $V$ contains a unique $\ca F$-conjugate of $R$, and $V\in R^{\ca F}$. 
By 4.5 there exist elements $x,y\in\ca L$ such that $U^x=R=V^y$ and such that 
$N_S(U)^{x\i}\leq N_S(R)\geq N_S(V)^{y\i}$. Imititating the proof of 4.6, we find that 
$(x,f,y\i)\in\bold D'$ via $R$ and that $h:=\Pi'(x,f,y\i)\in N_{\ca L'}(R)$. Then $h\in N_{\ca L}(R)$ 
by hypothesis, and $f=\Pi'(x\i,h,y)$. In particular, we have shown: 
\roster 

\item "{(*)}" Every element of $\ca L'$ is a product of elements of $\ca L$. 

\endroster 
We may now show that $\b^*$ maps $\bold D^+$ onto $\bold D'$. Thus, let 
$w'=(f_1,\cdots,f_n)\in\bold D'$ and set $X=S_{w'}$. If $X\in\D$ then $w'\in\bold D$ and 
$w'=w'\b^*$. So assume that $X\notin\bold D$. Then $X$ contains a unique $\ca F$-conjugate $U_0$ 
of $R$, and there is a sequence $(U_1,\cdots,U_n)\in\bold W(\D^+)$ given by $U_i=(U_{i-1})^f_i$. As 
seen in the preceding paragraph, we then have $U_i\in R^{\ca F}$ for all $i$, and there exists a 
sequence $\g=(\phi_1,\cdots,\phi_n)\in\bold W(\Phi)$ such that $\Pi'(\phi_i)=f_i$. Set 
$w_\g=[\phi_1]\circ\cdots\circ[\phi_n]$. Then $w_\g\in\bold D^+$ and $(w_\g)\b^*=w'$. Thus
$\b^*:\bold D^+\to\bold D'$ is surjective. That is, $\b$ is a projection, as defined in I.5.6. 

Let $g\in Ker(\b)$. If $g\in\ca L$ then $g=\1$ since $\b\mid_{\ca L}$ is the inclusion map. So assume 
that $g\notin\ca L$. Then $g=[\phi]$ for some $\phi=(x\i,h,y)\in\Phi$, and $[\phi]$ is a $\sim$-class. 
Set $U=U_\phi$ and $V=V_\phi$. As $\1=g\b=\w\Pi(\phi)$ it follows that $U=V=R$. But then $\phi\in\bold D$ 
as $N_{\ca L}(R)$ is a subgroup of $\ca L$, and so $[\phi]$ is not a $\sim$-class. Thus $Ker(\b)=\1$. 
As $\b$ is a projection, $\b$ is then an isomorphism by I.5.3(d). It follows from (*) that $\b$ is the 
unique homomorphism $\ca L^+\to\w{\ca L}$ which restricts to the identity on $\ca L$. 
\qed 
\enddemo 

Notice that Propositions 4.19 and 4.20 complete the proof of Theorem 4.3.

\vskip .2in 
\noindent 
{\bf Section 5: Elementary expansions, partial normal subgroups, and extensions of homomorphisms} 
\vskip .1in 

This section is a continuation of the preceding one, so we shall continue to assume Hypothesis 4.2. Our 
principal aim is to establish an inclusion-preserving bijection between the set of partial normal subgroups 
of $\ca L$ and the set of partial normal subgroups of the elementary expansion $\ca L^+$ of $\ca L$, 
as follows.  

\proclaim {Theorem 5.1} Assume Hypothesis 4.2, and let $(\ca L^+,\D^+,S)$ be the expansion of $\ca L$ given 
by Theorem 4.3. Let $\ca N\norm\ca L$ be a partial normal subgroup of $\ca L$, and let 
$\bold X:=\ca N^{\ca L^+}$ be the set of all $\ca L^+$-conjugates of elements of $\ca N$ (the set of all 
$f^g$ such that $f\in\ca N$, $g\in\ca L^+$ and $(g\i,f,g)\in\bold D^+$). Let $\ca N^+$ be 
partial subgroup of $\ca L^+$ generated by $\bold X$ (as in I.1.9). Then:   
\roster 

\item "{(a)}" $\ca N^+\norm\ca L^+$.  

\item "{(b)}" $\ca N^+\cap\ca L=\ca N$. 

\item "{(c)}" If $\ca M\norm\ca L^+$ is a partial normal subgroup of $\ca L^+$ such that 
$\ca M\cap\ca L=\ca N$ then $\ca M=\ca N^+$. 

\endroster 
\endproclaim

We will employ all of the notation from section 4. Thus, the reader will need to have in mind the meanings 
of $\ca L_0$, $\ca L_0^+$, $\ca L^+$, $\Phi$, $\sim$, and $\approx$. Set $T=S\cap\ca N$ and set 
$\bold X=\ca N^{\ca L^+}$. We refer the reader to sections 4 and 5 of Part I for the notion of maximal coset 
of $\ca N$ (I.4.13), and of the quotient locality $\ca L/\ca N$ (I.5.5).

\proclaim {Lemma 5.2} Let $w\in\bold D^+$, and suppose that $w\notin\bold W(\ca L)$. 
Then $S_w$ is an $\ca F$-conjugate of $R$, and $S_w=S_{\Pi^+(w)}$ if $\Pi^+(w)\notin\ca L$. 
\endproclaim 

\demo {Proof} As $w\notin\ca L$ it follows that $S_w\notin\D$. Then, as $w\in\ca L^+$, 4.2(1) implies that  
$S_w$ is an $\ca F$-conjugate of $R$. Write $w=(g_1,\cdots,g_n)$ with $g_i\in\ca L^+$, 
set $g=\Pi^+(w)$, and suppose that $g\notin\ca L$. Thus $S_g$ is an $\ca F$-conjugate of $R$, and since 
$S_w\leq S_g$ we obtain $S_w=S_g$.  
\qed 
\enddemo

\proclaim {Lemma 5.3} Assume that 5.1(b) holds. Then:  
\roster 

\item "{(a)}" Theorem 5.1 holds in its entirety. 

\item "{(b)}" Let $\r:\ca L\to\ca L/\ca N$ and $\r^+:\ca L^+/\ca N^+$ be the canonical projections, and 
let $\iota:\ca L\to\ca L^+$ be the inclusion map. Then there is a unique homomorphism 
$\eta:\ca L/\ca N\to\ca L^+/\ca N^+$ such that $\iota\circ\r^+=\r\circ\eta$. Moreover, $\eta$ is injective, 
and $\ca L^+/\ca N^+$ is an elementary expansion of $\ca L/\ca N$ via $\eta$. Further, if $S\cap\ca N\nleq R$ 
then $\eta$ is an isomorphism, while if $S\cap\ca N\leq R$ then 
$N_{\ca L/\ca N}(R\eta)=N_{\ca L^+/\ca N^+}(R\eta)$. 

\endroster 
\endproclaim 

\demo {Proof} We assume that $\ca L\cap\ca N^+=\ca N$ (5.1(b)). Set $\bold X_0=\bold X$, 
and recursively define $\bold X_n$ for $n>0$ by 
$$ 
\bold X_n=\{\Pi^+(w)\mid w\in\bold W(\bold X_{n-1})\cap\bold D^+\}. 
$$ 
Then $\ca N^+=\<\bold X\>$ is the union of the sets $\bold X_n$, by I.1.9. In order to show that 
$\ca N^+\norm\ca L^+$ it then suffices to show that each $\bold X_n$ is closed with respect to 
conjugation (whenever defined) in $\ca L^+$. 

Let $g\in\bold X_0$. Then there exists $f\in\ca N$ and $a\in\ca L^+$ such that $(a\i,f,a)\in\bold D^+$ 
and with $g=f^a$. Now let $b\in\ca L^+$ such that $(b\i,g,b)\in\bold D^+$. If $g\in\ca L$ then 
$g\in\ca N$ by 5.1(b), and then $g^b\in\bold X_0$. On the other hand, assume $g\notin\ca L$. Then 
$S_g\in R^{\ca F}$, and 5.2 yields $S_g=S_{(a\i,f,a)}$. Set $u=(b\i,g,b)$ and $v=(b\i,a\i,f,a,b)$. 
Then $S_u=S_v=(S_g)^b$, so $v\in\bold D^+$, and $g^b=\Pi^+(u)=\Pi^+(v)=f^{ab}$. Thus $g^b\in\bold X_0$ in 
any case, and so $\bold X_0$ is closed with respect to conjugation in $\ca L^+$. 

Assume now that there exists a largest integer $k$ such that $\bold X_k$ is closed with respect to 
conjugation in $\ca L^+$, and let now $f\in\bold X_{k+1}$ and $b\in\ca L^+$ such that $f^b$ is defined 
in $\ca L^+$ and is not in $\bold X_{k+1}$. By definition there exists 
$w=(f_1,\cdots,f_n)\in\bold W(\bold X_k)\cap\bold D^+$ such that $f=\Pi^+(w)$. Suppose $f\in\ca L$. Then 
$f\in\ca N$ by 5.2(b), and then $f^b\in\bold X$. As $\bold X\sub\bold X_m$ for all $m$, we have 
$f^b\in\bold X_{k+1}$ in this case. So assume that $f\notin\ca L$. Then $S_f=S_w\in R^{\ca F}$ by 
5.2. Set $w'=(b\i,f_1,b,\cdots,b\i,f_n,b)$. Then $w'\in\bold D^+$ via $(S_w)^b$, and we obtain 
$f^b=\Pi^+(w')=\Pi^+(f_1^b,\cdots,f_n^b)$. As $\bold X_k$ is closed with respect to conjugation in $\ca L^+$ 
we conclude that $f^b\in\bold X_{k+1}$, contrary to the maximality of $k$; and this contradiction 
completes the proof that $\ca N^+\norm\ca L^+$. That is, 5.1(a) holds. 

Next: let $\ca M\norm\ca L^+$ be any partial normal subgroup of $\ca L^+$ such that $\ca M\cap\ca L=\ca N$. 
Let $\phi$ be the composite of the inclusion map $\ca L\to\ca L^+$ followed by the quotient map 
$\r:\ca L^+\to\ca L^+/\ca M$. Then $Ker(\r)=\ca M$ (cf. I.5.3), and so $Ker(\phi)=\ca M\cap\ca L=\ca N$. 
Let $\r_S$ be the restriction of $\r$ (and hence also of $\phi$) to $S$. 

Set $\bar\D=\{P\phi\mid P\in\D\}$, and write $\bar P$ for $P\phi$. Let $\bar w\in\bold W(\ca L^+/\ca M)$, 
and set $\bar X=\bar S_{\bar w}$. That is (and in the usual way) $\bar X$ is the largest subgroup of 
$\bar S$ which is conjugated successively into $\bar S$ by the entries of $\bar w$. Write 
$\bar w=(\bar g_1,\cdots,\bar g_n)$, and for each $i$ let $g_i\in\ca L^+$ be a representative of $\bar g_i$ 
such that $g_i$ is $\up$-maximal with respect to $\ca M$. Set $w=(g_1,\cdots,g_n)$, and let  
$X$ be the preimage of $\bar X$ in $S$ via $\r_S$. Then $X\leq S_w$ by I.4.13(e). 
If $\bar X\in\bar\D$ then $X\in\D$ and $w\in\bold D(\ca L^+\mid_\D)$. As $\ca L^+\mid\D=\ca L$, we have 
$w\in\bold D(\ca L)$ in this case, and hence $\bar w\in\bold D((\ca L^+/\ca M)\mid_{\bar D})$. 
On the other hand, if $\bar w\in\bold D((\ca L^+/\ca M)\mid_{\bar D})$ then $\bar X\in\bar\D$. Thus: 
$$ 
Im(\phi)=(\ca L^+/\ca M)\mid_{\bar D}.  
$$ 
In particular, $Im(\phi)$ is a partial group (and a locality) in this way, and the mapping 
$$ 
\phi_0:\ca L\to Im(\phi) 
$$ 
given by $\phi$ is then a homomorphism of partial groups. One observes that the induced mapping 
$$ 
\phi_0^*:\bold W(\ca L)\to\bold W(Im(\phi)) 
$$ 
of free monoids carries $\bold D(\ca L)$ onto $\bold W(Im(\phi))$. That is, $\phi_0$ is a projection.  
Then, by I.5.6 there is an isomorphism $\bar\phi:\ca L/\ca N\to Im(\phi)$ induced by $\phi$. 

Let $\eta:\ca L/\ca N\to\ca L^+/\ca M$ be the composition of $\bar\phi$ followed by the inclusion map of 
$Im(\phi)$ into $\ca L^+/\ca M$. Our aim is to now show that $\ca L^+/\ca M$ is an expansion of $\ca L/\ca N$ 
via $\eta$. For this, what is required (cf. definition 4.1) is: (A) $\eta$ is injective (obvious, since 
$\eta$ is a composite of injective mappings); (B) $Im(\eta)$ is the restriction of $\ca L^+/\ca M$ to $\D\phi$ 
(which has been verified above);  and (C) $\D\phi$ is closed in the fusion system 
$\ca F'=\ca F_{\bar S}(\ca L^+/\ca M)$. In order to verify (C) we note first of all that 
$\ca F=\ca F_S(\ca L^+)$ by 4.3(a). Since $\r$ and $\phi$ are projections, (C) then follows from 1.10. 
Thus $\eta$ defines an expansion, as claimed. Moreover, $\eta$ defines an elementary expansion, since 1.10 
shows that the set $\{\bar P\mid P\in\D^+\}$ of objects of $\ca L^+/\ca M$ is the union of $\bar D$ with 
the set of $\ca F_{\bar S}(\ca L^+/\ca M)$-conjugates of $\bar R$. 

We have also to verify that Hypothesis 4.2 holds with respect to $\ca L/\ca M$ and $\bar R$, in order to 
apply Theorem 4.3 to the expansion $\ca L^+/\ca M$ of $Im(\phi)$. Set $T=S\cap\ca N$, and suppose first that 
$T\nleq R$. Then $TR\in\D$, so $\bar R\in\bar\D$, and then 4.2 holds trivially. On the other hand, suppose 
that $T\leq R$. Then every subgroup $\bar Q$ of $\bar S$ which properly contains an 
$\ca F_{\bar S}(\ca L^+/\ca M)$-conjugate of $\bar R$ is the image under $\phi$ of a subgroup $Q$ of $S$ 
which properly contains an $\ca F$-conjugate of $R$. That is, 4.2(1) obtains. Now let  
$\bar g\in N_{\ca L^+/\ca M}(\bar R)$, and let $g\in\ca L^+$ be a $\r$-preimage of $g$ such that $g$ is 
$\up$-maximal with respect to $\ca M$. Then $g\in N_{\ca L^+}(R)$ by I.4.13(e). Here 
$N_{\ca L^+}(R)=N_{\ca L}(R)$ by 4.3 applied to $\ca L^+$, so we have shown that $\r$ maps 
$N_{\ca L}(R)$ onto $N_{\ca L^+/\ca M}(\bar R)$. This verifies 4.2(2). 

The preceding analysis applies with $\ca N^+$ in the role of $\ca M$, where we now observe that 
$\ca N^+\leq\ca M$ since $\ca N\sub\ca M$. Let $\eta_0:\ca L/\ca N\to\ca L^+/\ca N^+$ be the homomorphism 
which defines the expansion of $\ca L/\ca N$ to $\ca L^+/\ca N^+$. By Theorem 4.3(b) there is then a  
unique homomorphism $\b:\ca L/\ca N^+\to\ca L/\ca M$ such that $\eta_0\circ\b=\eta$, and $\b$ is in fact 
an isomorphism. On the other hand, let $\g:\ca L/\ca N^+\to\ca L/\ca M$ be the mapping which sends a 
maximal coset $[g]_{\ca N^+}$ of $\ca N^+$ in $\ca L^+$ to the unique maximal coset $[g]_{\ca M}$ containing 
it. One easily verifies that $\g$ is a homomorphism and that $\eta_0\circ\g=\eta$, and thus $\g$ is an 
isomorphism. For $g\in\ca M$ we then have $[g]_{\ca N^+}=\ca M$, and this shows that $\ca M=\ca N^+$. 
That is, 5.1(c) holds, and we have obtained point (a) of the lemma. Point (b) was obtained in the course 
of proving (a).  
\qed 
\enddemo

The remainder of this section will be devoted to the proof of 5.1(b).

\definition {(5.4) Notation} Set $M=N_{\ca L}(R)$, and $K=M\cap\ca N$. As in the preceding section,  
for each $g\in\ca L^+$ set 
$$ 
\ca U_g=\{U\in R^{\ca F}\mid U\leq S_g\},  
$$ 
and for each $U\in R^{\ca F}$ set 
$$ 
\bY_U=\{y\in\ca L\mid R^y=U,\ N_S(U)^{y\i}\leq N_S(R)\}. 
$$ 
We continue to write $\bold X$ for $\ca N^{\ca L^+}$. 
\enddefinition

The notation (5.4) will remain fixed until the proof of Theorem 5.1 is complete. Note that $M$ is a 
subgroup of $\ca L$ by 4.2(2). Then $K$ is a normal subgroup of $M$ by I.1.8(c). 
\vskip .1in

\proclaim {Lemma 5.5} Suppose $T\leq R$. Then $\ca N^*=\bold X$, and $\bold N^*\cap\ca L=\ca N$. 
\endproclaim 

\demo {Proof} Let $f\in\ca N$, and suppose that $S_f$ contains an $\ca F$-conjugate of $R$. Then 
$T\leq S_f$ since, by I.4.1(a) $T$ is weakly closed in $\ca F$. Then I.4.1(b) yields the following result. 
\roster 

\item "{(1)}" Let $f\in\ca N$ such that $S_f$ contains an $\ca F$-conjugate $U$ of $R$. Then 
$P^f=P$ for each subgroup $P$ of $S_f$ containing $U$. In particular, $U^f=U$. 

\endroster 

Let $f'\in\bold X$, and let $f\in\ca N$ and $g\in\ca L^+$ such that $f'=f^g$. That is, assume that 
$v:=(g\i,f,g)\in\bold D^+$ and that $f'=\Pi^+(v)$. If $v\in\bold D$ then $f'=\Pi(v)\in\ca N$. On the other 
hand, suppose that $v\notin\bold D$. Then $S_v$ is an $\ca F$-conjugate of $R$ by 4.2(1). Set 
$U=(S_v)^{g\i}$. Then $U=U^f$ by (1). Now choose $a\in\bY_U$, and set $v'=(g\i,a\i,a,f,a\i,a,g)$. 
Then $v'\in\bold D^+$ via $U^g$, and $\Pi^+(v)=\Pi^+(v')$. Notice that (1) implies that 
$(a,f,a\i)\in\bold D$ via $P:=N_{S_f}(U)$, and that 
$$ 
T@>a>>U@>f>>U@>a\i>>T, 
$$ 
so that $f^{a\i}\in M$. Thus $f^{a\i}\in K$, and $\Pi^+(v)=\Pi^+(g\i a\i,f^{a\i},ag)$. This shows:   
\roster 

\item "{(2)}" $\bold X$ is the union of $\ca N$ with the set of all $\Pi^+(g\i,f,g)$ such that $f\in K$ 
and such that $(g\i,f,g)\in\bold D^+$. Moreover, for any such $v=(g\i,f,g)$, $\Pi^+(v)$ normalizes 
each $V\in R^{\ca F}$ such that $V\leq S_v$. 

\endroster 

Assume now that we have $v=(g\i,f,g)$ as in (2) (so that $f\in K$), and let $A$ be an $\ca F$-conjugate 
of $R$ contained in $S_v$. In order to analyze these things further, we shall need to be able to compute 
products in $\ca L^+$ in the manner described in 4.14 and 4.15. To that end, note first of all that since 
$A^g=T$ there exists a unique $h\in M$ and $y\in\bY_A$ such that $g\approx(\1,h,y)\in\Phi$, by 4.8. 
Set $\phi=(\1,h,y)$ and set $\psi=(\1,f,\1)$. Then $(\phi\i,\psi,\phi)$ is a $\G$-form of $(g\i,f,g)$, 
as defined in 4.13. We then compute via 4.14 that 
$$ 
f'=\Pi^+(g\i,f,g)=[y\i,h\i,\1][\1,f,\1][\1,h,y]=[y\i,f^h,y].  
$$ 
\roster 

\item "{(3)}" Let $f'\in\bold X$, such that $S_{f'}$ contains an $\ca F$-conjugate of $R$. Then $f'$ is 
an $\approx$-class $[y\i,k,y]$ with $k\in K$. 

\endroster 
If $f'\in\ca L$ then $(y\i,k,y)\in\bold D$ by 4.11, and so $f'\in\ca N$. Thus: 
\roster 

\item "{(4)}" $\bold X\cap\ca L=\ca N$. 

\endroster

Now let $w=(f_1',\cdots,f_n')\in\bold W(\bold X)\cap\bold D^+$, and set $B=S_w$. 
Suppose that $\Pi^+(w)\notin\bold X$. Then $\Pi^+(w)\notin\ca N$, so $\Pi^+(w)\notin\ca L$ by (4). 
Thus $w\notin\bold D$, so $B\in R^{\ca F}$, and then (2) shows that each $f_i'$ normalizes $B$. 
Fix $b\in\bold Y_B$. Then (3) implies that there exist elements $k_i\in K$ such that 
$f_i'=[b\i,k_i,b]$. One observes that the sequence of elements $(b\i,k_i,b)$ of $\Phi$ is a 
$\G$-form for $w$, and then 4.14 yields $\Pi^+(w)=[b\i,k,b]$ where $k=\Pi(k_1,\cdots,k_n)\in K$. 
This simply means that $\Pi^+(w)=k^b$, since $b\i=[b\i,\1,\1]$, $k=[\1,k,\1]$, and $b=[\1,\1,b]$. 
Thus $\Pi^+(w)\in\bold X$. Since $\bold X$ is closed under inversion, we have thus shown that 
$\bold X$ is a partial subgroup of $\ca L$. 

Finally, let $c\in\ca L^+$ be given so that $(c\i,f',c)\in\bold D^+$ (and where $f'\in\bold X$). Suppose  
that $\Pi^+(c\i,f',c)\notin\bold X$. Then $f'\notin\ca N$, so $f'\notin\ca L$ by (4), and then 
$f'=\Pi^+(g\i,f,g)$ for some $f\in\ca N$ and some $g\in\ca L^+$, and where $S_{f'}=S_{(g\i,f,g)}$. 
Set $u=(c\i,g\i,f,g,c)$. Then $u\in\bold D^+$ via $(S_{f'})^c$, and $\Pi^+(u)=\Pi^+(c\i g\i,f,gc)$. 
Thus $\Pi^+(c\i,f',c)\in\bold X$ after all, and $\bold X\norm\ca L^+$. Thus $\bold X=\ca N$, and the proof 
is complete. 
\qed 
\enddemo

Let $\bar L$ be the quotient locality $\ca L/\ca N$, let $\r:\ca L\to\bar{\ca L}$ be the 
quotient map, and let $\r^*$ be the induced homomorphism $\bold W(\ca L)\to\bold W(\bar{\ca L})$ of 
free monoids.  For any subset or element $X$ of $\ca L$, $\bar X$ shall denote the image of $X$ 
under $\r$. We extend this convention to subsets and elements of $\bold W(\bar{\ca L})$ in the obvious 
way. Set $\bar\D=\{\bar P\mid P\in\D\}$. 

\proclaim {Lemma 5.6} Assume that $T\nleq R$, and set $\bar F=\ca F_{\bar S}(\bar{\ca L})$. Then 
$\bar R^{\bar{\ca F}}\sub\bar\D$. 
\endproclaim 

\demo {Proof} Let $U\in R^{\ca F}$. As $T\nleq R$ we then have $T\nleq U$ by I.3.1(a). Then 
$U$ is a proper subgroup of $UT$, so $UT\in\D$ by 4.2(1). Then $\bar U=\bar{UT}\in\bar\D$.  
\qed 
\enddemo

\proclaim {Lemma 5.7} Assume that $T\nleq R$. There is then a homomorphism $\s:\ca L^+\to\bar{\ca L}$ 
such that the restriction of $\s$ to $\ca L$ is the quotient map $\r$. 
\endproclaim 

\demo {Proof} Set $\bar{\bold D}=\bold D(\bar{\ca L})$ and let $\bar\Pi:\bar{\bold D}\to\bar{\ca L}$ be 
the product in $\bar{\ca L}$. As $\bar R^{\bar{\ca F}}\sub\bar\D$ by 5.6, we have 
$\bar\Phi\sub\bar{\bold D}$, and so there is a mapping $\l:\Phi\to\bar{\ca L}$ given by 
$\phi\l=\bar\Pi(\phi\r^*)$, where $\r^*:\bold W(\ca L)\to\bold W(\bar{\ca L})$ is the homomorphism 
of free monoids induced by $\r$. That is: $\phi\l=\bar\Pi(\bar\phi)$. 

Let $\phi_1,\phi_2\in\Phi$ with $\phi_1\sim\phi_2$, and write 
$\phi_i=(x_i\i,h_i,y_i)$. Then $(x_2 x_1\i)h_1=h_2(y_2 y_1\i)$ (in $M$), and so 
$$ 
(\bar x_2\bar x_1\i)\bar h_1=\bar h_2(\bar y_2\bar y_1\i)\tag 1
$$ 
in $\bar M$. As $(x_2,x_1\i,h_1)$ and $(h_2,y_2,y_1\i)$ are in $\bold D^+$ via the appropriate conjugates 
of $R$, we have $(\bar x_2,\bar x_1\i,\bar h_1)$ and $(\bar h_2,\bar y_2,\bar y_1\i)$ in $\bar{\bold D}$, 
and then 
$$ 
\bar\Pi(\bar x_2,\bar x_1\i,\bar h_1)=\bar\Pi(\bar h_2,\bar y_2,\bar y_1\i) 
$$ 
by (1) and $\bar{\bold D}$-associativity. A standard cancellation argument (cf. I.2.4(a)) then yields 
$\bar\Pi(\bar\phi_1)=\bar\Pi(\bar\phi_2)$, and thus $\l$ is constant on $\sim$-classes. 

Now suppose that $\phi\in\Phi_g$ for some $g\in\ca L$. That is, suppose that $\phi\in\Phi\cap\bold D$ 
and that $g=\Pi(\phi)$. Then $\phi\l=\bar g$ since $\r$ is a homomorphism of partial groups, and 
so $\l$ is constant on $\approx$-classes. This shows that there is a (well-defined) mapping 
$\s:\ca L^+\to\bar{\ca L}$ given by $\r$ on $\ca L$ and by $\l$ on $\ca L_0^+$. It remains only to 
check that $\s$ is a homomorphism. 

Let $w=(f_1,\cdots,f_n)\in\bold D^+$. If $w\in\bold D$ then $w\s^*=w\r^*\in\bar{\bold D}$ and 
$\bar\Pi(w\s^*)=(\Pi(w))\s$. On the other hand, suppose that $w\notin\bold D$. Then $f_i=[\phi_i]$ 
for some $\phi_i\in\Phi$, and there is a $\G$-form $\g=(\phi_1,\cdots,\phi_n)$ of $w$. Thus the word 
$w_\g=\phi_1\circ\cdots\circ\phi_n$ has the property that $S_{w_\g}$ is an $\ca F$-conjugate of $R$, 
and hence $\bar{w_\g}\in\bar{\bold D}$. Then 
$$ 
\bar\Pi((w_\g)\s^*)=\bar\Pi(\bar{w_\g}). 
$$ 
Write $\phi_i=(x_i\i,h_i,y_i)$. Then 
$$ 
(\Pi^+(w_\g))\s=[x_1\i,\Pi(w_0),y_n]\s 
$$ 
where $w_0\in\bold D(M)$ is given by the formula in 4.14(b). One then obtains 
$\bar\Pi((w_\g)\s^*)=(\Pi^+(w_\g))\s$ via $\bar{\bold D}$-associativity, and so $\s$ is a homomorphism, 
as desired. 
\qed 
\enddemo

\proclaim {Lemma 5.8} $\ca N^*\cap\ca L=\ca N$. 
\endproclaim 

\demo {Proof} Let $\s:\ca L^+\to\bar{\ca L}$ be a homomorphism, as in 5.7, whose restriction to 
$\ca L$ is the quotient map $\r:\ca L\to\bar{\ca L}$. Then $Ker(\s)\cap\ca L=Ker(\r)$, where 
$Ker(\r)=\ca N$, by I.5.3. As $Ker(\s)\norm\ca L^+$, by I.1.14, the lemma follows. 
\qed 
\enddemo

Theorem 5.1 now follows from lemmas 5.3, 5.5, and 5.8.

\vskip .2in 
\noindent 
{\bf Section 6: Proper localities} 
\vskip .1in

Recall from I.3.6 that if $G$ is a group and $p$ is a prime then $Syl_p(G)$ is the (possibly empty) set 
of maximal $p$-subgroups $S$ of $G$ such that every $p$-subgroup of $G$ is conjugate in $G$ to a subgroup 
of $S$. For any group $G$ with $Syl_p(G)\neq\nset$, define $O_p(G)$ to be the subgroup $\bigcap Syl_p(G)$ of 
$G$. Thus $O_p(G)$ is the largest normal $p$-subgroup of $G$. If $(\ca L,\D,S)$ is a locality then there is 
another usage, whereby $O_p(\ca L)=\1^\star$, where $(\Omega,\star)$ is the stratification on 
$\ca F_S(\ca L)$ induced from $\ca L$. If it so happens that $\ca L$ is a group then these two definitions 
of $O_p(\ca L)$ agree by I.3.7, so there is no cause for confusion.

\definition {Definition 6.1} A group $G$ is of {\it characteristic $p$} if $Syl_p(G)\neq\nset$ and 
$C_G(O_p(G))\leq O_p(G)$. 
\enddefinition 

For a fusion system $\ca F$ on a group $S$ there is a largest subgroup $X$ of $S$ such that 
$\ca F=N_{\ca F}(X)$, by 1.9. In the case that $S$ is a $p$-group we write $O_p(\ca F)$ for $X$. 

\proclaim {Lemma 6.2} Let $(\ca L,\D,S)$ be a locality on $\ca F=\ca F_S(\ca L)$, and assume that 
$N_{\ca L}(P)$ is of characteristic $p$ for all $P\in\D$. Let $R\leq S$ be a subgroup of $S$. Then 
$R\norm\ca L$ if and only if $R\norm\ca F$. In particular, we have $O_p(\ca F)=O_p(\ca L)$. 
\endproclaim 

\demo {Proof} As $\ca F$ is generated by the conjugation maps $c_g:S_g\to S$ for $g\in\ca L$, we have 
$R\norm\ca F$ if $R\norm\ca L$. So assume $R\norm\ca F$ and, 
by way of contradiction, that $R$ is not normal in $\ca L$. Among all 
elements of $\ca L$ not in $N_{\ca L}(R)$, choose $g$ so that $dim(S_g)$ is as large as possible. Then 
$R\nleq S_g$, since $c_g:S_g\to S$ is an $\ca F$-homomorphism. In particular, $S_g\neq S$. As $S_g\in\Omega$, 
$S_g$ is then a proper subgroup of $N_S(S_g)$ by I.3.2. 

Set $P=S_g$, $P'=P^g$, and let $Q\in P^{\ca F}$ be fully normalized in $\ca F$. By  
by 3.2(a) there exists $x\in\ca L$ such that $P^x=Q$ and such that $N_S(P)\leq S_x$. As $P\in\Omega$ 
we have $dim(P)<dim(N_S(P))$ by I.3.2, and so $R\leq S_x$. Since $Q\in(P')^{\ca F}$ there exists also 
$y\in\ca L$ such that $(P')^y=Q$ and such that $N_S(P')R\leq S_y$. 

Note that $(x\i,g,y)\in\bold D$ via $Q$, and that $f:=\Pi(x\i,g,y)\in N_{\ca L}(Q)$. Note also that 
$(x,x\i,g,y,y\i)\in\bold D$ via $P$, and that 
$$ 
\Pi(x,f,y\i)=\Pi(x,x\i,g,y,y\i)=\Pi(g) 
$$ 
by $\bold D$-associativity (I.1.4). If $R\leq S_f$ then $R\leq S_{(x,f,y\i)}$, and then $R\leq S_g$. Thus 
$R\nleq S_f$, and we may therefore replace $g$ with $f$, and $P$ with $Q$. That is, we may assume that  
$P$ is fully normalized in $\ca F$ and that $g$ is in $N_{\ca L}(P)$. 

Set $H=N_{\ca L}(P)$, set $D=N_R(P)$, and set $\ca E=\ca F_{N_S(P)}(H)$. Then $D\norm\ca E$, and so the 
conjugation map $c_g:P\to P$ extends to an $\ca E$-automorphism of $PD$. Thus, there exists $h\in N_H(PD)$ 
such that $c_h:PD\to PD$ restricts to $c_g$ on $P$. Then $gh\i\in C_H(P)$. Since $H$ is is of 
characteristic $p$ by hypothesis, we then have $gh\i\in P$. Then $gh\i\in N_H(D)$. 
Since $c_h$ extends to an $\ca F$-homomorphism which fixes $R$, we have $D^h=D$, and so also 
$D^g=D$. Thus $D\leq P$, and then $R\leq P$ by I.3.2. This contradicts the choice of $g$, and 
completes the proof. 
\qed 
\enddemo 

By a {\it $p'$-group} we mean a torsion group $G$ all of whose elements have order relatively prime to $p$. 
In I.3.2 it was seen that the finite-dimensionality condition on a locality $(\ca L,\D,S)$ has an important 
consequence for the structure of the $p$-group $S$. Namely, it was shown that if $X\in\Omega$ and $Y\leq S$ 
is a subgroup of $S$ with $X<Y$ (proper inclusion), then $X<N_Y(X)$. In order to begin to analyze the 
interplay between $p$-subgroups and $p'$-subgroups of $\ca L$, we shall need to generalize Thompson's 
$A\times B$-Lemma (concerning finite groups). First, we recall the definition I.6.2:  
A $p$-group $S$ has the {\it normalizer-increasing property} if, for every 
pair $X$ and $Y$ of subgroups of $S$ with $X<Y$, we have $X<N_Y(X)$. (The remark following I.6.2 
provides an example of a countable, locally finite $p$-group which does not have the normalizer-increasing 
property.

\proclaim {Lemma 6.3} Let $G$ be a torsion group, and let $A$ and $B$ be subgroups of $G$ such 
that the order of every element of $A$ is relatively prime to the order of every element of $B$. 
Suppose that $[A,B]\leq C_B(A)$. Then $[A,B]=\1$. 
\endproclaim 

\demo {Proof} Let $a\in A$ and $b\in B$. Then $a\i a^b=[a,b]$ commutes with $a$, by hypothesis, and 
so $a^b$ commutes with $\<a\>$. As $|a|=|a^b|$ is relatively prime to the order of any element of $B$, the 
same is then true of $a\i a^b$. As $[a,b]\in B$ by hypothesis, we conclude that $[a,b]=1$, and thus 
$[A,B]=\1$. 
\qed 
\enddemo

\proclaim {Lemma 6.4 (Thompson's ``$A\times B$" Lemma)} Let $G$ be a torsion group, and assume:  
\roster 

\item "{(i)}" $G=AP$ where $P$ is a normal $p$-subgroup of $G$ having the normalizer-increasing property, 
and where $A$ is a $p'$-subgroup of $G$. 

\item "{(ii)}" There exists a subgroup $B$ of $C_P(A)$ such that $[C_P(B),A]=1$. 

\endroster 
Then $G$ is isomorphic to the direct product $A\times P$. 
\endproclaim 

\demo {Proof} As $P\cap A=\1$ it suffices to show that $[P,A]=1$. Suppose false, so that 
$C_P(A)$ is a proper subgroup of $P$. Set $Q=N_P(C_P(A))$. $B\leq C_P(A)\norm Q$. Thus $[Q,B]\leq C_P(A)$; 
and thus $[[Q,B],A]=1$. Also $[[B,A],Q]=1$ since $B\leq C_P(A)$, and 
then the Three Subgroups Lemma yields $[[A,Q],B]=1$. That is, we have $[A,Q]\leq C_Q(B)$, and hence 
$[A,Q]\leq C_Q(A)$ by (ii). Then $Q\leq C_P(A)$ by 6.3. The normalizer-increasing hypothesis for $P$ 
then yields $C_P(A)=P$, as required.  
\qed 
\enddemo

By a {\it $p$-local subgroup} of a group $G$ we mean a subgroup of the form $N_G(P)$, where 
$P$ is a $p$-subgroup of $G$.

\proclaim {Lemma 6.5} Let $G$ be a countable, locally finite group of characteristic $p$.  
\roster 

\item "{(a)}" Every normal subgroup of $G$ is of characteristic $p$. 

\item "{(b)}" If a Sylow $p$-subgroup of $G$ has the normalizer-increasing property  
then every $p$-local subgroup of $G$ is of characteristic $p$. 

\item "{(c)}" Let $V$ be a normal $p$-subgroup of $G$, and let $X$ be the set of elements $x\in C_G(V)$ 
such that $[O_p(G),x]\leq V$. Then $X$ is a normal $p$-subgroup of $G$. 

\endroster 
\endproclaim 

\demo {Proof} Let $K\norm G$ be a normal subgroup of $G$ and set $R=O_p(K)$. Then 
$[O_p(G),K]\leq R$, and so $[O_p(G),C_K(R),C_K(R)]=\1$. Then the group $O^p(C_K(R))$ generated by the 
$p'$-elements of $C_K(R)$ is contained in $C_G(O_p(G))$ by 6.3, and thus $O^p(C_K(R))\leq O_p(G)$. 
Then $C_K(R)$ is a normal $p$-subgroup of $K$, and so $C_K(R)\leq R$. Thus (a) holds. 

Next, let $U$ be a $p$-subgroup of $G$, and set $H=N_G(U)$, $P=O_p(H)$, and $Q=O_p(G)$. Then 
$N_Q(U)\leq P$. Let $A$ be a $p'$-subgroup of $C_H(P)$. Then $[N_Q(U),A]=\1$. Assuming now that 
$Q$ has the normalizer-increasing property, it follows from 6.4 that $[Q,A]=\1$. 
Then $A\leq Q$, and thus $A=\1$, proving (b). 

For the proof of (c), notice that $C_G(V)\norm G$, and that $X$ is the intersection of $C_G(V)$ with 
the preimage in $G$ of the normal subgroup $C_{G/V}(O_p(G)/V)$ of $G/V$. Thus $X\norm G$. 
Each $p'$-subgroup of $X$ centralizes $O_p(G)$ by 6.3, so $X$ is a $p$-group. 
\qed 
\enddemo 

Recall that, for any fusion system $\ca F$ on $S$, $\ca F^c$ is the set of $\ca F$-centric 
subgroups of $S$; i.e. subgroups $P$ of $S$ such that $C_S(Q)\leq Q$ for all $Q\in P^{\ca F}$.  

\definition {Definition 6.6} Let $(\ca F,\Omega,\star)$ be a stratified fusion system on the 
countable, locally finite $p$-group $S$. A subgroup $P\leq S$ is: 
\roster 

\item "{$\cdot$}" {\it $\ca F$-radical} if there exists an $\ca F$-conjugate $Q$ of $P$ such 
$Q$ is fully normalized in $(\ca F,\Omega,\star)$ and such that $Q=O_p(N_{\ca F}(P))$, 

\item "{$\cdot$}" {\it $\ca F$-quasicentric} if there exists an $\ca F$-conjugate $Q$ of $P$ such that $Q$ 
is fully normalized in $(\ca F,\Omega,\star)$ and such that $N_{\ca F}(Q)$ is the trivial fusion system on 
$N_S(Q)$ (cf. 1.4(2)), and 

\item "{$\cdot$}"  {\it $\ca F$-subcentric} if there exists an $\ca F$-conjugate $Q$ of $P$ such that $Q$ is 
fully normalized in $(\ca F,\Omega,\star)$ and such that $O_p(N_{\ca F}(Q))\in\ca F^c$. 

\endroster 
Denote by $\ca F^{cr}$ the set of all subgroups $P\leq S$ such that $P$ is both $\ca F$-centric and 
$\ca F$-radical, by $\ca F^q$ the set of all subgroups $P\leq S$ such that $P$ is $\ca F$-quasicentric, 
and by $\ca F^s$ the set of all subgroups $P\leq S$ such that $P$ is $\ca F$-subcentric. 
\enddefinition 

\definition {Remark} The definition given above of $\ca F$-radical subgroup differs from the 
standard one (in the case where $S$ is finite). Namely, in the standard definition 
a subgroup $P$ of $S$ is $\ca F$-radical if $Inn(P)=O_p(Aut_{\ca F}(P))$, but we 
find the definition given above to be more useful. We mention also that the notion of 
$\ca F$-subcentric subgroup was first explored by Ellen Henke [He2]. 
\enddefinition 

\definition {Definition 6.7} Let $(\ca L,\D,S)$ be a locality. Then $\ca L$ is {\it proper} if: 
\roster 

\item "{(PL1)}" $\ca F^{cr}\sub\D$. 

\item "{(PL2)}" Each of the groups $N_{\ca L}(P)$, for $P\in\D$, is of characteristic $p$. 

\item "{(PL3)}" $S$ has the normalizer-increasing property. 

\endroster 
\enddefinition

\proclaim {Lemma 6.8} Let $(\ca L,\D,S)$ be a locality on $\ca F$. Assume 
\roster 

\item "{(1)}" $S$ has the normalizer-increasing property, and 

\item "{(2)}" $\ca F^{cr}\sub\D\sub\ca F^c$. 

\endroster 
Then there exists a partial normal subgroup $\Theta\norm\ca L$ such that $S\cap\Theta=\1$ and such that,  
upon identifying $S$ with its image in $\ca L/\Theta$, the quotient $\ca L/\Theta$ is a proper locality 
on $\ca F$. 
\endproclaim 

\demo {Proof} As $\D\sub\ca F^c$ it follows from  
I.6.1 that for each $P\in\D$ there is a largest normal 
$p'$-subgroup $\Theta(P)$ contained in $N_{\ca L}(P)$, and that then $N_{\ca L}(P)/\Theta(P)$ is of 
characteristic $p$. Define $\Theta$ to be the union of the groups $\Theta(P)$ for $P\in\D$. With (1) 
the hypotheses of I.6.3 are fulfilled, and we obtain  
$\Theta\norm\ca L$, $\ca F=\ca F_S(\ca L/\Theta)$, and $\ca L/\Theta$ is a proper locality 
on $\ca F$.  
\qed 
\enddemo

\proclaim {Lemma 6.9} Let $(\ca L,\D,S)$ be a proper locality on $\ca F$. Then $\D\sub\ca F^s$, 
and for any $P\in\D$: 
\roster 

\item "{(a)}" $P\in\ca F^{cr}$ if and only if $P=O_p(N_{\ca L}(P))$. 

\item "{(b)}" $P$ is centric in $\ca F$ if and only if $C_{\ca L}(P)=Z(P)$. 

\item "{(c)}" $P$ is quasicentric in $\ca F$ if and only if $C_{\ca L}(P)\leq O_p(N_{\ca L}(P))$. 

\endroster 
\endproclaim 

\demo {Proof} Let $P\in\D$, set $M=N_{\ca L}(P)$, and let $Q\in P^{\ca F}$ with $Q$ fully normalized in 
$\ca F$. Then $Q=P^g$ for some $g\in\ca L$ by I.3.5, and then the conjugation map $c_g:M\to M^g$ is an 
isomorphism by I.2.3(b). As $\ca F^{cr}$, $\ca F^c$, $\ca F^q$, and $\ca F^s$ are $\ca F$-invariant by 1.9, 
it suffices to establish that $P\in\ca F^s$ and that (a) through (c) hold, under the assumption (which we 
now make) that $P$ itself is fully normalized in $\ca F$. 

Regard $M$ as a locality whose set of objects is the set of all overgroups of 
$O_p(M)$ in $N_S(P)$. Since the group $M$ is of characteristic $p$ by (PL2), the locality $M$ is proper as a 
consequence of 6.6(b). As $N_{\ca F}(P)=\ca F_{N_S(P)}(M)$ by 3.2(b), 6.2 implies that 
$O_p(N_{\ca F}(P))=O_p(M)$. Here $C_S(O_p(M))\leq O_p(M)$ (as $M$ is of characteristic $p$), so 
3.4(a) yields $O_p(M)\in\ca F^c$. Thus $P\in\ca F^s$, and we have shown moreover that 
if $P\in\ca F^{cr}$ then $P=O_p(M)$. On the other hand, if $P=O_p(M)$ then 
$P=O_p(N_{\ca F}(P))$, and so (a) holds. 

Set $K=C_{\ca L}(P)$. Then $C_{\ca F}(P)=\ca F_{C_S(P)}(K)$ by 3.2(b). If $P\in\ca F^c$ then 
$C_S(P)=Z(P)\leq Z(K)$, and then $O^p(K)$ is a normal $p'$-subgroup of $K$ by I.6.1. Since $K$ is of 
characteristic $p$ by 6.6(a), $K$ has no non-identity normal $p'$-subgroups, so we conclude that $O^p(K)=1$.  
Thus $K=Z(P)$. Conversely, if $K=Z(P)$ then $C_S(P)\leq P$. As $P$ is fully normalized in $\ca F$, 3.4(a) 
then yields $P\in\ca F^c$. This establishes (b). 

Suppose next that $P\in\ca F^q$, so that $C_{\ca F}(P)$ is the trivial fusion system on $C_S(P)$. Then 
$N_K(U)/C_K(U)$ is a $p$-group for every subgroup $U$ of $C_S(P)$. Take $U=O_p(K)$. Thus  
$K/Z(U)$ is a $p$-group, so $K$ is a $p$-group, and then $K=C_S(P)$ is a normal 
$p$-subgroup of $M$. Conversely, if $K\leq O_p(M)$ then $C_{\ca F}(P)$ is trivial, so we have (c). 
\qed 
\enddemo

\vskip .2in 
\noindent 
{\bf Section 7: Theorem A} 
\vskip .1in 

This section concerns expansions of proper localities and of homomorphic images of proper 
localities. Throughout, $(\ca L,\D,S)$ will be a proper locality on $\ca F$, $(\Omega,\star)$ will be 
the stratification induced from $\ca L$, and the assertion that a subgroup $X$ of $S$ is ``fully normalized 
in $\ca F$" will mean that $X$ is fully normalized with respect to $(\Omega,\star)$. We shall see 
somewhat later (in 8.3) that $\ca F$ is inductive, so by 2.9 the notion of being fully normalized in $\ca F$ 
is independent of the given stratification.

\proclaim {Lemma 7.1} Let $(\ca L,\D,S)$ be a proper locality on $\ca F$ and let $R\in\ca F^s$. Assume:  
\roster 

\item "{(1)}" $R$ is fully normalized in $\ca F$, and 

\item "{(2)}" every strict overgroup of $R$ in $S$ is in $\D$. 

\endroster 
Set $M=N_{\ca L}(R)$. Then: 
\roster 

\item "{(a)}" $M$ is a subgroup of $\ca L$, of characteristic $p$. 

\item "{(b)}" $N_S(R)\in Syl_p(M)$, and $N_{\ca F}(R)=\ca F_{N_S(R)}(M)$. 

\endroster 
In particular, if $R\in\Omega_S(\ca L)$ then $\ca L$ and $R$ satisfy the conditions of Hypothesis 4.2. 
\endproclaim 

\demo {Proof} We begin with the observation that the normalizer in $\ca L$ of any subgroup of $S$ 
is a partial subgroup of $\ca L$. Set $\D_R=\{P\in\D\mid R\norm P\}$. We will show that 
$(M,\D_R,N_S(R))$ is a proper locality on $N_{\ca F}(R)$. 

Let $Q\in\D$ with $R\leq Q$. Then either $R=Q$ or $R<N_Q(R)$, and in either case (2) implies that 
$N_Q(R)\in\D$. Thus: 
\roster 

\item "{(3)}" $\D_R=\{N_Q(R)\mid R\leq Q\in\D\}$. 

\endroster 
Set $\ca F_R=\ca F_{N_S(R)}(M)$. Then $\D_R$ is $\ca F_R$-closed since 
$\D_R$ is closed with respect to overgroups in $N_S(R)$, and since $\D_R$ is closed with 
respect to conjugation maps $c_g:P\to P^g$ with $g\in M$ and with $P\leq S_g$. 
Let $w\in\bold W(M)\cap\bold D(\ca L)$, and set $Q=S_w$. Then (3) yields $w\in\bold D(\ca L)$ via $N_Q(R)$, 
and we have thus verified that $(M,\D_R)$ is an objective partial group (I.2.1), and  
that $(M,\D_R,N_S(R))$ is a pre-locality. By I.2.15, $\Omega_{N_S(R)}(M)$ is finite-dimensional. 

Set $V=N_S(R)$. An argument that should by now be familiar will show that $V$ is a 
maximal $p$-subgroup of $M$. Namely, let $V\leq X$ where $X$ is a $p$-subgroup of $M$. By I.3.8 
there exists $g\in\ca L$ with $X^g\leq S$. Then $dim(V^g)=dim(X^g)$ since $R$ is fully normalized,  
and so $(V^g)^*=(X^g)^*$. Since $V=N_{V^*}(R)$ we get $V^g=N_{(V^g)^*}(R^g)$, so $V^g\geq X^g$, and 
$V=X$. Thus $N_S(R)$ is a maximal $p$-subgroup of $M$, and $(M,\D_R,V)$ is a locality on 
$\ca F_R$. 

For any $P\in\D_R$ we have 
$$ 
N_{\ca L_R}(P)=N_{N_{\ca L}(P)}(R), 
$$ 
and 6.5(b) then implies that $N_{\ca L_R}(P)$ is of characteristic $p$. Thus $\ca L_R$ satisfies the 
condition (PL2) in the definition 6.7 of proper locality. 

If $R\in\D$ then $M$ is a subgroup of $\ca L$ of characteristic $p$ (since $\ca L$ is proper), and then  
$N_{\ca F}(R)=\ca F_R$ by 3.2. Thus (a) and (b) hold in this case, and we may assume that $R\notin\D$. 
Set $Q=O_p(N_{\ca F}(R))$ and suppose next that $R=Q$. Then $R\in\ca F^c$ by (1), and so $R\in\ca F^{cr}$ 
(cf. definition 6.6). Since $\ca F^{cr}\sub\D$, and since we have eliminated the case where $R\in\D$, we 
conclude that $R\neq Q$. As $\ca F_R$ is a subsystem of $N_{\ca F}(R)$ we have $Q\leq O_p(\ca F_R)$, so 
$R<O_p(\ca F_R)$. Thus $O_p(\ca F_R)\in\D$, and so $(\ca F_R)^{cr}\sub\D_R$.  
That is: $\ca L_R$ satisfies the condition (PL1) in definition 6.7, 
Condition (PL3) (the normalizer-increasing property) is inherited by $N_S(R)$ from $S$, and thus 
$(\ca L_R,\D_R,N_S(R))$ is a proper locality. Then $O_p(\ca F_R)=O_p(M)$ by 6.2, and then $M$ is a group 
since $O_p(\ca F_R)\in\D$. As $Q\leq O_p(\ca F_R)$, 
(1) implies that $Q$ is centric in $\ca F$, hence also in $\ca F_R$, and so $M$ is of characteristic $p$. 
Thus (a) holds. 
 
As $N_S(R)$ is a maximal $p$-subgroup of $M$, I.3.11 yields $N_S(R)\in Syl_p(M)$. Since $Q\in\D$, each 
$N_{\ca F}(R)$-isomorphism is the restriction of a conjugation map $c_g:X\to Y$ with $Q\leq X\cap Y$, 
and for some $g\in L$ such that $R^g=R$. Thus $N_{\ca F}(R)$-isomorphisms are $\ca F_R$-isomorphisms, 
and this completes the proof of (b). 
\qed 
\enddemo

Our aim now is to extend the results of sections 4 and 5 beyond elementary expansions. The arguments 
amount to an exercise with Zorn's Lemma, and will be based on the restricted notion of 
expansion in definition 4.1 (where the homomorphism $\eta$ is assumed to be an inclusion map).

\proclaim {Theorem 7.2} Let $(\ca L,\D,S)$ be a proper locality on $\ca F$ and let $\D^+$ be an 
$\ca F$-closed collection of subgroups of $S$ such that $\D\sub\D^+\sub\ca F^s$. Then the following hold. 
\roster 

\item "{(a)}" There exists an expansion $(\ca L^+,\D^+,S)$ of $\ca L$ such that $\ca L^+$ is a 
proper locality on $\ca F$. Moreover, for any such expansion we then have 
$\ca F_S(\ca L^+)=\ca F$, and $\ca L^+$ is generated by $\ca L$ as a partial group. 

\item "{(b)}" For any expansion $(\w{\ca L},\D^+,S)$ having the same set $\D^+$ of objects as $\ca L^+$, and 
such that $\w{\ca L}$ is a proper locality, there is a unique homomorphism $\b:\ca L^+\to\w{\ca L}$ such that 
$\b$ restricts to the identity map on $\ca L$; and $\b$ is then an isomorphism. 

\endroster 
\endproclaim 

\demo {Proof} Let $\w{\Bbb E}$ be the class of all expansions $(\ca L',\D',S)$ of $\ca L$ such that 
$\ca L$ is a proper locality, with $\ca F_S(\ca L')=\ca F$, and such that $\ca L^+$ is generated 
by $\ca L$ as a partial group. For $\ca L_1,\ca L_2\in\w{\Bbb E}$, we say that a homomorphism 
$\phi:\ca L_1\to\ca L_2$ is {\it rigid} if $\phi$ restricts to the identity map on $\ca L$. 
Let $\approx$ be the equivalence relation on $\w{\Bbb E}$ given by 
$\ca L_1\approx\ca L_2$ if there exists a rigid isomorphism $\ca L_1\to\ca L_2$. 

For any locality $(\ca L,\D,S)$, each of the groups $N_{\ca L}(P)$ for $P\in\D$ is countable. Then 
$S$ is countable, and so the cardinality of $\D$ is at most $\aleph_1$. For any $P,Q\in\D$ the set 
$N_{\ca L}(P,Q)$ is countable, and since $\ca L$ is the union of the sets $N_{\ca L}(P,Q)$, the 
cardinality of $\ca L$, and then also of $\bold W(\ca L)$ and $\bold D(\ca L)$, is at most $\aleph_1$. 
There is then a {\it set} of isomorphism classes of localities, and a set $\Bbb E$ of $\approx$-classes 
of expansions $\ca L'\in\w{\Bbb E}$ of $\ca L$. Further, we may view $\Bbb E$ as a poset, whereby  
$X_1\preceq X_2$ if there exist representatives $\ca L_1\in X_1$ and $\ca L_2\in X_2$ such that 
there exists a rigid injective homomorphism $\phi\ca L_1\to\ca L_2$. Since $\ca L_1$ and $\ca L_2$ are 
generated by $\ca L$, $\phi$ is then the unique rigid homomorphism $\ca L_1\to\ca L_2$. 

Let $\{X_i\}_{i\in I}$ be a totally ordered subset of $\Bbb E$, and choose $\ca L_i\in X_i$ for all $i$. 
Then $\{\ca L_i\}_{i\in I}$ comes equipped with rigid homomorphisms $\phi_{i,j}:\ca L_i\to\ca L_j$, and 
the uniqueness of $\phi_{i,j}$ implies that $\{\ca L_i\}_{i\in I}$ is a directed system in this way. 
The category of partial groups possesses limits (Appendix B in Part I), so the direct limit of 
this directed system is 
a partial group $\ca L'$. For convenience, we may then identify each $\ca L_i$ with its image in $\ca L'$,  
and take the homomorphisms $\phi_{i,j}$ to be inclusion maps. 

Let $\D'$ be the union of the sets $\D_i$. Since each $(\ca L_i,\D_i,S)$ is a locality on $\ca F$, $\D'$ is 
$\ca F$-closed. For any $w\in\bold D(\ca L')$ we have $w\in\bold D(\ca L_i)$ for some $i$, and then 
$(\ca L',\D')$ is an objective partial group, and $(\ca L',\D',S)$ is a pre-locality. Since the 
stratification on $\ca F$ induced from $\ca L_i$ is the same as that which is induced from $\ca L$, 
$\Omega_S(\ca L')$ is finite-dimensional. Evidently $S$ is a maximal $p$-subgroup of $\ca L'$, so 
$(\ca L',\D',S)$ is a locality. For any $P\in\D'$ we have $P\in\D_i$ for some $i$, and then 
$N_{\ca L'}(P)=N_{\ca L_i}(P)$ is of characteristic $p$. Thus $\ca L'$ is a proper locality, and  
$\ca L'$ is an expansion of $\ca L$. The $\approx$-class of 
$\ca L'$ is thus an upper bound for $\ca C$ in $\Bbb E$, and so $\Bbb E$ contains a maximal element 
$Y$ by Zorn's Lemma. Rather than introduce further notation, let 
us now write $(\ca L',\D',S)$ for a representative of $Y$. 

Set $\w D=\{P\leq S\mid P^\star\in\D\}$. Then $(\ca L,\w D,S)$ is an expansion of $\ca L$ by I.3.15, and 
there is then no loss of generality in assuming that $\D=\w\D\cap\D^+$. Suppose now that $\D'\neq\D^+$, 
choose $R\in\D^+-\D'$ so that $dim(R)$ is as large as possible, and set $\D^*=\D'\cup R^{\ca F}$. 
Then $\D^*$ is $\ca F$-closed. Since $\D^+\sub\ca F^s$, we 
may assume that $R$ has been chosen also so that $R$ is fully normalized in $\ca F$ and so that 
$O_p(N_{\ca F}(R))\in\ca F^c$. Then 7.1 implies that the normalizer $M=N_{\ca L}(R)$ is a subgroup of 
$\ca L$ of characteristic $p$ and that $N_S(R)\in Syl_p(M)$. Theorem 4.3 then yields the existence of 
a proper locality $(\ca L^*,\D^*,S)$ on $\ca F$ whose restriction to $\D'$ is $\ca L'$, 
with $\ca F_S(\ca L^*)=\ca F$, and such that $\ca L'$ 
generates $\ca L^*$ as a partial group. This contradicts the maximality of $Y$, and completes the proof. 
\qed 
\enddemo

\proclaim {Theorem 7.3} Let $(\ca L,\D,S)$ be a proper locality on $\ca F$, and let $(\ca L^+,\ca D^+,S)$ 
be an expansion of $\ca L$ to an $\ca F$-closed subset $\D^+$ of $F^s$. Let $\frak N$ be the set 
of partial normal subgroups of $\ca L$, and let $\frak N^+$ be the set of partial normal subgroups of 
$\ca L^+$. For each $\ca N\in\frak N$ there is then a unique $\ca N^+\in\frak N^+$ such that 
$\ca N=\ca N^+\cap\ca L$. In particular, we have $S\cap\ca N=S\cap\ca N^+$, and the mapping 
$\ca N\maps\ca N^+$ is an inclusion-preserving bijection $\frak N\to\frak N^+$. 
\endproclaim 

\demo {Proof}  Fix $\ca N\norm\ca L$ and let $\w{\Bbb E}_{\ca N}$ to be the class of expansions 
$(\ca L',\D',S)$ of $\ca L$ such that $\D'\sub\D^+$ and such that there exists a unique partial normal 
subgroup $\ca N'\norm\ca L'$ with $\ca L\cap\ca N'=\ca N$. 
As in the proof of 7.2 there is then a set $\Bbb E_{\ca N}$ of rigid isomorphism-classes of members of 
$\w{\Bbb E}_{\ca N}$, $\Bbb E_{\ca N}$ is a poset, and any set of representatives for $\Bbb E_{\ca N}$ forms 
a directed set. As in the proof of 7.1, we may assume that $P\in\D$ whenever $P\leq S$ and $P^\star\in\D$. 

Let $\ca C=\{Y_i\}_{i\in I}$ be a totally ordered subset of $\Bbb E_{\ca N}$, and choose representatives 
$(\ca L_i,\D_i,S)\in Y_i$. As in the proof of 7.2, these representatives form a directed system 
whose direct limit is then a proper locality $(\ca L',\D',S)$ which is an expansion of $\ca L$. Also 
as in 7.2, we may take $\ca L_i$ to be the restriction of $\ca L'$ to $\D_i$, and take the connecting 
homomorphisms for the directed system to be inclusion maps. 

Let $\ca N_i$ be the unique partial normal subgroup of $\ca L_i$ whose intersection with $\ca L$ is $\ca N$. 
and then let $\ca N'$ be the union 
of the $\ca N_i$ in $\ca L'$. Plainly, $\ca N'$ is a partial subgroup of $\ca L'$, and is closed 
under $\ca L'$-conjugation (when conjugation is defined). Thus $\ca N'$ is a partial normal subgroup of 
$\ca L'$. One observes that 
$$ 
\ca N'\cap\ca L=\bigcup(\ca N'\cap\ca L_i)\cap\ca L=\bigcup(\ca N_i)\cap\ca L=\ca N.  
$$ 
In order to conclude that the $\approx$-class of $\ca L'$ is an upper bound for $\ca C$ in $\Bbb E$ it 
remains to show that if $\ca M$ is any partial normal subgroup of $\ca L'$ which meets $\ca L$ in $\ca N$ 
then $\ca M=\ca N'$. Indeed: 
$$ 
\ca M=\bigcup(\ca M\cap\ca L_i)=\bigcup(\ca N'\cap\ca L_i)=\ca N'.  
$$ 
Thus $\ca C$ has an upper bound, and Zorn's Lemma implies that $\Bbb E$ contains a maximal element $Y$.  

Continuing as in the proof of 7.2, let $(\ca L',\D',S)$ now denote a representative of $Y$. If 
$\D'\neq\D^+$ then there is an elementary expansion $(\ca L^*,\D^*,S)$ of $\ca L'$ which is not in 
$\w{\Bbb E}_{\ca N}$, and then Theorem 5.1 yields a contradiction. 
\qed 
\enddemo 

The preceding two results are the Theorems A1 and A2, as stated in the introduction, and we think of them 
as comprising the two parts of a single Theorem A. 

\vskip .1in 

Recall from I.6.8 that the product of any non-empty set of partial normal subgroups of a locality is 
again a partial normal subgroup. That is: if $\frak N\neq\nset$ is a set of partial normal subgroups of 
$\ca L$, then the partial subgroup $\ca M=\<\bigcup\frak N\>$ of $\ca L$ is normal in $\ca L$, and for 
each total ordering $\preceq$ of $\frak N$ and each element $x\in\ca M$, there exists a finite sequence 
$(\ca N_1,\cdots,\ca N_k)$ of elements of $\frak N$, and elements $x_i\in\ca N_i$ $(1\leq i\leq k)$, 
such that $\ca N_1\preceq\cdots\preceq\ca N_k$, $(x_1,\cdots,x_k)\in\bold D$, and $x=Pi(x_1,\cdots,x_k)$.

\proclaim {Corollary 7.4} Let $\frak N$ be a non-empty set of partial normal subgroups of the proper 
locality $(\ca L,\D,S)$, and let $(\ca L^+,\D^+,S)$ be a proper expansion of $\ca L$. Let $\ca M$ 
be the partial subgroup of $\ca L$ generated by $\bigcup\frak N$, and for any $\ca N\norm\ca L$ 
let $\ca N^+$ be the unique partial normal subgroup of $\ca L^+$ such that $\ca N^+\cap\ca L=\ca N$. Then 
$\ca M^+$ is generated by $\bigcup\{\ca N^+\}_{\ca N\in\frak N}$, and if $\bigcap\frak N\leq S$ then 
$(\bigcap\frak N)^+=\bigcap\frak N$. 
\endproclaim 
 
\demo {Proof} We provide the details only for the case that $\frak N$ consists of a pair 
$\{\ca N_1,\ca N_2\}$ of partial normal subgroups of $\ca L$ (this being the only 
case that will be required for later 
applications), and leave the general case to the reader who is interested in carrying out the exercise.  
For any pair $X_1$ and $X_2$ of non-empty subsets of $\ca L$ it follows from I.1.9 that 
$\<X_1\>\<X_2\>$ is a subset of $\<X_1X_2\>$. Then 
$\<\ca N_1^{\ca L^+}>\<\ca N_2^{\ca L^+}\>\sub\<\ca N_1^{\ca L^+}\ca N_2^{\ca L^+}\>$, and so 
$\ca N_1^+\ca N_2^+\leq(\ca N_1\ca N_2)^+$. The reverse inclusion follows from the observation that 
$\ca N_1\ca N_2\leq \ca L\cap \ca N_1^+\ca N_2^+$. 

Set $R=\ca N_1\cap\ca N_2$ and suppose $R\leq S$. We have $R\norm\ca L$, so 6.2 yields $R\norm\ca L^+$. 
Observing that 
$$ 
(\ca N_1^+\cap\ca N_2^+)\cap\ca L=(\ca N_1^+\cap\ca L)\cap(\ca N_2^+\cap\ca L)=\ca N_1\cap\ca N_2 
$$ 
and that $R\cap\ca L=\ca N_1\cap\ca N_2$, we conclude from 7.3 that $R=\ca N_1^+\cap\ca N_2^+$.    
\qed 
\enddemo

We wish also to have a version of Theorem A for localities which are homomorphic images of proper 
localities. 

\definition {Hypothesis 7.5} $(\ca L,\D,S)$ is a locality on $\ca F$, and 
$\ca N\norm\ca L$ is a partial normal subgroup of $\ca L$. Set $T=S\cap\ca N$. Let 
$(\bar{\ca L},\bar\D,\bar S)$ be the quotient locality $\ca L/\ca N$, and let $\r:\ca L\to\bar{\ca L}$ 
be the quotient map. Set $\bar{\ca F}:=\ca F_{\bar S}(\bar{\ca L})$.
\enddefinition

\vskip .1in 
The following lemma clarifies the relationship between the fusion system $\ca F$ of $\ca L$ and the 
fusion system of $\bar{\ca L}$. 

\proclaim {Lemma 7.6} Assume the setup of 7.5. 
\roster 

\item "{(a)}" $\r$ restricts to a surjective homomorphism $\s:S\to\bar S$, and $\s$ is then a homomorphism 
$\ca F\to\bar{\ca F}$ of fusion systems. Moreover, if $X$ and $Y$ are subgroups of $S$ containing $T$ 
then the induced map 
$$ 
\s_{X,Y}:Hom_{\ca F}(X,Y)\to Hom_{\bar{\ca F}}(\bar X,\bar Y)
$$ 
is surjective. 

\item "{(b)}" Let $X\leq S$ be a subgroup of $S$ containing $T$. Then $X$ is fully normalized in 
$\ca F$ if and only if $\bar X$ is fully normalized in $\bar{\ca F}$. 

\item "{(c)}" Let $\bar X\in\bar{\ca F}^c$, and let $X$ be the preimage of $\bar X$ in $S$. Then 
$X\in\ca F^c$. 

\item "{(d)}" Let $\bar X\in\bar{\ca F}^{cr}$, and let $X$ be the preimage of $\bar X$ in $S$. Then 
$X\in\ca F^{cr}$. 

\item "{(e)}" If $\ca F^{cr}\sub\D$ then $\bar{\ca F}^{cr}\sub\bar\D$. 

\endroster 
\endproclaim 

\demo {Proof} As $\r$ is a projection we have $\bar\D=\{P\r\mid P\in\D\}$ by definition I.5.4. Thus 
$\s:S\to\bar S$ is a surjective homomorphism of groups. For each $w\in\bold W(\ca L)$ and each subgroup 
$X$ of $S_w$, write $X^w$ for the image of $U$ under the composite $c_w$ of the conjugation maps given 
sequentially by the entries of $w$. Similarly define $\bar X^{\bar w}$ when $\bar w\in\bold W(\bar{\ca L})$ 
and $\bar X\leq\bar S_{\bar w}$. By definition, $Hom(\ca F)$ is the set of all mappings $c_w:X\to Y$ with 
$X^w\leq Y\leq S$, and similarly for $Hom(\bar{\ca F})$. If $w=(g_1,\cdots,g_n)$ there is then a 
commutative diagram of homomorphisms: 
$$ 
\CD 
X@>c_w>>Y \\ 
@V\s VV   @VV\s V  \\ 
\bar X@>>c_{w\r^*}>\bar Y 
\endCD 
$$ 
and thus $\s$ is fusion-preserving. Suppose now that $X$ and $Y$ contain $T$, and let 
$\bar w\in\bold W(\bar{\ca L})$ with $\bar X^{\bar w}\leq\bar Y$. By I.3.11 there exists a 
$\r^*$-preimage $w$ of $\bar w$ with $w\in\bold W(N_{\ca L}(T))$. Then $X^w\leq Y$, and 
$\s\mid_X\circ c_{\bar w}=c_w\circ\s\mid Y$ as maps from $X$ to $\bar Y$. That is, the 
$\bar{\ca F}$-homomorphism $c_{\bar w}:\bar X\to\bar Y$ is in the image of $\s_{X,Y}$, and thus 
(a) holds. 

Point (b) is given by the observation that if $X$ is a subgroup of $S$ containing $T$ then 
$\bar{N_S(X)}=N_{\bar S}(\bar X)$. Now let $\bar X\in\bar{\ca F}^c$, let $X$ be the preimage of 
$\bar X$ in $S$, and let $Y$ be an $\ca F$-conjugate of $X$. Then $\bar Y$ is an $\bar{\ca F}$-conjugate 
of $\bar X$, so $C_{\bar S}(\bar Y)\leq\bar Y$, and hence $C_S(Y)\leq Y$. This proves (c). Now 
assume further that $\bar X\in\bar{\ca F}^{cr}$, and let $Y\in X^{\ca F}$ with $Y$ fully normalized in 
$\ca F$. Then $\bar Y$ is fully normalized in $\bar{\ca F}$ by (b). Set $Q=O_p(N_{\ca F}(Y))$. Then 
$\bar Q\leq O_p(N_{\bar{\ca F}}(\bar Y))$ by (a), so $\bar Q=\bar Y$, and hence $Q=Y$. This shows that 
$X\in\ca F^{cr}$ and proves (d). Point (e) is immediate from (d). 
\qed 
\enddemo

If $(\ca L,\D,S)$ is a proper locality and $\ca N\norm\ca L$ then the quotient locality $\ca L/\ca N$ 
need not be proper, just as a homomorphic image of a group of characteristic $p$ need not be of 
characteristic $p$. Thus Theorem A not directly applicable to homomorphic images of proper localities.

\proclaim {Theorem 7.7} Let $(\ca L,\D,S)$ be a proper locality on $\ca F$, let $\ca N\norm\ca L$, let 
$\bar{\ca L}=\ca L/\ca N$ be the quotient locality, and let $\r:\ca L\to\bar{\ca L}$ be the canonical 
projection. Let $\D^+$ be an $\ca F$-closed set of subgroups of $S$ such that $\D\sub\D^+\sub\ca F^s$, 
and let $(\ca L^+,\D^+,S)$ be the expansion of $\ca L$ to $\D^+$. Set 
$\bar{\D}^+=\{P\r\mid P\in\D^+\}$, and set $\bar{\ca F}=\ca F_{\bar S}(\bar{\ca L})$.   
\roster 

\item "{(a)}" There is a locality $(\bar{\ca L}^+,\bar\D^+,\bar S)$ on $\bar{\ca F}$ whose restriction to 
$\bar\D$ is $\bar{\ca L}$, and $\bar{\ca L}^+$ is unique up to a unique isomorphism which restricts to 
the identity map on $\bar{\ca L}$. 

\item "{(b)}" There is a unique projection $\r^+:\ca L^+\to\bar{\ca L}^+$ whose restriction to $\ca L$ 
is $\r$. Moreover, $Ker(\r^+)=\ca N^+$, where $\ca N^+$ is the unique partial normal subgroup of 
$\ca L^+$ whose intersection with $\ca L$ is equal to $\ca N$. 

\endroster 
\endproclaim 

\demo {Proof} Note first of all that since $\ca F^{cr}\sub\D$ we have $\bar{\ca F}^{cr}\sub\bar\D$ by 
7.6(e). Notice also that if $\D=\D^+$ then there is nothing to prove, so we may assume that $\D$ is 
properly contained in $\D^+$.  

Among all $\ca F$-closed sets $\D_1$ with $\D\sub\D_1\sub\D^+$, let $\D_1$ be maximal subject to the 
condition that (a) and (b) hold with $\D_1$ in the role of $\D^+$. Then $\D_1$ is a proper subset of 
$\D^+$, and by replacing $\D$ with $\D_1$ we reduce (as in the proof of Theorem A1) to the case 
where $\D^+=\D\cup R^{\ca F}$ for some $R\leq S$. Thus $\ca L^+$ is an elementary expansion of $\ca L$. 

Take $R$ to be fully normalized in $\ca F$, and suppose that $T\nleq R$. Then $RT\in\D$, and then 
$\D^+\r=\D\r=\bar\D$. Then $(\bar{\ca L},\bar\D,\bar S)$ is the unique locality on $\bar{\ca F}$ whose 
restriction to $\bar\D$ (namely, itself) is $\bar{\ca L}$, and thus (a) holds in this case. In order 
to verify (b) in this case we appeal to Theorem 4.3(c). Thus $\ca L^+$ is the ``free amalgamated product" 
in the category of partial groups of the ``amalgam" given by the inclusion maps $\ca L_0\to\ca L$ and 
$\ca L_0\to\ca L_0^+$. Point (b) will follow if: 
\roster 

\item "{(1)}" There is a unique homomorphism $\r_0^+:\ca L_0^+\to\bar{\ca L}$ such that $\r_0^+$  
agrees with $\r$ on $\ca L_0$. 

\item "{(2)}" $\r_0+$ induces a surjection $\bold D(\ca L_0^+)\to\bold D(\bar{\ca L}_0)$, where 
$\bar{\ca L}_0$ is the set of all $\bar f\in\bar{\ca L}$ such that $\bar S_{\bar f}$ contains a  
$\bar{\ca F}$-conjugate of $\bar R$. 

\item "{(3)}" $Ker(\r^+)=\ca N^+$. 

\endroster 
Versions of the same three points will need to be verified in the case where $T\leq R$, and our 
approach will be to merely sketch the proofs in each case, leaving some of the entirely mechanical details 
to the reader. Thus, let $\Phi$ be the set of triples $(x\i,g,y)$ defined following 4.5, 
let $\approx$ be the relation on $\Phi$ defined in 4.11, and let 
$\r^*:\bold W(\ca L)\to\bold W(\bar{\ca L})$ be the $\r$-induced homomorphism of free monoids. Then 
$\r^*$ maps $\Phi$ into $\bold D(\bar{\ca L})$, and it may be routinely verified that the composition 
$$ 
\Phi@>\r^*>>\bold D(\bar{\ca L})@>\bar\Pi>>\bar{\ca L} 
$$
is constant on $\approx$-classes. Thus, there is a mapping $\r_0^+:\ca L_0^+\to\bar{\ca L}$ given by 
$[x\i,g,y]\maps\bar\Pi(\bar x\i,\bar g,\bar y)$. The verification that $\r_0^+$ is a homomorphism is also 
routine, and yields (1). Since $\r^*$ maps $\bold D(\ca L)$ onto $\bold D(\bar{\ca L})$, (2) is then 
immediate. As $Ker(\r^+)\cap\ca L=\ca N$, (3) is immediate from Theorem A2. Thus, the theorem holds if 
$T\nleq R$.  

Assume henceforth that $T\leq R$. Then $N_{\bar{\ca L}}(\bar R)=\bar{N_{\ca L}(R))}$, since each element of 
$\bar{\ca L}$ has a $\r$-preimage in $N_{\ca L}(T)$ by I.3.11. Then $N_{\bar{\ca L}}(\bar R)$ is a 
subgroup of $\bar{\ca L}$ by I.1.15. Thus Hypothesis 4.2 is satisfied, with $\bar{\ca L}$ and $\bar R$ in 
the roles of $\ca L$ and $R$. Theorem 4.3 then yields (a). For the same reasons as in the preceding 
case (and because the verification of (3) did not in fact make use of the hypothesis that 
$T\nleq R$) it now suffices to verify the following two points.  
\roster 

\item "{(4)}" There is a unique homomorphism $\r_0^+:\ca L_0^+\to\bar{\ca L}^+$ such that $\r_0^+$  
agrees with $\r$ on $\ca L_0$. 

\item "{(5)}" $\r_0+$ induces a surjection $\bold D(\ca L_0^+)\to\bold D(\bar{\ca L}_0^+)$, where 
$\bar{\ca L}_0^+$ is the set of all $\bar f\in\bar{\ca L}^+$ such that $\bar S_{\bar f}$ contains a  
$\bar{\ca F}$-conjugate of $\bar R$.  

\endroster 
Let $\Phi$ and $\r^*$ be as above, and define $\bar\Phi$ to be the set of triples 
$\bar\phi=(\bar x\i,\bar g,\bar y)$ such that one has 
$$ 
\bar U@>\bar x\i>>\bar R@>\bar g>>\bar R@>\bar y>>\bar V 
$$ 
(a sequence of conjugation isomorphisms between subgroups of $\bar S$, labeled by the conjugating elements), 
with $N_{\bar S}(\bar U)\leq\bar S_{\bar x\i}$, and with $N_{\bar S}(\bar V)\leq\bar S_{\bar y\i}$. Then 
$\r^*$ maps $\Phi$ onto $\bar\Phi$ by 5.4(a). It need cause no confusion to denote also by $\sim$ and 
$\approx$ the two equivalence relations on $\bar\Phi$ given by direct analogy with 4.6 and 4.11. 

Again by means of 7.6(a) one verifies that the restriction of $\r^*$ to $\Phi$ preserves these equivalence 
relations, and hence induces a surjective mapping $\r_0^+\ca L_0^+\to\bar{\ca L}_0^+$. Here $\ca L_0$ is 
by definition the set of elements $f\in\ca L$ such that $S_f$ contains an $\ca F$-conjugate of $R$, and 
any such $f$ is identified with its $\approx$-class $[f]$, consisting of all those 
$\phi\in\Phi\cap\bold D(\ca L)$ 
such that $\Pi(\phi)=f$. The analogous definition of $\bar{\ca L}_0$ leads to the conclusion that 
the restriction of $\r_0^+$ to $\ca L_0$ is the (surjective) homorphism $\r_0:\ca L_0\to\bar{\ca L}_0$ 
induced by $\r$. The verification that $\r_0^+$ is a homomorphism, and hence that (4) holds, is then 
straightforward. 

Set $\bold D_0^+=\bold D(\ca L_0^+)$ and $\bar{\bold D}_0^+=\bold D(\bar{\ca L_0}^+)$, let  
$\bar w\in\bar{\bold D}_0^+$, and write $\bar w=([\bar\phi_1],\cdots,[\bar\phi_n])$. Here 
$[\bar\phi_i]$ is the $\approx$-class of an element $\bar\phi_i=(\bar x_i\i,\bar g_i,\bar y_i)$ in 
$\bar\Phi$, and the representatives $\bar\phi_i$ may be chosen so that the sequence 
$\bar\g=(\bar\phi_1,\cdots,\bar\phi_n)$ is a $\bar\G$-form of $\bar w$, as defined in 3.11. That is, 
the word 
$$ 
\bar w_{\bar g}=\bar\phi_1\circ\cdots\circ\bar\phi_n 
$$ 
has the property that $\bar S_{\bar w_{\bar\g}}$ contains a $\bar{\ca F}$-conjugate of $\bar R$. Let 
$\phi_i$ be a $\r^*$-preimage of $\bar\phi_i$ in $\Phi$, set $\g=(\phi_1,\cdots,\phi_n)$, let 
$[\phi_i]$ be the $\approx$-class of $\phi_i$ (in $\ca L_0^+$), and set  
$w=([\phi_1],\cdots,[\phi_n])$. One verifies that $\g$ is a $\G$-form of $w$, and hence that 
$w\in\bold D_0^+$. Since $\r^*$ maps $w_\g$ to $\bar w_{\bar\g}$, it follows that $\r_0^+$ induces a 
surjection as required in (5). Thus (5) holds, and the proof is complete. 
\qed 
\enddemo

The preceding result is essentially a generalization of Theorem A1 to homomorphic images of proper 
localities. Here is the corresponding version of Theorem A2.

\proclaim {Theorem 7.8} Let the hypothesis and the setup be as in the preceding theorem, let 
$\frak N$ be the set of partial normal subgroups of $\ca L$ containing $Ker(\r)$, and let 
$\frak N^+$ be the set of partial normal subgroups of $\ca L^+$ containing $Ker(\r^+)$. Also, let
$\bar{\frak N}$ be the set of partial normal subgroups of $\bar{\ca L}$, and let $\bar{\frak N}^+$ be the 
set of partial normal subgroups of $\bar{\ca L}^+$. Then there is a commutative diagram of 
inclusion-preserving bijections:  
$$ 
\CD 
\frak N       @>\eta>>      \frak N^+ \\ 
@V{\b}VV                      @VV{\b^+}V  \\ 
\bar{\frak N} @>>\bar\eta> \bar{\frak N}^+ 
\endCD 
$$ 
where $\eta$ is the map $\ca N\maps\ca N^+$, $\b$ is the map $\ca N\maps\ca N\r$, and $\b^+$ is the 
map $\ca N^+\maps\ca N^+\r^+$. Moreover:  
\roster 

\item "{(a)}" $Ker(\r^+)=Ker(\r)\eta$, and 

\item "{(b)}" for each $\ca N\in\frak N$, $\bar{\ca N}^+$ is the unique 
partial normal subgroup of $\bar{\ca L}^+$ whose intersection with $\bar{\ca L}$ is $\bar{\ca N}$. 

\endroster 
\endproclaim 

\demo {Proof} By I.5.7 there are inclusion-preserving bijections 
$$ 
\b:\frak N\to\bar{\frak N}\ \ \text{and}\ \ \b^+:\frak N^+\to\bar{\frak N}^+, 
$$ 
by which a partial normal subgroup is sent to its image in $\bar{\ca L}$ (or in $\bar{\ca L}^+$) under 
the canonical projection $\r$ (or $\r^+$). Theorem A2 yields the inclusion-preserving bijection 
$\eta:\frak N\to\frak N^+$ which sends $\ca N\in\frak N$ to the unique $\ca N^+\norm\ca L^+$ 
such that $\ca N^+$ intersects $\ca L$ in $\ca N$. We obtain the above commutative square by 
taking $\bar\eta$ to be $\b\i\circ\eta\circ\b^+$, so it remains to prove (a) and (b), 

We have $Ker(\r)=\ca L\cap Ker(\r^+)$ since $\r$ is the restriction of $\r^+$ to $\ca L$. 
Thus $Ker(\r^+)$ is a partial normal subgroup of $\ca L^+$ which intersects $\ca L$ in $Ker(\r)$. As 
this condition characterized $Ker(\r)\eta$, we have (a). 
 
Fix $\ca N\in\frak N$, set $\bar{\ca N}=\ca N\r$, and set $\bar{\ca N}^+=\ca N^+\r^+$. Thus 
$\bar{\ca N}^+=\bar{\ca N}\bar\eta$. Identify $\bar{\ca L}^+$ with $\ca L^+/Ker(\r^+)$, let 
$g\in\ca N^+$, and let $[g]$ be the maximal coset of $Ker(\r^+)$ in $\ca L^+$ containing $g$. 
Then $[g]=Ker(\r^+)h$ for some $h\in[g]$ such that $h$ is $\up$-maximalin $\ca L^+$ with respect to 
$Ker(\r^+)$. The Splitting Lemma (I.4.11) yields $x\in Ker(\r^+)$ such that $(x,h)\in\bold D(\ca L^+)$, 
$xh=g$, and $S_{(x,h)}=S_g$. Then $(x\i,x,h)\in\bold D(\ca L^+)$, and $\bold D$-associativity yields 
$h=x\i g$. As $Ker(\r)\leq\ca N$, (a) implies that $Ker(\r^+)\leq\ca N^+$. Thus $h\in\ca N^+$, and 
then $[g]\sub\ca N^+$. Now suppose that also $[g]\in\bar{\ca L}$. As $\bar{\ca L}$ is the image of 
$\ca L$ under $\r^+$, we may then assume that $g$ was chosen so that $g\in\ca L$. That is, we may 
take $g\in\ca L\cap N^+$. Then $g\in\ca N$, and this shows that $\bar{\ca N}^+$ intersects 
$\bar{\ca L}$ in $\bar{\ca N}$. The uniqueness of $\bar{\ca N}^+$ for this property is given by the 
commutativity of the displayed diagram, and so (b) holds. 
\qed 
\enddemo

\vskip .2in 
\noindent 
{\bf Section 8: The fusion system of a proper locality} 
\vskip .1in

This section will build on the results of section 3, concerning fusion systems of centric localities, 
using Theorem A1 to make the transition from proper centric localities to proper localities in general. 
As in section 3, many of the ideas here are taken from [5A]. We have relied also on [He2] for the 
devolopment of the ``subcentric" locality associated with a proper locality. 

\vskip .1in 
First, we need the following 
fundamental result, which is a version (for discrete localities) of Theorem 2.2 in [5A].

\proclaim {Theorem 8.1} Let $(\ca L,\D,S)$ be a locality on $\ca F$, and let $\bold P(\ca F)$ be the 
set of all $\ca F$-centric subgroups $P$ of $S$ such that for some $\ca F$-conjugate $Q$ of $P$ 
with $Q$ fully normalized in $\ca F$ $($with respect to the stratification induced from $\ca L)$, 
we have: 
$$ 
Syl_p(Aut_{\ca F}(Q))\neq\nset\ \ \text{and}\ \ Inn(Q)=O_p(Aut_{\ca F}(Q))\cap Aut_S(Q). \tag*   
$$ 
Assume: 
\roster 

\item "{(1)}" $\bold P(\ca F)\sub\D$, and 

\item "{(2)}" $N_{\ca L}(P)$ is of characteristic $p$ for all $P\in\D$. 

\endroster 
Then $\bold P(\ca F)=\ca F^{cr}$, and if $S$ has the normalizer-increasing property then $\ca L$ is a 
proper locality. 
\endproclaim 

\demo {Proof} The first step will be to show: 
\roster 

\item "{(3)}" $\bold P(\ca F)=\ca F^{cr}\cap\D$. 

\endroster 
Let $Q\in\bold P(\ca F)$ be given. In order to show that $Q\in\ca F^{cr}$ we may assume (since both 
$\ca F^{cr}$ and $\bold P(\ca F)$ are $\ca F$-invariant) that $Q$ is fully normalized in $\ca F$ and 
satisfies (*). Set $M=N_{\ca L}(Q)$ and set $R=N_S(Q)$. Here $Q\in\D$ by (1), so $\ca F_R(M)=N_{\ca F}(Q)$,  
and I.3.10 yields $R\in Syl_p(M)$. As $Q$ is $\ca F$-centric, (2) yields $C_M(Q)=Z(Q)$, and thus 
$Aut_{\ca F}(Q)$ may be identified with $M/Z(Q)$. As $R$ is a Sylow $p$-subgroup of $M$ we obtain 
$Aut_R(Q)\in Syl_p(Aut_S(Q)$, and then (*) yields $Q=O_p(M)$. Then $Q=O_p(N_{\ca F}(Q))$as $M$ is of 
characteristic $p$, and thus $Q\in\ca F^{cr}$. This shows $\bold P(\ca F)\sub\ca F^{cr}\cap\D$. 

Next, let $Q\in\ca F^{cr}\cap\D$ be given. In order to show that $Q\in\bold P(\ca F)$ we may assume 
that $Q$ is fully normalized in $\ca F$. With $M$ and $R$ defined as above, we again have 
$Q=O_p(M)=O_p(N_{\ca F}(R))$, and again $Aut_{\ca F}(Q)$ may be identified with $M/Z(Q)$. Then 
$Inn(Q)=O_p(Aut_{\ca F}(Q))$, so $Q$ satisfies (*), and the proof of (3) is complete.  

If $\bold P(\ca F)=\ca F^{cr}$ then (1), (2), and the condition that $S$ have the normalizer-increasing 
property are the same as the conditions (PL1) through (PL3) in the definition (6.7) of proper 
locality, and there is then nothing more to prove. Thus, we may assume that 
$\bold P(\ca F)\neq\ca F^{cr}$, and then (3) yields $\ca F^{cr}\nsub\D$. In particular, 
$\ca F^c\nsub\D$. 

By Zorn's Lemma there exists an expansion $(\ca L^+,\D^+,S)$ of $\ca L$ which is maximal subject to 
$\ca F_S(\ca L^+)=\ca F$, and such that all $\ca L^+$-normalizers of members of $\D^+$ are 
of characteristic $p$. We may then take $\ca L=\ca L^+$. Let $\w{\D}$ be the set of all 
subgroups $X\leq S$ such that $X^*\in\D$. Then $\bold D(\ca L)$ is the set of all 
$w\in\bold W(\ca L)$ such that $S_w$ contains a member of $\w{\D}$, so we may take $\D=\w{\D}$. 
Among all $Q\in\ca F^c$ with $Q\notin\D$, take $Q$ so that $dim(Q)$ is as large as possible. 
take $Q$ to be fully normalized in $\ca F$. As $\D=\w{\D}$ we then have 
$Q\in\Omega$. Any strict overgroup of $Q$ in $S$ is then of greater dimension than $Q$, and is 
therefore in $\D$. We are free to replace $Q$ with any $\ca F$-conjugate of $Q$, so we may assume 
that $Q$ is fully normalized in $\ca F$. 

Set $R=N_S(Q)$, set $\ca M=N_{\ca L}(Q)$, and let $\S$ be the set of strict overgroups of $Q$ in $R$. 
For any $w\in\bold W(\ca M)\cap\bold D(\ca L)$ we have $Q<S_w$ as $S_w\in\D$, and then 
$Q<N_{S_w}(Q)$ by I.3.2. Thus $N_{S_w}(Q)\in\S$. As $\ca M$ is a partial subgroup of $\ca L$ 
we obtain: 
$$ 
\bold D(\ca M)=\bold W(\ca M)\cap\bold D(\ca L)=\{w\in\bold W(\ca M)\mid R\cap S_w\in\S\}. 
$$ 
Since $\S$ is closed under $\ca M$-conjugation into $R$ and is closed with respect to overgroups in $R$, 
we conclude that $(\ca M,\S)$ is an objective partial 
group. As $Q$ is fully normalized, I.3.9 implies that $R$ is a maximal $p$-subgroup of $\ca M$, and 
so $(\ca M,\S,R)$ is a pre-locality. The poset $\Omega_R(\ca M)$ is finite-dimensional by I.2.15, so 
$(\ca M,\S,R)$ is a locality. 

Let $g\in C_{\ca L}(Q)$ and set $X=S_g\cap R$. Then $XX^g\sub S$, and so $[X,g]$ is defined and is a 
$p$-subgroup of $C_S(Q)$. As $Q$ is centric in $\ca F$ we get $[X,g]\leq Z(Q)$, so $X$ is $g$-invariant. 
As $X\in\S\sub\D$, $g$ is an element of the locally finite group $N_{\ca L}(X)$. of characteristic $p$. 
Write $g=ab$ with $a$ of order prime to $p$ and with $b$ a $p$-element. Then $a$ is a power of 
$g$, so $[X,a]\leq Z(Q)$ and then $[X,a]=\1$ by 6.3. Set $Y=O_p(N_{\ca L}(X))$. Then 
$[Y,a]\leq C_Y(X)$, and so $[Y,a]\leq Z(X)$ as $X$ is $\ca F$-centric. Thus $[Y,a]=\1$, and then 
$a=\1$ since $N_{\ca L}(X)$ is of characteristic $p$. Thus $Z(Q)\<g\>$ is a $p$-subgroup of 
$C_{\ca L}(Q)$. Here $C_R(Q)=Z(Q)$, and since $C_{\ca L}(Q)\norm\ca M$ it follows from 
I.4.1(c) that $Z(Q)$ is a maximal $p$-subgroup of $C_{\ca L}(Q)$. We conclude that $g\in Z(Q)$, and so: 
\roster 

\item "{(4)}" $C_{\ca L}(Q)=Z(Q)$. 

\endroster 
We aim next to show:  
\roster 

\item "{(5)}" For each $P\in Q^{\ca F}$ there exists $g\in\ca L$ such that $N_S(P)\leq S_g$ and such that  
$P^g=Q$. 

\endroster 

Since $N_S(P)\in\D$ for $P\in Q^{\ca F}$, (5) is equivalent to saying that $Q$ is normalizer-inductive 
in $\ca F$ (as defined in 2.8). The proof of (5) will then be given by 3.17, in the form in which 
3.17 was formulated via the conditions (A) through (D), once these conditions have been verified to hold 
in the present context. For ease of reference, we repeat those conditions here. 
\roster 

\item "{(A)}" $(\ca L,\D,S)$ is a locality on $\ca F$, with stratification $(\Omega,\star)$ 
induced from $\ca L$, and with $\D\sub\ca F^c$.  

\item "{(B)}" We are given an $\ca F$-closed set $\G$ containing $\D$, such that $\ca F$ is 
$\G$-inductive. 

\item "{(C)}" For each $X\in\G$ such that $X$ is fully normalized in $\ca F$ we are given a 
centric locality $(\ca L_X,\D_X,N_S(X))$ on $N_{\ca F}(X)$, where $\D_X$ is the set of all 
$P\in\D$ such that $X\norm P$, and where $\ca L_X$ is the set of all $g\in N_{\ca L}(X)$ such 
that $N_{S_g}(X)\in\D_X$. Moreover, if $dim(\Omega_X)< dim(\Omega)$ then $\ca L_X$ satisfies the 
conclusion of 3.9. 

\item "{(D)}" We are given $V\in\Omega$ such that $V$ is fully normalized in $\ca F$ and such that 
the set $\G^+:=\G\cup V^{\ca F}$ is $\ca F$-closed. 

\endroster 
Of course, (A) is immediate from the prevailing hypotheses. Taking $\G$ to be $\D$ and taking $V=Q$, 
we then obtain (B) from 3.1, and again (D) is immediate. It thus remains to verify (C). So, let 
$X\in\D$ be given with $X$ fully normalized in $\ca F$, and set $R=N_S(X)$. Then $\D_X$ is the set of 
all overgroups of $X$ in $R$. In particular we have $X\in\D_X$, so $\ca L_X$ is the group $N_{\ca L}(X)$, 
whence $\ca F_R(\ca L_X)=N_{\ca F}(X)$, and thus $(\ca L_X,\D_X,R)$ is a locality on $N_{\ca F}(X)$. 
Since $X\in\ca F^c$ it follows that for each $w\in\bold W(\ca L_X)$ the group $N_{S_w}(X)$ is centric 
in $\ca F$, hence centric also in $N_{\ca F}(X)$. Thus for $X\in\D$ we may view the group $\ca L_X$ as 
a centric locality on $N_{\ca F}(X)$. Theorem 3.9 then yields (C), and so (5) holds. 

The last main step is: 
\roster 

\item "{(6)}" We have $\ca F_R(\ca M)=N_{\ca F}(Q)$. 

\endroster 
Suppose for the moment that (6) holds. 
The mapping $\ca M\to Aut_{\ca F}(Q)$ given by $x\maps c_x$ is then a surjective homomorphism, and the 
quotient locality $\ca M/C_{\ca L}(Q)$ is therefore a group. That is, $\ca M/Z(Q)$ is a group, by (5).  
Since $\bold D(\ca M)$ is the set of all $w\in\bold W(\ca M)$ such that $w$ projects to a member of 
$\bold D(\ca M/Z(Q)$, we obtain $\bold D(\ca M)=\bold W(\ca M)$, and thus $\ca M$ is a subgroup of $\ca L$. 
Hypothesis 4.2 then holds with $Q$ in the role of $R$. By Theorem 4.3 there is an  
expansion $(\ca L^+,\D^+,S)$ of $\ca L$ to a locality on $\ca F$, with $\D^+=\D\cup Q^{\ca F}$ and 
with $N_{\ca L^+}(Q)=\ca M$. 
As this contradicts the maximality of $\D$, it therefore remains only to prove (6).  

Set $\ca F_Q=\ca F_R(\ca M)$. Assuming now that (6) is false, there exists an 
$N_{\ca F}(Q)$-isomorphism $\phi:A\to B$ such that $\phi$ is not an $\ca F_Q$-isomorphism. As 
$\phi$ extends to an $N_{\ca F}(Q)$-isomorphism $AQ\to BQ$ we may assume $Q\leq A$. If $Q<A$ then 
$A\in\D$, so $\phi$ given by conjugation by an element of $\ca M$, and $\phi$ is then an 
$\ca F_Q$-homomorphism. Thus $Q=A$ and $\phi\in Aut_{\ca F}(Q)$. 

Let $w\in\bold W(\ca L)$ with $Q\leq S_w$ and with $\phi=c_w$, and write $w=(g_1,\cdots,g_n)$. 
Set $Q_0=Q$, and for all $i$ with $1\leq i\leq n$ set $Q_i=(Q_{i-1})^{g_i}$. 
Thus $Q_n=Q^w=Q$. By (5) there exist elements $a_i\in\ca L$ such that $N_S(Q_i)\leq S_{a_i}$ and 
with $(Q_i)^{a_i}=Q$, and where we may take $a_0=a_n=\1$. For each $i$ with $1\leq i\leq n$ set 
$u_i=(a_{i-1}\i,g_i,a_i)$, and set 
$$ 
X_i=N_{S_{g_i}}(Q_{i-1})^{a_{i-1}}. 
$$
Then $Q<X_i$, so $X_i\in\D$, and thus $u_i\in\bold D(\ca L)$ via $X_i$. Set $h_i=\Pi(u_i)$ and set 
$w'=(h_1,\cdots,h_n)$. Then $w'\in\bold W(N_{\ca L}(Q))$ and $c_{w'}$ is equal to $c_w$ as 
automorphisms of $Q$. This shows that $\phi$ is an $\ca F_Q$-automorphism, so the proof of (6), and 
of the theorem, is now complete. 
\qed 
\enddemo

\definition {Definition 8.2} Let $\ca F$ be a fusion system on the countable, locally finite $p$-group $S$, 
and assume that for each subgroup $X\leq S$ the group $Aut_{\ca F}(X)$ is countable and locally finite. Let 
$Y$ be a subgroup of $S$. For any $\ca F$-isomorphism $\a:X\to Y$ define $N_\a$ to be the set of all 
$g\in N_S(X)$ such that the automorphism $\a\i\circ c_g\circ\a$ of $Y$ is in $Aut_S(Y)$.
\roster 

\item "{(1)}" $Y$ is {\it fully automized} in $\ca F$ if $Aut_S(Y)$ is a Sylow $p$-subgroup of 
$Aut_{\ca F}(Y)$.  

\item "{(2)}" $Y$ is {\it receptive} in $\ca F$ if for each $X\in Y^{\ca F}$ and each $\ca F$-isomorphism 
$\a:X\to Y$ there exists an extension of $\a$ to an $\ca F$-homomorphism $\bar{\a}:N_{\a}\to N_S(Y)$.  

\endroster 
$\ca F$ is {\it saturated} if:
\roster 

\item "{(3)}" For each subgroup $X\leq S$ there exists an $\ca F$-conjugate 
$Y$ of $X$ such that $Y$ is both fully automized and receptive in $\ca F$, and 

\item "{(4)}" $\ca F$ is stratified. 

\endroster 
\enddefinition

\proclaim {Theorem 8.3} Let $(\ca L,\D,S)$ be a proper locality on $\ca F$, let $A\leq S$ be fully 
normalized in $\ca F$, and let $B$ be fully centralized in $\ca F$.   
\roster 

\item "{(a)}" $\ca F$ is inductive. 

\item "{(b)}" For each subgroup $X\leq S$ the group $Aut_{\ca F}(X)$ is countable and locally finite. 

\item "{(c)}" $A$ is fully automized and receptive in $\ca F$, and $\ca F$ is saturated. 

\item "{(d)}" There exists a proper locality $(\ca L_A,\D_A,N_S(A))$ on $N_{\ca F}(A)$. Moreover, if 
$\D_A=N_{\ca F}(A)^c$ and $\D=\ca F^c$ then $\ca L_A$ may be taken to be contained in $N_{\ca L}(A)$, and 
the inclusion map $\ca L_A\to N_{\ca L}(A)$ to be a homomorphism of partial groups.  

\item "{(e)}" There exists a proper locality $(\ca C_B,\S_B,C_S(B))$ on $C_{\ca F}(B)$. Moreover, if 
$\S_B=C_{\ca F}(B)^c$ and $\D=\ca F^c$ then $\ca C_B$ may be taken to be contained in $C_{\ca L}(B)$, and 
the inclusion map $\ca C_B\to C_{\ca L}(B)$ to be a homomorphism of partial groups.  

\endroster 
\endproclaim 

\demo {Proof} By Theorem A1, the restriction of $\ca L$ to the $\ca F$-closure of $\ca F^{cr}$ is a proper 
locality $\ca L_0$ on $\ca F$, as is the expansion of $\ca L_0$ to any locality $(\ca L_1,\D_1,S)$ with 
$\D_1=\ca F^c$. Without loss of generality we may then take $\ca L=\ca L_1$ and $\D=\ca F^c$. Then 
Theorem 3.9 yields (a), and yields also a centric locality $(\ca L_A,\D_A,N_S(A))$ with an 
inclusion homomorphism $\ca L_A\to N_{\ca L}(A)$, and a centric locality $(\ca C_B,\S_B,C_S(B))$ on 
$C_{\ca F}(B)$ with an inclusion homomorphism $\ca C_B\to C_{\ca L}(B)$. As in 3.7 and 3.8 we 
may take $\D_A=\{P\in\D\mid A\norm P\}$ and $\S_B=\{Q\leq C_S(B)\mid Z(B)\leq Q\ \text{and}\ \ BQ\in\D\}$. 

In order to complete the proofs of (d) and (e) it remains to show that $\ca L_A$ and $\ca C_B$ are proper. 
Since $N_S(A)$ and $C_S(B)$ inherit the normalizer-increasing property from $S$ it then suffices to show 
that $N_{\ca L_A}(P)$ is of characteristic $p$ for all $P\in\D_A$ and that $N_{\ca C_B}(Q)$ is of 
characteristic $p$ for all $Q\in\S_B$. So, let $P\in\D_A$ be given. Then $P\in\D$, and $N_{\ca L_A}(P)$ is 
the group $N_{N_{\ca L}(P)}(A)$, which is of characteristic $p$ by 6.5(b). 
Similarly, if $Q\in\S_B$ then $BQ\in\D$ then $N_{\ca C_B}(Q)$ is the group $C_{N_{\ca L}(BQ)}(B)$. 
Then $N_{\ca C_B}(Q)\norm N_{N_{\ca L}(BQ)}(B)$, and then 6.5(a) shows that $N_{\ca C_B}(Q)$ is of 
characteristic $p$. Thus (d) and (e) hold. 

Next, set $\ca K=C_{\ca L_A}(A)$, $P=C_S(A)A$, and $H=N_{\ca L_A}(P)$. 
Then $\ca K\norm\ca L_A$, and the Frattini Lemma (I.4.10) yields $\ca L_A=H\ca K$. As $A$ is 
fully normalized, (a) implies that $A$ is fully centralized, and then $P\in\ca F^c$ by 3.5(c). 
Thus $P\in\D_A$, so $H$ is a subgroup of $\ca L$, and then 
$H$ is countable and locally finite. Here $Aut_{\ca F}(A)\cong\ca L_A/\ca K\cong H/C_H(A)$, so 
$Aut_{\ca F}(A)$ is countable and locally finite. Since $Aut_{\ca F}(A)\cong Aut_{\ca F}(A')$ for 
every $\ca F$-conjugate $A'$ of $A$, we have point (b). 

We may view $H$ as a locality $(H,\D_A,N_S(A))$ (where finite-dimensionality is given by I.2.15(a)). 
Then $N_S(A)\in Syl_p(H)$ by 3.8, and so $A$ is fully automized in $\ca F$. To show that $A$ is receptive, 
let $A'\in A^{\ca F}$ and let $\a:A'\to A$ be an $\ca F$-isomorphism. By (a) there exists an 
$\ca F$-homomorphism $\phi:N_S(A')\to N_S(A)$ with $A'\phi=A$. Let $\psi$ be the restriction of 
$\phi$ to $A'$, and set $\eta=\psi\i\circ\a$. Thus $\eta\in Aut_{\ca F}(A)$. 
Let $x\in N_S(A')$. Then $c_{x\phi}\in Aut_S(A)$. Set $\b=\psi\circ c_{x\phi}\circ\psi\i$. 
Then $\b\in Aut(A')$, and for any $b\in A'$ we compute: 
$$ 
b\b=((x\i\psi)(b\psi)(x\psi))\psi\i=x\i b x. 
$$ 
Thus $\b=c_x$. Now assume that $x\in N_{\a}$, and compute 
$$ 
\eta\i\circ c_{x\phi}\circ\eta=\a\i\circ(\psi\circ c_{x\phi}\circ\psi\i)\circ\a
=\a\i\circ c_x\circ\a\in Aut_S(A), 
$$ 
which shows that $N_{\a}\phi\leq N_{\eta}$. If $\eta$ extends to an $\ca F$-homomorphism 
$\bar{\eta}:N_{\eta}\to N_S(A)$ then 
$\phi\mid_{N_{\a}}\circ\bar{\eta}$ is an extension of $\a$ to an $\ca F$-homomorphism 
$N_{\a}\to N_S(A)$. Thus, it now suffices to show that $A$ is receptive in $N_{\ca F}(A)$, 
in order to show that $A$ is receptive in $\ca F$. Thus, in proving receptivity for $A$ in $\ca F$ 
we may assume that $A\norm\ca F$, $A\norm\ca L$, and $\a$ is an automorphism of $A$. 

Recall that $Aut_{\ca F}(A)=Aut_H(A)$, where $H=N_{\ca L}(C_S(A)A)$. We may therefore assume that 
$\ca L=H$, a group of characteristic $p$. Here $\a=c_h$ for some $h\in H$, and for $x\in N_{\a}$ 
we have $h\i c_x h=c_y$ for some $y\in S$. Thus $c_{x^h}=c_y$, so $x^h y\i\in C_H(A)$. Thus  
$(N_{\a})^h\leq C_H(A)S$. Here $S\in Syl_p(C_H(A)S)$ by I.3.13(b), so there exists $k\in C_H(A)$ with 
$(N_{\a})^{hk}\leq S$. Then $c_{hk}:N_{\a}\to S$ is an extension of $\a$, and thus $A$ is receptive. 

The requirement (4) in 8.2, that $\ca F$ be stratified, is given by the stratification induced from $\ca L$. 
The proof is thereby complete. 
\qed 
\enddemo

\proclaim {Lemma 8.4} Let $(\ca L,\D,S)$ be a proper locality on $\ca F$ let $\ca N\norm\ca L$ be a 
partial normal subgroup, set $T=S\cap\ca N$, and let $\ca E$ be the fusion system $\ca F_T(\ca N)$. 
Then $P\cap T\in\ca E^c$ for each $P\in\ca F^{cr}$. 
\endproclaim 

\demo {Proof} Fix $P\in\ca F^{cr}$, set $U=P\cap T$, and let $V\in U^{\ca F}$ be fully normalized in 
$\ca F$. As $\ca F$ is inductive, and since $P\leq N_S(U)$, there is an $\ca F$-conjugate $Q$ of $P$ 
with $Q\cap T=V$. Set $H=N_{\ca L}(Q)$. Then $H$ is a subgroup of $\ca L$ of 
characteristic $p$ (as $\ca F^{cr}\sub\D$), and $Q=O_p(H)$ by 6.9(a). 

Set $D=N_{C_T(V)}(Q)$. Then $[Q,D]\leq Q\cap T=V$, and so $[Q,D,D]=\1$. As $Q=O_p(H)$ we then have 
$D\leq Q$ by 6.6(c), and so $D\leq V$. As $C_T(V)$ is $Q$-invariant, and since $C_T(V)Q$ inherits the 
normalizer-increasing property from $S$, we conclude that $C_T(V)\leq V$. For any $\ca F$-conjugate 
$X$ of $V$ there exists an $\ca F$-homomorphism $N_S(X)\to N_S(V)$ which sends $X\to V$, 
and hence $C_T(X)\leq X$. Thus $V$ is centric in $\ca E$, as is every $\ca F$-conjugate of $V$; and 
so $P\cap T\in\ca E^c$ as desired. 
\qed 
\enddemo

The next few results concern the set $\ca F^s$ of $\ca F$-subcentric subgroups of the fusion 
system $\ca F$ of a proper locality. The proofs will be variations of proofs of corresponding 
results in [He2], where $S$ is finite.

\proclaim {Lemma 8.5} Let $\ca F$ be the fusion system on $S$ of a proper locality $\ca L$, and let $V\leq S$ 
be fully centralized in $\ca F$. Then $V\in\ca F^s$ if and only if $O_p(C_{\ca F}(V))$ is centric in 
$C_{\ca F}(V)$. 
\endproclaim 

\demo {Proof} Let $U\in V^{\ca F}$ be fully normalized in $\ca F$. As $\ca F$ is inductive $U$ is 
then fully centralized in $\ca F$, and 1.6 yields an isomorphism $C_{\ca F}(U)\cong C_{\ca F}(V)$ of fusion 
systems.  We may therefore assume to begin with that $V$ is fully normalized in $\ca F$. Set 
$\ca F_V=N_{\ca F}(V)$ and $R=O_p(\ca F_V)$. Also, set $\ca C_V=C_{\ca F}(V)$ and $Q=O_p(\ca C_V)$. 

Suppose first that the lemma holds with $\ca F_V$ in place of $\ca F$. That is, suppose that $R$ is 
centric in $\ca F_V$ if and only if $Q$ is centric in $\ca C_V$. Since $R\in(\ca F_V)^c$ if and only if 
$R\in\ca F^c$ by 3.6(b), we then have the lemma in general. Thus, we are reduced to the case where 
$V\norm\ca F$, and $R=O_p(\ca F)$. 

We are given a proper locality $(\ca L,\D,S)$ on $\ca F$, and by Theorem A1 we may assume that 
$\ca F^c\sub\D$. Suppose that $V\in\ca F^s$. Then $R\in\ca F^c$, so $R\in\D$, and $\ca L$ is the 
group $N_{\ca L}(R)$. Notice that $[R,C_S(VQ)]\leq C_R(VQ)\leq Q$, and that both 
$R$ and $Q$ are normal subgroups of $\ca L$. Then $C_S(VQ)\leq R$ by 6.5(c). Then also 
$$ 
C_{C_S(V)}(Q)=C_S(VQ)\leq C_R(V)\leq Q, 
$$ 
and thus $Q$ is centric in $C_{\ca F}(V)$. On the other hand, assuming that $Q$ is centric in 
$C_{\ca F}(V)$, we obtain 
$$ 
C_S(R)\leq C_S(VQ)\leq Q\leq R, 
$$ 
and so $R\in\ca F^c$, as required. 
\qed 
\enddemo 

\proclaim {Corollary 8.6} Let $\ca F$ be the fusion system of a proper locality. Then 
$$ 
\ca F^{cr}\sub\ca F^c\sub\ca F^q\sub\ca F^s.   
$$
\endproclaim 

\demo {Proof} We have $\ca F^{cr}\sub\ca F^c$ by definition. Let $(\ca L,\D,S)$ be a proper locality on 
$\ca F$, and $P\leq S$ be fully normalized in $\ca F$. 
If $P\in\ca F^c$ then $C_{\ca F}(P)$ is the trivial fusion system on $Z(P)$, and so $P\in\ca F^q$. Now 
suppose instead that we are given $P\in\ca F^q$. Then $O_p(C_{\ca F}(Q))=C_S(Q)$, and then $P\in\ca F^s$ 
by 8.5. 
\qed 
\enddemo

\proclaim {Corollary 8.7} Let $\ca F$ be the fusion system on $S$ of a proper locality, and let  
$V\leq S$ with $V$ fully normalized in $\ca F$. Then 
$$
\{Q\in N_{\ca F}(V)^s\mid V\leq Q\}\sub\ca F^s. 
$$
\endproclaim 

\demo {Proof} One need only observe that $C_{\ca F}(V)=C_{C_{\ca F}(V)}(V)$, in order to obtain the 
desired result from 8.5. 
\qed 
\enddemo

\proclaim {Theorem 8.8} Let $(\ca L,\D,S)$ be a proper locality on $\ca F$. 
\roster 

\item "{(a)}" $\ca F^s$ is $\ca F$-closed. 

\item "{(b)}" If $X\leq S$ with $X^*\in\ca F^s$, then $X\in\ca F^s$. 

\endroster 
\endproclaim 

\demo {Proof} The points (a) and (b) concern only $\ca F$ and not $\ca L$, so by Theorem A1 we may take 
$\D=\ca F^c$. 

We first prove (b). Let $X\leq S$ be given with $X^*\in\ca F^s$. Since 
$\ca F^s$ is $\ca F$-invariant by definition we are free to 
replace $X^*$ by any $\ca F$-conjugate of $X^*$, and thus we may assume to begin with that $X^*$ is 
fully normalized in $\ca F$ and that $O_p(N_{\ca F}(X^*))$ is $\ca F$-centric. 
By 2.4(c) we have $X^*=Y^*$ for any $N_{\ca F}(X^*)$-conjugate $Y$ of $X$, 
so we are free to replace $X$ by any $N_{\ca F}(X^*)$-conjugate of $X$. As $\ca F$ is inductive, 2.11 
shows that we may then take $X$ to be fully normalized in $\ca F$. Set $\ca E=N_{\ca F}(X)$ and set 
$\ca E^*=N_{\ca F}(X^*)$. Then $\ca E=N_{\ca E^*}(\ca E)$. In order to show that $X$ is 
subcentric in $\ca F$ it suffices, by 3.6(b), to show that $O_p(\ca E)$ is centric in $\ca E$, so 
we may assume that $\ca L$ is the locality $\ca L_{X^*}$ given by Theorem 3.9(b). That is, we may 
assume $X^*\norm\ca L$. Then $X^*=O_p(\ca L)\in\ca F^c$, and thus $\ca L$ is a group $G$ of characteristic 
$p$. Then also $N_G(X)$ is of characteristic $p$, by 6.5(b). As $\ca E=\ca F_{N_S(X)}(N_G(X))$ we then have 
$O_p(\ca E)\in\ca E^c$, as desired. We thus have (b).  

In proving (a), note that by definition $\ca F^s$ is $\ca F$-invariant and that 
$S\in\ca F^s$. Thus it only remains to show that $\ca F^s$ is closed with respect to overgroups in $S$. 
Among all $V\in\ca F^s$ such that some overgroup $Q$ of $V$ in $S$ is not subcentric in $\ca F$, choose 
$V$ so that $V$ is fully normalized in $\ca F$, and so that $dim(V)$ is as large as possible. 

Suppose first that $V^*\notin\ca F^s$. Thus, we may take $Q=V^*$. Then $Q$ is not subcentric in 
$N_{\ca F}(V)$, by 8.7. But $V$ is subcentric in $N_{\ca F}(V)$ by 8.6. This shows that the locality 
$\ca L_V$ given by 8.3(d) is a counterexample to (a), and so we may assume $V\norm\ca L$. Then 
$V^*=O_p(\ca L)=O_p(\ca F)=O_p(N_{\ca F}(V))$, so $V^*\in\ca F^c$, contrary to the hypothesis 
that $V^*\notin\ca F^s$. Thus $V^*\in\ca F^s$ after all. Here $Q^\star\notin\ca F^s$ by (b), 
we may take $V=V^\star$ and $Q=Q^\star$. 

As $V<Q$ and $V\Omega$ we have $dim(N_Q(V))>dim(V)$, and the maximality of $dim(V)$ in our choice of $V$ 
then implies that $N_Q(V)\notin\ca F^s$. Thus we may now replace $Q$ with $N_Q(V)$. That is, we may assume 
$V\norm Q$. We are then free to replace $Q$ with any $N_{\ca F}(V)$-conjugate of $Q$, and so we may assume 
that $Q$ is fully normalized in $N_{\ca F}(V)$. 

Suppose that $Q\in N_{\ca F}(V)^s$. As $C_{\ca F}(Q)=C_{N_{\ca F}(V)}(Q)$ we then obtain $Q\in\ca F^s$ from 
8.6. Thus $Q$ is not subcentric in $N_{\ca F}(V)$, and we again reduce to the situation where $V\norm\ca L$ 
and where $Q$ is fully normalized in $\ca F$. 
Here $O_p(\ca F)=O_p(\ca F)\in\ca F^c$ and hence $\ca L$ is a group of characteristic $p$. 
Then $N_{\ca L}(Q)$ is of characteristic $p$ by 6.5(b), and since $Q$ is fully normalized in $\ca F$ 
we have also $N_{\ca F}(Q)=\ca F_{N_S(Q)}(N_{\ca L}(Q))$. 
Then $O_p(N_{\ca F}(Q))\in N_{\ca F}(Q)^c$, and hence $O_p(N_{\ca F}(Q))\in\ca F^c$ by 3.6(b). 
Thus $Q\in\ca F^s$, and we have (a). 
\qed 
\enddemo

The following result is immediate from 8.8 and Theorem A1.

\proclaim {Corollary 8.9} Let $(\ca L,\D,S)$ be a proper locality on $\ca F$. Then there exists an 
expansion of $\ca L$ to a proper locality $(\ca L^s,\ca F^s,S)$ on $\ca F$. Moreover, up to a unique 
isomorphism which restricts to the identity map on $S$, every expansion of $\ca L$ to a proper locality 
on $\ca F$ is a restriction of $\ca L^s$. 
\qed 
\endproclaim 

\definition {Definition 8.10} If $(\ca L,\D,S)$ be a proper locality on $\ca F$ then the locality 
$(\ca L^s,\ca F^s,S)$ is the {\it subcentric closure} of $\ca L$. 
\enddefinition

\proclaim {Lemma 8.11} Let $\ca F$ be the fusion system of a proper locality $(\ca L,\D,S)$, and let 
$P\leq S$ be a subgroup of $S$ such that $O_p(\ca F)P\in\ca F^s$. Then $P\in\ca F^s$. 
\endproclaim 

\demo {Proof} Among all counterexamples, choose $P$ so that $dim(P)$ is as large as possible. Set 
$Q=O_p(\ca F)$, and let $\phi:P\to P'$ be an $\ca F$-isomorphism such that $P'$ is fully normalized in 
$\ca F$. Then $\phi$ extends to an $\ca F$-isomorphism $PQ\to P'Q$, so $P'Q\in\ca F^s$. Thus we may 
assume to begin with that $P$ is fully normalized in $\ca F$. Note that $P^*\notin\ca F^s$ by 8.8(b), while 
$P^*Q\in\ca F^s$ by 8.8(a). Thus $P^*$ is a counterexample to the lemma, and so we may take $P\in\Omega$. 

Set $\ca E=N_{\ca F}(P)$ and set $D=N_Q(P)$. Then $D\leq O_p(\ca E)$. If $D\leq P$ then $P=PR$, and $P$ 
is not a counterexample. Thus $P<DP$, and since $P\in\Omega$ we then have $DP\in\ca F^s$ by the maximality 
of $dim(P)$ in our choice of $P$. Since 
$N_S(P)\leq N_S(PD)$ we may assume that $P$ was chosen so that also $PD$ is fully normalized in $\ca F$. 
Set $R=N_S(PD)$. Also, by 8.10 we may assume that $\D=\ca F^s$, so that $PD\in\D$. 
Then $G:=N_{\ca L}(PD)$ is a group of characteristic $p$, $R\in Syl_p(G)$, and 
$\ca F_{R}(G)=N_{\ca F}(PD)$. Since each $\ca E$-homomorphism extends to 
an $N_{\ca F}(PD)$-homomorphism, $\ca E$ is then the fusion system of $N_G(P)$ at $N_S(P)$. As 
$G$ is of characteristic $p$, so is $N_G(P)$ by 6.5(b). Thus $O_p(\ca E)$ is centric in $\ca E$, 
and hence centric in $\ca F$ by 3.6(b). Thus $P\in\ca F^s$. 
\qed 
\enddemo 

We end this section by showing that our notion of radical centric subgroup is equivalent to a 
more familiar notion from fusion systems over finite $p$-groups, in the case of the fusion 
system of a proper locality. 

\proclaim {Lemma 8.12} Let $(\ca L,\D,S)$ be a proper locality on $\ca F$ and let $P\leq S$ be a 
subgroup of $S$. Assume that $P$ is centric in $\ca F$. Then the following are equivalent. 
\roster 

\item "{(1)}" There exists $Q\in P^{\ca F}$ such that $Q$ is fully normalized in $\ca F$ with respect 
to the stratification induced from $\ca L$ and such that $Q=O_p(N_{\ca F}(Q))$. 

\item "{(2)}" $P=O_p(Aut_{\ca F}(P))$. 

\endroster 
\endproclaim 

\demo {Proof} By Theorem A we may assume without loss of generality that $\ca F^c\sub\D$. Thus $P\in\D$, 
and then $C_{\ca L}(P)=Z(P)$ as $\ca L$ is proper. Thus $Aut_{\ca F}(P)$ may be identified with 
$N_{\ca L}(P)/Z(P)$, and then $O_p(Aut_{\ca F}(P))=O_p(N_{\ca L}(P))/Z(P)$. 
For any $Q\in P^{\ca F}$ there exists $g\in\ca L$ with $Q=P^g$. Then conjugation by $g$ is an 
isomorphism $N_{\ca L}(P)\to N_{\ca L}(Q)$, and we obtain thereby the equivalence of (1) and (2).  
\qed 
\enddemo

\vskip .2in 
\noindent 
{\bf Section 9: $O^p_{\ca L}(\ca N)$ and $O^{p'}_{\ca L}(\ca N)$} 
\vskip .1in 

Throughout this section $(\ca L,\D,S)$ be a locality, $\ca N\norm\ca L$ is a partial normal subgroup, 
and we set $T=S\cap\ca N$. 

\vskip .1in 
\definition {Definition 9.1} Set 
$$ 
\Bbb K=\{\ca K\norm\ca L\mid\ca KT=\ca N\}\ \ \text{and}\ \ \Bbb K'=\{\ca K'\norm\ca L\mid T\leq\ca K\}.   
$$ 
Then set 
$$ 
O^p_{\ca L}(\ca N)=\bigcap\Bbb K\ \ \text{and}\ \ O^{p'}_{\ca L}(\ca N)=\bigcap\Bbb K'. 
$$ 
Write $O^p(\ca L)$ for $O^p_{\ca L}(\ca L)$, and $O^{p'}(\ca L)$ for $O^{p'}_{\ca L}(\ca L)$. 
\enddefinition

\proclaim {Lemma 9.2} Let $(\ca L,\D,S)$ be a locality and let $\ca N\norm\ca L$ be a partial normal 
subgroup, such that every subgroup of $\ca N$ is a $p$-group. Then $\ca N$ is a subgroup of $S$. 
\endproclaim 

\demo {Proof} A partial subgroup of a group is a group, by I.1.8(c), so it suffices to show that $\ca N$ is 
contained in $S$. Suppose false, set $T=S\cap\ca N$, and among all $g\in\ca N$ with $g\notin T$ choose $g$ 
so that $dim(S_g\cap T)$ is as large as possible. Set $P=S_g$ and set $X=P\cap T$. Then 
$N_{\ca N}(P)$ is a subgroup of $\ca N$, and so $N_{\ca N}(P)$ is a $p$-group. By I.3.8 there exists 
$a\in\ca L$ with $N_{\ca N}(P)^a\leq T$, and the Frattini Lemma (I.4.10) shows that we may choose 
$a\in\ca N$.  

Suppose first that 
$X=T$. Then $g\in N_{\ca N}(P)$ by I.4.1(b). As $N_{\ca N}(P)\leq N_{\ca N}(T)$ amnd $N_{\ca N}(P)$ is a 
subgroup of $\ca N$, so $g$ is a $p$-element of $N_{\ca N}(T)$. But then $g\in T$ since $T$ is a maximal 
$p$-subgroup of $\ca N$ by I.4.1(c), so in fact $X<T$. As $P\in\Omega$ we then have $dim(X)<dim(N_T(P))$ by 
I.3.2. Thus $dim(S_a\cap T)>dim(X)$, so $a\in T$, and thus $N_{\ca N}(P)\leq T$. A similar argument 
with $g\i$ in place of $g$ shows that $N_{\ca N}(P^g)\leq T$, so 
$X<N_T(P)\leq S_g\cap T=X$. This contradiction proves the lemma. 
\qed 
\enddemo

\proclaim {Proposition 9.3} Let $(\ca L,\D,S)$ be a locality, let $\ca N\norm\ca L$ be a partial 
normal subgroup, and set $T=S\cap\ca N$. Define $\Bbb K$ and $\Bbb K'$ as in 9.1. Then 
$O^p_{\ca L}(\ca N)\in\Bbb K$ and $O^{p'}_{\ca L}(\ca N)\in\Bbb K'$. 
\endproclaim 

\demo {Proof} It is obvious that $O^{p}_{\ca L}(\ca N)$ and $O^{p'}_{\ca L}(\ca N)$ are partial 
normal subgroups of $\ca L$ contained in $\ca N$, and that $T\leq O^{p'}_{\ca L}(\ca N)$. Thus 
$O^{p'}_{\ca L}(\ca N)\in\Bbb K'$, and it only remains to show that $O^{p}_{\ca L}(\ca N)T=\ca N$. 

Set $\ca K_0=O^p_{\ca L}(\ca N)$, set $\bar{\ca L}=\ca L/\ca K_0$, and adopt the usual ``bar"-convention 
for images of subsets of $\ca L$ under the canonical projection $\ca L\to\bar{\ca L}$. Now let $H$ be a 
subgroup of $\ca N$ and let $x$ be a $p'$-element of $H$. Then $x\in\ca K$ for each $\ca K\in\Bbb K$, 
and so $x\in\ca K_0$. Thus $\bar H$ is a $p$-group. Every subgroup of $\bar{\ca N}$ is of the 
form $\bar H$ for some subgroup $H$ of $\ca N$ by I.2.16 and I.5.3(c), so $\bar{\ca N}=\bar T$ by 9.2. 
Thus $\ca N=\ca K_0T$, as required. 
\qed 
\enddemo

\proclaim {Lemma 9.4} Let ``$*$" be either of the symbols ``$p$" or ``$p'$, let $(\ca L,\D,S)$ be a 
proper locality, and let $(\ca L^+,\D^+,S)$ be an expansion of $\ca L$. Let $\ca N\norm\ca L$, and for 
any partial normal subgroup $\ca K\norm\ca L$ let $\ca K^+$ be the corresponding partial normal subgroup 
of $\ca L^+$ given by Theorem A2. Then $O^{*}_{\ca L}(\ca N)^+=O^{*}_{\ca L^+}(\ca N^+)$. 
\endproclaim 

\demo {Proof} Write $\Bbb K^+$ for the set of all $\ca K^+$ with $\ca K\in\Bbb K$. Then 
$$ 
(\bigcap\Bbb K)^+\leq\bigcap(\Bbb K^+)
$$ 
as $\bigcap(\Bbb K^+)$ is a partial normal subgroup of $\ca L^+$ containing $\bigcap\Bbb K$. The 
reverse inclusion is given by Theorem A2, along with the observation that 
$$ 
\ca L\cap(\bigcap(\Bbb K^+))=\bigcap\{\ca L\cap\ca K^+\}_{\ca K\in\Bbb K}=\bigcap\Bbb K. 
$$
Thus the lemma holds for ``$p$", and the same argument applies to ``$p'$". 
\qed 
\enddemo 

\proclaim {Lemma 9.5} Let $\ca N$ and $\ca M$ be partial normal subgroups of $\ca L$, with $\ca N\leq\ca M$. 
Let ``$*$" be either of the symbols ``$p$" or ``$p'$. Then $O^*_{\ca L}(\ca N)\leq O^*_{\ca L}(\ca M)$.   
\endproclaim 

\demo {Proof} Since $O^{p'}_{\ca L}(\ca M)$ is a partial normal subgroup of $\ca L$ containg $S\cap\ca N$, 
we have $O^{p'}_{\ca L}(\ca N)\leq O^{p'}_{\ca L}(\ca M)$ by definition. Set $\ca K=O^p_{\ca L}(\ca M)$, 
set $\bar{\ca L}=\ca L/\ca K$, and let $\r:\ca L\to\bar{\ca L}$ be the canonical projection. Then 
$(S\cap\ca M)\r=\ca M\r\geq\ca N\r$, and hence $\ca N\r=(S\cap\ca N)\r$. Subgroup correspondence (I.4.7) 
then yields $\ca N\leq\ca K(S\cap\ca N)$, and then $O^p_{\ca L}(\ca N)\leq\ca K$ by definition. 
\qed 
\enddemo

\vskip .2in 
\noindent 
{\bf Section 10: Henke's Theorem on extending homomorphisms} 
\vskip .1in 

In this brief section we adapt [Lemma 3.1 in He3] to our present context. The hypothesis given below is 
in some respects weaker than that in [He3], but Henke's proof (essentially repeated in 10.3, below) remains 
valid in every detail. We have included it for the sake of clarity and completeness. 

\definition {Hypothesis 10.1} We are given localities $(\ca L,\D,S)$ and $(\ca L^+,\D^+,S)$ such that 
$\D\sub\D^+$ and such that $\ca L$ is the restriction of $\ca L^+$ to $\D$. Set 
$\ca F=\ca F_S(\ca L^+)$, and write $(*,\Omega)$ for the stratification on $\ca F$ induced from 
$\ca L^+$.  Assume: 
\roster 

\item "{(1)}" $\D^+=\D\cup R^{\ca F}$, where $R\in\Omega$, and where $R$ is fully normalized in 
$\ca F$ with respect to $(*,\Omega)$. 

\endroster 
Set $M=N_{\ca L^+}(R)$. Assume that there is given a partial group $\w{\ca L}$, a homomorphism 
$\a:\ca L\to\w{\ca L}$, and a homomorphism $\g_M:M\to\w{\ca L}$, such that: 
\roster 

\item "{(2)}" $\a\mid_{N_{\ca L}(R)}=\g_M\mid_{N_{\ca L}(R)}$. 

\endroster 
Notice that (2) yields a mapping $\b:\ca L\cup M\to\w{\ca L}$, and hence also a mapping 
$\b^*:\bold W(\ca L\cup M)\to\bold W(\w{\ca L})$. Assume finally: 
\roster 

\item "{(3)}" If $w\in\bold W(\ca L\cup M)$ and $S_w\in\D^+$, then $w\b^*\in\bold D(\w{\ca L})$. 

\endroster 
\enddefinition 

\proclaim {Theorem 10.2} Assume Hypothesis 10.1. Then there exists a unique homomorphism 
$\g:\ca L^+\to\w{\ca L}$ such that $\g\mid_{\ca L}=\a$ and $\g\mid_M=\g_M$. 
\endproclaim 

\demo {Proof} Write $\Pi:\bold D\to\ca L$, $\Pi^+:\bold D^+\to\ca L^+$, and $\w\Pi:\w{\bold D}\to\w{\ca L}$ 
for the products on $\ca L$, $\ca L^+$, and $\w{\ca L}$, respectively. 

By 3.1, $\ca F$ is $\D^+$-inductive. Thus, for any $P\in R^{\ca F}$ there exists $x_P\in\ca L^+$ such 
that $N_S(P)\leq S_{x_P}$ and with $P^{x_P}=R$. We have $P\in\Omega$ by 10.1(1), so $P<N_S(P)$, and 
then $N_S(P)\in\D$ since $\D^+=\D\cup R^{\ca F}$. Thus $x_P\in\ca L$ for all $P\in R^{\ca F}$. In 
choosing the elements $x_P$, we take $x_R=\1$. 

Let $f\in\ca L^+$ and assume that there exists $P\in R^{\ca F}$ with $P\leq S_f$. Set $Q=P^f$ and set 
$u=(x_P,x_P\i,f,x_Q,x_Q\i)$. Then $u\in\bold D^+$ via $P$, and we have 
$$ 
f=\Pi^+(u)=x_P g x_Q\i, 
$$ 
where $g=\Pi^+(x_P\i,f,x_Q)\in M$. Then $((x_P)\a,g\g_M,(x_Q\i)\a)\in\w{\bold D}$ by 10.1(3). Set 
\roster 

\item "{(4)}" $f\g=\w\Pi((x_P)\a,g\g_M,(x_Q\i)\a)$. 

\endroster  
If $P<S_f$ then $S_f\in\D$ and $u\in\bold D$ via $N_{S_f}(P)$. In that case we have 
$$ 
f\g=\w\Pi((x_P)\a,g\a,(x_Q\i)\a)=\w\Pi(u\a^*)=\Pi(u)\a, 
$$ 
as $\a:\ca L\to\w{\ca L}$ is a homomorphism. Thus $f\g=f\a$ if $S_f$ properly contains a conjugate of $R$. 
For every element $f\in\ca L$ we have $S_f\in\D^+$, so either $f\in\ca L$ or $S_f$ contains a conjugate 
of $R$. We may therefore extend $\g$ to a mapping $\ca L^+\to\w{\ca L}$ by taking $\g\mid_{\ca L}=\a$. 
If $f\in M$ then we may take $P=Q=R$, and since $x_R=\1$ we obtain $f\g=f\g_M$ from (4). Thus,  
we now have $\g\mid_{\ca L}=\a$ and $\g\mid_M=\g_M$. 

We next show that $\g$ is a homomorphism. That is, we will show: 
\roster 

\item "{(5)}" If $w\in\bold D^+$ then $w\g^*\in\w{\bold D}$, and $\w\Pi(w\g^*)=\Pi^+(w)\g$. 

\endroster 
Let $w=(f_1,\cdots,f_n)\in\bold D^+$ be given. If $S_w\in\D$ then $w\in\bold D$, $w\g^*=w\a^*$, 
and then (5) holds since $\a$ is a homomorphism. On the other hand, suppose $S_w\notin\bold D$, so that 
$S_w\in R^{\ca F}$. Set $P_0=S_w$ and recursively define $P_i$ for $1\leq i\leq n$ by 
$P_i=(P_{i-1})^{f_i}$. For $0\leq i\leq n$ set $x_i=x_{P_i}$, and for $1\leq i\leq n$ set 
$u_i=(x_{i-1}\i,f_i,x_i)$. Noting that $u_i\in\bold D^+$ via $P_{i-1}$, set $g_i=\Pi^+(u_i)$, 
and set $v_i=(x_{i-1},g_i,x_i\i)$. Thus $g_i\in M$ and $\Pi^+(v_i)=f_i$. Set 
$v=v_1\circ\cdots\circ v_n$. Then $v\in\bold D^+$ via $P_0$, and $v\in\bold W(\ca L\cup M)$. Then 
10.1(3) yields $v\g^*\in\w{\bold D}$. An appeal to $\w{\bold D}$-associtivity then yields 
$$ 
(\w\Pi(v_1\g^*)\circ\cdots\circ\w\Pi(v_n\g^*))\in\w{\bold D}. \tag6 
$$ 
Here $\w\Pi(v_i\g^*)=(\Pi^+(v_i))\g$ by (4), and $\Pi^+(v_i)=f_i$. Thus (6) may be written as: 
$w\g^*\in\w{\bold D}$. 

We have $\w\Pi(w\g^*)=\w\Pi(v\g^*)$, as we have seen via $\w{\bold D}$-associativity. As $\g\mid_{\ca L}$ 
is a homomorphism, a further application of $\w{\bold D}$-associativity yields 
$$ 
\w\Pi(w\g^*)=\w\Pi(x_0\g,g_1\g,\cdots,g_n\g,x_n\i\g).   
$$ 
Set $g=\Pi^+(g_1,\cdots,g_n)$. As $\g\mid_M$ is a homomorphism we obtain 
$$ 
\w\Pi(w\g^*)=\w\Pi(x_0\g,g\g,x_n\i)=\w\Pi(x_0\a,g\g_M,x_n\i), 
$$ 
and so $\w\Pi(w\g^*)=\Pi^+(x_1,g,x_n\i)\g$ by (4). Since $\Pi^+(x_1,g,x_n\i)=\Pi^+(v)=\Pi^+(w)$ by 
$\bold D^+$-associativity, we have established (5). That is, $\g$ is a homomorphism. We have seen 
that every $f\in\ca L^+$ is a product $\Pi^+(u)$ with $u\in\bold W(\ca L\cup M)$, so $\ca L^+$ is 
generated as a partial group (I.1.9) by $\ca L\cup M$. This observation establishes the uniqueness of 
$\g$, and completes the proof. 
\qed 
\enddemo 

\definition {Definition 10.3} Let $\ca L$ be a partial group, and let $X$ and $Y$ be subsets of 
$\ca L$. Then $X$ is {\it $Y$-invariant} (or $X$ is {\it invariant under} $Y$) if for each $x\in X$ and 
each $y\in Y$ we have $(y\i,x,y)\in\bold D(\ca L)$ and $x^y\in Y$. 
\enddefinition 

\proclaim {Corollary 10.4} Let $(\ca L,\D,S)$ be a proper locality on $\ca F$, and let $(\ca L',\D',S)$ 
be an expansion of $\ca L$ to a proper locality on $\ca F$. Let $(\w{\ca L},\w{\D},\w S)$ be a locality, 
and let $\a:\ca L\to\w{\ca L}$ be a homomorphism. Assume that there exists a subgroup $\w V\leq\w S$ such 
that $\w V$ is invariant under $Im(\a)$ and such that: 
$$ 
\{(P\a)\w V\mid P\in\D'\}\sub\w{\D}. \tag*
$$ 
Then $\a$ extends in a unique way to a homomorphism $\a'\ca L'\to\w{\ca L}$. 
\endproclaim 

\demo {Proof} Set $\Omega=\Omega_S(\ca L)$, and let $(*,\Omega)$ be the stratification on $\ca F$ induced 
from $\ca L'$. For $X\leq S$ write $dim(X)$ for $dim_{\Omega}(X)$, and let $m$ be the maximum, taken over 
all $X\in\D'-\D$, of $dim(X)$. Among all counter-examples to 10.4, assume that $(\ca L,\ca L',\a)$ 
has been chosen so that $m$ is as small as possible. 

Let $\S$ be the set of all subgroups $Y\leq S$ such that $Y^*\in\D$. Then $\D\sub\S$, and $\S$ is a set of 
objects for $\ca L$ by I.3.16. Similarly, the set $\S'$ of all $Y\leq S$ such that $Y^*\in\D'$ is a set of 
objects for $\ca L'$, and we have $\S\sub\S'$. We may therefore take $\D=\S$ and $\D'=\S'$. Thus, we 
have reduced to the case where $X\in\D$ for all subgroups $X\leq S$ such that $X^*\in\D$. 

Choose $R\in\D'$ with $R\notin\D$, and so that $dim(R)=m$. Then $dim(R^*)=m$ (as otherwise $R\in\D$), 
and thus $R\in\Omega$. As $dim(R)=dim(P)$ for $P\in R^{\ca F}$, we may take $R$ to be fully normalized 
in $\ca F$ relative to $(*,\Omega)$. For any $P\in R^{\ca F}$ and any subgroup $X$ of $S$ with $P<S$ we 
have $dim(X)>m$, and so $X\in\D$. Thus the set $\D^+:=\D\cup R^{\ca F}$ is $\ca F$-closed.  
Let $(\ca L^+,\D^+,S)$ be the restriction of $\ca L'$ to $\D^+$. We now have the setup of 10.1(1). 

Set $M=N_{\ca L^+}(R)$. Then $M$ is a subgroup of $\ca L$, by 7.1. Thus 10.1(2) obtains. In order to 
verify 10.1(3), take $w\in\bold W(\ca L)$, set $P=S_w$, and assume $P\in\D^+$. Then $P\in\D'$, and so 
$P\a\w V\in\w{\D}$ by (*). As $\w V$ is $Im(\a)$-invariant, we have $P\a\w V\leq S_{w\a^*}$, and thus 
$w\a^*\in\bold D(\w{\ca L})$. That is, 10.1(3) holds. Now 10.2 yields a unique extension of $\a$ to 
a homomorphism $\g:\ca L^+\to\w{\ca L}$. 

By I.2.3(c) $\w V$ is invariant under all conjugation maps $c_u$, where $u\in\bold W(Im(\a)$. As 
$\ca L^+$ is generated as a partial group by $\ca L$ it follows that $\w V$ is invariant under $Im(\g)$. 
We therefore have the hypothesis of 10.4 fulfilled with $\ca L^+$ and $\g$ in place of $\ca L$ and $\a$. 
As $dim(R)<dim(X)$ for $X\in\D$, we now obtain the desired conclusion via the minimality of $m$. 
\qed 
\enddemo

\enddocument

\vskip .2in 
\noindent 
{\bf Appendix A: $p$-local compact groups as localities} 
\vskip .1in 

Our aim in this appendix is to show that the $p$-local compact groups, introduced in [BLO3] and 
further developed in [BLO3], [JLL] and elsewhere, may be viewed as proper localities of a certain kind, 
(and which we shall call compact localities). As a corollary to a result of Ran Levi and Assaf Libman 
[LL] (and with an improvement due to Molinier [M], removing the reliance of [LL] on the Classification of 
the finite simple groups), we will show that if $\ca F$ is the fusion system of a 
compact locality $(\ca L,\D,S)$ then, up to an 
isomorphism of partial groups, $\ca L$ is the unique compact locality on $\ca F$ having $\D$ as its set of 
objects. Thus, our aim is to form a bridge between the theory being developed here and a theory 
having a much more homotopy-theoretic flavor. The appendix will end with a question concerning the 
sort of conditions that can be placed on the isomorphism establishing the uniqueness of $\ca L$. 
As mentioned in the Introduction, the results obtained here were first suggested by Alex Gonzales, 
and the proof of uniqueness (in A.24) is in large part due to him. 

We shall be closely following the arguments in the appendix 
to [Ch1], where an equivalence was established between the ``$p$-local finite groups" 
introduced in [BLO2] and proper, finite, centric localities. 

Recall from the Appendix A in Part I (or from any of the above references) that a {\it $p$-torus} $T$ is by 
definition the direct product of a finite number of copies of the Pr\" ufer group $\Bbb Z/{(p^\infty)}$. 
A group $G$ is {\it virtually $p$-toral} if there exists a $p$-torus of finite index in $G$. 
A virtually $p$-toral $p$-group is a {\it discrete $p$-toral group}. 

\definition {Definition A.1} The discrete locality $(\ca L,\D,S)$ on $\ca F$ is {\it compact} 
(for the prime $p$) if:  
\roster 

\item "{(1)}" $\ca L$ is proper, and 

\item "{(2)}" the groups $N_{\ca L}(P)$ for $P\in\D$ are virtually $p$-toral. 

\endroster 
\enddefinition 

We will show via A.7, A.22, and A.23 below that there is an equivalence between the notions of  
compact locality and $p$-local compact group. 

\vskip .1in 
The set of elements $x$ in a $p$-torus $T$ such that 
$x^p=\1$ is an elementary abelian $p$-group of finite order $p^k$ where $k$ is equal to 
the number of factors in any decomposition of $T$ as a direct product of Pr\" ufer groups. We 
refer to $k$ as the {\it rank} of $T$, and write 
$$
rk(T)=k.  
$$ 
The identity group is the $p$-torus of rank $0$. 

\proclaim {Lemma A.2} Let $G$ be a virtually $p$-toral group and let $P$ be a discrete $p$-toral group. 
\roster 

\item "{(a)}" There is a unique $p$-torus $T$ such that $T$ has finite index in $G$; and then $T$ contains 
every $p$-toral subgroup of $G$. 

\item "{(b)}" All subgroups and homomorphic images of $G$ are virtually $p$-toral,  
and all subgroups and homomorphic images of $P$ are discrete $p$-toral. 

\item "{(c)}" If $X$ and $Y$ are subgroups of $P$ with $X<Y$ ($X$ is a proper subgroup of $Y$) then 
$X< N_Y(X)$.  

\endroster 
\endproclaim 

\demo {Proof} Let $X$ be a $p$-torus contained in $G$. Then $T\cap X$ has finite index in $X$, and then 
$T\cap X=X$ since $X$ is $p$-divisible. This proves (a). Now let $H$ be a subgroup of $G$. As   
$T\norm G$ by (a), there is an isomorphism $HT/T\cong H/(H\cap T)$, and thus $H$ is discrete 
$p$-toral. Let $N\norm G$. Then $TN/N$ has finite index in $G/N$, and $TN/N\cong T/(N\cap T)$ where 
$T/(N\cap T)$ by $p$-divisibility. This proves (b). Point (c) is given by Lemma A.3(a) in 
the Appendix to Part I. 
\qed 
\enddemo 

We refer to the unique $p$-torus of finite index in the virtually $p$-toral group $G$ as the 
maximal torus of $G$. If $P$ is a discrete $p$-toral group with maximal torus $T$ then (following  
[BLO3]) the {\it order} of $P$ is defined to be the pair 
$$ 
|P|=(rk(P),|P/T|). 
$$ 
If $Q$ is a discrete $p$-toral group with maximal torus $U$ then we write  
$$ 
|P|<|Q| 
$$ 
if $(rk(P),|P/T|)<(rk(Q),|Q/U|)$ lexicographically (i.e. $rk(P)<rk(Q)$, or $rk(P)=rk(Q)$ and $|P/T|<|Q/T|$).

\proclaim {Lemma A.3} Let $P$ and $Q$ be discrete $p$-toral groups. 
\roster 

\item "{(a)}" If $P\cong Q$ then $|P|=|Q|$. 

\item "{(b)}" If $P\leq Q$ then $|P|\leq |Q|$, with equality if and only if $P=Q$. 

\item "{(c)}" If $P\leq Q$ and $Q$ is a $p$-torus, then either $P=Q$ or $rk(P)<rk(Q)$. 

\item "{(d)}" Assume $P\leq Q$, set $P_0=P$, and recursively define $P_k$ for $k>0$ by 
$P_k=N_Q(P_{k-1})$. Set $B=\bigcup\{P_k\}_{k\geq 0}$. Then either $B=Q$ or $rk(B)<rk(P)$. 

\endroster 
\endproclaim  

\demo {Proof} Points (a) and (b) are straightforward. Now let $P\leq Q$ and assume that $Q$ is a 
$p$-torus. For any abelian $p$-group $A$ and any $m>0$ let $A_m$  
be the subgroup of $A$ generated by elements of order dividing $p^m$. Assuming that $P$ is 
a proper subgroup of $Q$, we find $Q_m\nleq P$ for $m$ sufficiently large, and then the rank of 
$P_m$ is less than the rank of $Q_m$. This yields $rk(P)<rk(Q)$, and thus (c) holds. 

In proving (d) let $V$ be the maximal torus of $B$, and let $V_0$ the 
maximal torus of $P_0$. Then $V_0\leq V$. Suppose that $rk(B)=rk(P)$. Then $V_0=V$, so $|B:P|$ is 
finite, and so there exists $k$ with $P_k=P_{k+1}=B$. Then $B=P$ by A.2(c), and thus (d) holds.  
\qed 
\enddemo 

Before we can begin working towards a correspondence between compact localities and $p$-local 
compact groups (and before we can state the definition of $p$-local compact group), we have first to sort 
out some potentially conflicting terminology 
concerning fusion systems. The difficulty here stems from our having already established some terminology in 
section 2 above, relating to stratified fusion systems. This terminology will be shown to agree with  
that of the cited references on $p$-local compact groups 
once the context has been narrowed down to the fusion systems associated with linking 
systems, but until that point has been reached there is a very real problem of confusion.  
For example, if $(\ca F,\Omega,\star)$ is a stratified fusion system on a $p$-group $S$ 
then we have the notion (2.7) of a subgroup $P\leq S$ being fully normalized or fully centralized  
in $\ca F$, but there is a quite different definition of these terms in [BLO3]; and there is a 
similar difficulty concerning saturation of fusion systems. Our solution 
is to slightly alter the terminology from [BLO3], in the manner of the following definition.

\definition {Definition A.4} Let $\ca F$ be a fusion system over the discrete $p$-toral group $S$. 
\roster 

\item "{$\cdot$}" A subgroup $P\leq S$ is {\it fully order-centralized in $\ca F$} if
$|C_S(P)|\geq|C_S(Q)|$ for all $Q\in P^{\ca F}$. 

\item "{$\cdot$}" A subgroup $P\leq S$ is {\it fully order-normalized in $\ca F$} if 
$|N_S(P)|\geq|N_S(Q)|$ for all $Q\in P^{\ca F}$. 

\item "{$\cdot$}" $\ca F$ is {\it order-saturated} if the following three conditions hold. 

\item "{(I)}" For each $P\leq S$ the group $Out_{\ca F}(P)$ is finite. Moreover, if 
$P$ is fully order-normalized in $\ca F$, then $P$ is fully 
order-centralized in $\ca F$, and 
$$Out_S(P)\in Syl_p(Out_{\ca F}(P)).$$  

\item "{(II)}" If $P\leq S$ and $\phi\in Hom_{\ca F}(P,S)$ are such that $P\phi$ is fully order-centralized 
in $\ca F$, and if we set 
$$ 
N_{\phi}=\{g\in N_S(P)\mid \phi\i\circ c_g\circ\phi\in Aut_S(P\phi)\}, 
$$ 
then there exists $\bar\phi\in Hom_{\ca F}(N_\phi,S)$ such that $\phi=\bar\phi\mid_P$. 

\item "{(III)}" If $P_1<P_2<P_3<\cdots$ is an increasing sequence of subgroups of $S$, with 
$P_\infty=\bigcup\{P_n\}_{n=1}^\infty$, and if $\phi:P_\infty\to S$ is a homomorphism such that 
$\phi\mid_{P_n}\in Hom_{\ca F}(P_n,S)$ for all $n$, then $\phi\in Hom_{\ca F}(P_\infty,S)$. 

\endroster 
\enddefinition

By [BLO3, Lemma 1.6], if $\ca F$ is a fusion system over a discrete $p$-toral group $S$, then for any 
subgroup $P\leq S$ there are upper bounds for $|N_S(Q)|$ and for $|C_S(Q)|$ taken over $Q\in P^{\ca F}$. 
Thus $P$ has at least one fully order-normalized $\ca F$-conjugate and at least one fully order-centralized 
$\ca F$-conjugate.

\proclaim {Lemma A.5} Let $(\ca L,\D,S)$ be a compact locality on the order-saturated 
fusion system $\ca F$, and let $P\leq S$ be fully normalized in $\ca F$ with respect to the 
stratification induced from $\ca L$. 
\roster 

\item "{(a)}" $P$ is fully order-normalized in $\ca F$. 

\item "{(b)}" $Inn(P)=O_p(Aut_{\ca F}(P))$ if and only if $P=O_p(N_{\ca F}(P))$. 

\endroster 
\endproclaim 

\demo {Proof} As remarked above, there exists $Q\in P^{\ca F}$ such that 
$Q$ is fully order-normalized in $\ca F$. Set $U=N_S(P)$ and $V=N_S(Q)$. As $\ca F$ is 
order-saturated there then exists an $\ca F$-homomorphism $\phi:U\to V$ with $P\phi=Q$. 

Set $\Omega=\Omega_S(\ca L)$ and let $(*,\Omega)$ be the stratification on $\ca F$ induced from $\ca L$.  
For $X\leq S$ write $dim(X)$ for $dim_{\Omega}(X)$. As $P$ is fully normalized we have 
$dim(U)\geq dim(V)$, and then equality holds since $U\phi\leq V$. Let $\psi$ 
be an extension of $\phi$ to an $\ca F$-homomorphism $U^*\to V^*$. Then $U^*\psi=V^*$ since 
$dim(U^*)=dim(U)$ and $dim(V^*)=dim(V)$. Then also $U\phi=V$ since $U=N_{U^*}(U)$ and 
$V=N_{V^*}(V)$. Thus (a) holds. 

As $\ca F$ is order-saturated there exists a subgroup $\bar P$ of $N_S(P)$ such that $P\leq\bar P$ and 
$Aut_{\bar P}(P)=O_p(Aut_{\ca F}(P))$. Condition (II) in definition A.5 then yields 
$\bar P\norm N_{\ca F}(P)$. So, if $P=O_p(N_{\ca F}(P))$ then $P=\bar P$, and 
thus $Inn(P)=O_p(Aut_{\ca F}(P))$. Now set $R=O_p(N_{\ca F}(P))$ and notice that 
$Aut_R(P)\norm Aut_{\ca F}(P)$. Then $P=R$ if $P=\bar P$, and this completes the proof of (b). 
\qed 
\enddemo

From now on, when speaking of the fusion system $\ca F=\ca F_S(\ca L)$ of a  compact locality, if we  
say that a subgroup $P\leq S$ is fully normalized in $\ca F$ then we mean that $P$ is fully normalized 
with respect to the stratification induced from $\ca L$; and similarly for ``fully centralized".

\proclaim {Lemma A.6} Let $(\ca L,\D,S)$ be a compact locality on $\ca F$. Then $\ca F$ is 
order-saturated. Further: 
\roster 

\item "{(a)}" A subgroup $P\leq S$ is fully order-normalized in $\ca F$ if and only if $P$ is 
fully normalized in $\ca F$. 

\item "{(b)}" A subgroup $P\leq S$ is fully order-centralized in $\ca F$ if and only if $P$ is 
fully centralized in $\ca F$. 

\item "{(c)}" Let $T$ be the maximal torus of $S$.   
Let $R$ be any subgroup of $S$ such that $R$ is a $p$-torus, let $R'$ be an $\ca F$-conjugate of $R$, 
and let $\a:R\to R'$ be an $\ca F$-isomorphism. Then $R$ and $R'$ are subgroups of $T$, and $\a$ 
extends to an $\ca F$-automorphism of $T$. 

\endroster 
\endproclaim 

\demo {Proof} By 8.3 $\ca F$ is saturated, and so $\ca F$ is inductive. Let $P\leq S$ be fully 
order-normalized in $\ca F$, and let $Q\in P^{\ca F}$ be fully normalized. Then there exists an 
$\ca F$-homomorphism $\phi:N_S(P)\to N_S(Q)$ with $P\phi=Q$. 
As $|N_S(P)|\geq|N_S(Q)|$ it follows from A.3 that $N_S(P)\phi=N_S(Q)$. Thus $P$ is fully normalized in 
$\ca F$, and $Q$ is fully order-normalized in $\ca F$. This yields (a), and (b) follows in similar fashion. 

Let $P$ be fully normalized in $\ca F$, with $P\in\ca F^{cr}$. That is, we have $C_S(P)\leq P$ and 
$P=O_p(N_{\ca F}(P))$. 
As $(\ca L,\D,S)$ is proper we may appeal to Theorem 7.2 to obtain a proper expansion 
$(\ca L^+,\D^+,S)$ on $\ca F$, with $\ca F^c\sub\D^+$. By construction, subgroups of $\ca L^+$ are 
conjugates of subgroups of $\ca L$, and thus $\ca L^+$ is compact. We may therefore replace $\ca L$ 
with $\ca L^+$ in the remainder of the proof. That is, we may assume that we have $\ca F^c\sub\D$. 

By 8.3(c) $\ca F$ is receptive, which is to say that $\ca F$ satisfies the condition (II) in definition A.4. 
We now show that $\ca F$ satisfies the condition A.4.I. Let $P$ be as above (i.e. fully order-normalized and 
hence fully normalized). Then $P$ is fully order-centalized. 
Set $X=C_S(P)P$. Then $X$ is $\ca F$-centric, and so $X\in\D$. Set $H=N_{\ca L}(X)$. Every 
$\ca F$-automorphism of $P$ extends to an $\ca F$-automorphism of $X$ by receptivity, so 
$$ 
Aut_{\ca F}(P)=Aut_{N_H(P)}(P).   
$$   
Thus $Aut_{\ca F}(P)$ is a homomorphic image of the virtually $p$-toral group $N_H(P)$, and in particular 
every element of $Aut_{\ca F}(P)$ is of finite order. By [BLO3, Lemma 1.5(b)] every torsion subgroup of 
$Out(P)$ is finite, so we conclude 
that $Out_{\ca F}(P)$ is finite. Here $P$ is fully automized in $\ca F$ 
by 8.3(c), so $Out_S(P)$ is a Sylow $p$-subgroup of $Out_{\ca F}(P)$, and thus A.4(I) holds. 

Let $R$ be any $p$-torus contained in $S$, let $T$ be the maximal torus of $S$, 
and let $\a:R\to R'$ be an $\ca F$-isomorphism. Then $R$ and $R'$ are subgroups of $T$ by A.2(a). 
Let $V\in R^{\ca F}$ be fully centralized in $\ca F$ and since $T\leq C_S(R)$ there exists an 
$\ca F$-automorphism $\b:T\to T$ with $R'\b=V$. Then $\a\circ\b$ is an $\ca F$-isomorphism 
$R\to V$, and $\a\circ\b$ extends to an $\ca F$-automorphism $\g$ of $T$ since $V$ is receptive. 
Then $\g\circ\b\i$ is an extension of $\a$ to $T$, and so (c) holds. 

It now remains to verify (III) in A.4. Let $\s=(P_i)_{i=1}^\infty$ be an increasing sequence of subgroups of 
$S$ and let $\t=(\phi_i)_{i=1}^\infty$ be a sequence of $\ca F$-homomorphisms, $\phi_i:P_i\to S$, such that 
$\phi_i=\phi_{i+1}\mid_{P_i}$ for all $i$. Let $P$ be the union of the groups $P_i$ and let 
$\phi:P\to S$ be the union of the mappings $\phi_i$. For each infinite subset $I$ of natural numbers, 
$\phi$ is then the union of the mappings $\phi_i$ for $i\in I$; so in order to prove that $\phi$ is an 
$\ca F$-homomorphism we are free to replace $\s$ and $\t$ by the sequences corresponding to $I$, and 
then to assume that $I=\Bbb N$. Thus, 
by finite-dimensionality we have $(P_i)^*=P^*$ for $i$ sufficiently large, and so we may assume that  
$(P_i)^*=P^*$ for all $i$. Then each $\phi_i$ extends to an $\ca F$-homomorphism $P^*\to S$, which 
then restricts to an $\ca F$-homomorphism $\psi_i:P\to S$. 

Let $U$ be the maximal torus of $P$. Then $P=UP_n$ for $n$ sufficiently large, and so we may assume 
$P=UP_1$. As $Aut_{\ca F}(T)=Out_{\ca F}(T)$ is finite, (c) implies that there are only 
finitely many $\ca F$-conjugates of $U$, and indeed that there exists an $\ca F$-homomorphism 
$\xi:U\to S$ and an infinite set $I$ of natural numbers such that $\psi_i\mid_U=\xi$ for all $i\in I$. 
As $P=UP_1$ the set $\{\psi_i\}_{i\in I}$ then has a single element $\psi$, and $\psi$ is then the 
union of the homomorphisms $\phi_i$. Thus $\psi=\phi$, so $\phi$ is an $\ca F$-homomorphism, and 
(III) holds. 
\qed 
\enddemo

This completes the preliminaries concerning fusion systems. We next show how a compact locality 
$(\ca L,\D,S)$ on $\ca F$ gives rise to a ``transporter system" of a certain kind, and then that this 
transporter system is a $p$-local compact group provided that $\D$ is the set $\ca F^c$
of $\ca F$-centric subgroups of $S$. The definition of 
transporter system will be taken from [BLO4].

\vskip .1in 
Let $(\ca L,\D,S)$ be any locality. For each $(P,Q)\in\D\times\D$ let $N_{\ca L}(P,Q)$ be the set of 
all $g\in\ca L$ such that $P\leq S_g$ and $P^g\leq Q$. There is then a category $\ca T=\ca T_{\D}(\ca L)$ 
whose set of objects is $\D$, and whose morphisms $g:P\to Q$ are triples $(g,P,Q)$ such that 
$g\in N_{\ca L}(P,Q)$, with composition defined by 
$$ 
(g,P,Q)\circ (h,Q,R)=(gh,P,R). 
$$  
In practice, the role of the objects $P$ and $Q$ will always be clear from the context, and 
we may therefore identify $Mor_{\ca T}(P,Q)$ with $N_{\ca L}(P,Q)$. 

Let $\ca F$ be the fusion system $\ca F_S(\ca L)$, i.e. the fusion system on $S$ generated by the conjugation 
maps $c_g:P\to Q$ with $P,Q\in\D$ and with $g\in N_{\ca L}(P,Q)$. There is then a functor 
$$ 
\r:\ca T\to\ca F 
$$ 
such that $\r$ is the inclusion map $\D\to Sub(S)$ on objects, and such that 
$\r_{P,Q}:N_{\ca L}(P,Q)\to Hom_{\ca F}(P,Q)$ is the map which sends $g$ to the conjugation 
homomorphism $c_g:P\to Q$. We write $\r_P$ for the homomorphism 
$\r_{P,P}:Aut_{\ca T}(P)\to Aut_{\ca F}(P)$. 

Since $(S,\D,S)$ is a locality we have the category $\ca T_{\D}(S)$; and there is a functor 
$$ 
\e:\ca T_{\D}(S)\to\ca T 
$$ 
which is the identity map $\D\to\D$ on objects, and where $\e_{P,Q}$ is the inclusion map 
$N_S(P,Q)\to N_{\ca L}(P,Q)$. We write $\e_P$ for the inclusion $N_S(P)\to N_{\ca L}(P)$. Also, for 
$P,Q\in\D$ with $P\leq Q$, write $\iota_{P,Q}$ for $(\1)\e_{P,Q}$. 

\vskip .1in 
The definition from [BLO4] of a transporter system over a discrete $p$-toral group is 
embedded in the statement of the following result. We remind the reader that we always understand 
composition of morphisms in a category to be taken from left to right.

\proclaim {Proposition A.7} Let $(\ca L,\D,S)$ be a compact locality on $\ca F$ and let 
$$ 
\ca T_{\D}(S)@>\e>> \ca T @>\r>>\ca F 
$$  
be the pair of functors defined above. Then the following hold.  
\roster 

\item "{(A1)}" $\e$ is the identity on objects and $\r$ is the inclusion on objects. 

\item "{(A2)}" For each $P,Q\in\D$ the the group $Ker(\r_P)$ acts freely on $Mor_{\ca T}(P,Q)$ 
from the left (by composition), and $\r_{P,Q}$ is the orbit map for this action. Also, $Ker(\r_Q)$ 
acts freely on $Mor_{\ca T}(P,Q)$ from the right. 

\item "{(B)}" For each $P,Q\in\D$ the map $\e_{P,Q}$ is injective, and $\e_{P,Q}\circ\r_{P,Q}$ sends 
$g\in N_S(P,Q)$ to $c_g\in Hom_{\ca F}(P,Q)$. 

\item "{(C)}" For all $\phi\in Mor_{\ca T}(P,Q)$ and all $x\in P$, the following square commutes in $\ca T$. 
$$
\CD 
P                      @>\phi>>              Q   \\
@V{(x)\e_P}VV                       @VV{(x)(\r(\phi)\circ\e_Q)}V   \\ 
P                      @>\phi>>              Q 
\endCD   
$$

\item "{(I)}" Each $\ca F$-conjugacy class of subgroups in $\D$ contains a subgroup $P$ such that the 
image of $N_S(P)$ under $\e_P$ is a Sylow $p$-subgroup of $Aut_{\ca T}(P)$ (i.e. a subgroup of finite index 
relatively prime to $p$). 

\item "{(II)}" Let $\phi:P\to Q$ be a $\ca T$-isomorphism, and regard conjugation by $\phi$ as a 
mapping $c_\phi:Aut_{\ca T}(P)\to Aut_{\ca T}(Q)$. Let $P\norm\bar P\leq S$ and $Q\norm\bar Q\leq S$, 
be given, and suppose that $c_\phi$ maps $(\bar P)\e_P$ into $(\bar Q)\e_Q$. Then there exists 
$\bar\phi\in Mor_{\ca T}(\bar P,\bar Q)$ such that $\iota_{P,\bar P}\circ\bar\phi=\phi\circ\iota_{Q,\bar Q}$. 

\item "{(III)}" Let $P_1\leq P_2\leq P_3\leq\cdots$ be an increasing sequence of members of $\D$, and 
for each $i$ let $\psi_i:P_i\to S$ be a $\ca T$-homomorphism. Assume that 
$\psi_i=\iota_{P_i,P_{i+1}}\circ\psi_{i+1}$ for all $i$, and set $P=\bigcup\{P_i\}_{i=1}^\infty$. 
Then there exists $\psi\in Mor_{\ca T}(P,S)$ such that $\psi_i=\iota_{P_i,S}\circ\psi$ for all $i$. 

\endroster 
\endproclaim 

\demo {Proof} The condition (A1) is immediate from the definition of the functors $\e$ and $\r$. 
Under the identification of $Aut_{\ca T}(P)$ with $N_{\ca L}(P)$ we have composition in 
$Aut_{\ca T}(P)$ given by group multiplication in $N_{\ca L}(P)$, $Ker(\r_P)=C_{\ca L}(P)$, 
and similarly $Ker(\r_Q)=C_{\ca L}(Q)$. The actions defined in (A2) are then obviously free, and 
since $Aut_{\ca F}(P)\cong N_{\ca L}(P)/C_{\ca L}(P)$ we obtain the conclusion of (A2). 

The condition (B) is again immediate from the definition of $\e$. Now let $g\in N_{\ca L}(P,Q)$ and 
let $x\in P$. Regard $g$ as a $\ca T$-homomorphism $\phi:P\to Q$. Then $(x)\e_P\circ\phi$ is 
simply the product $xg$, while the composition $\phi\circ((x)(\r(\phi)\circ\e_Q))$ is the 
product $gx^g$. As $xg=gx^g$ we have the required commutativity of the diagram in (C). 

Each $\ca F$-conjugacy class of subgroups in $\D$ contains a subgroup $P$ such that $N_S(P)$ 
is a Sylow $p$-subgroup of $N_{\ca L}(P)$, by I.3.10. Thus (I) holds.  

Again let $g\in N_{\ca L}(P,Q)$ be a $\ca T$-isomorphism. Then $P^g=Q$, and $c_g$ is an isomorphism 
$N_{\ca L}(P)\to N_{\ca L}(Q)$. If $P\norm\bar P\leq S$ and $Q\norm\bar Q\leq S$ with 
$\bar P^g\leq\bar Q$, then $g\in N_{\ca L}(\bar P,\bar Q)$, and in this way $g$ is a 
$\ca T$-homomorphism $\bar\phi:\bar P\to\bar Q$. That is, (II) holds. 

Let $(P_i)_{i=1}^\infty$, $(\psi_i)_{i=1}^\infty$, and $P$ be given as in (III). Then $\psi_i$, written 
in full detail, is a triple $(g_i,P_i,S)$ where $P_i\leq S_{g_i}$. The ``inclusion morphism" 
$\iota_{P_i,P_{i+1}}$ is the triple $(\1,P_i,P_{i+1})$, and thus 
$$ 
\iota_{P_i,P_{i+1}}\circ\psi_{i+1}=(g_{i+1},P_i,S). 
$$ 
The hypothesis of (III) therefore translates into the statement that the sequence $(g_i)$ is a constant 
sequence $(g)$ where $P\leq S_g$. Taking $\psi=(g,P,S)$ then yields the conclusion of (III). 
\qed 
\enddemo

\definition {Definition A.8} Let $\ca F$ be a fusion system over the discrete $p$-toral group $S$ and 
let $\D$ be any $\ca F$-closed set of subgroups of $S$. A {\it transporter system} associated to $\ca F$ 
consists of a category $\ca T$ with $Ob(\ca T)=\D$, together with a pair of functors  
$$ 
\ca T_{\D}(S)@>\e>>\ca T\@>\r>>\ca F, 
$$ 
which satisfy the conditions (A1), (A2), (B), (C), (I), (II), and (III) from the preceding proposition 
(where $\r_P:Aut_{\ca T}(P)\to Aut_{\ca F}(P)$ is an abbreviation for $\r_{P,P}$, and where 
$\iota_{P,Q}$ is an abbreviation for $(\1)\e_{P,Q}$ if $P,Q\in\D$ with $P\leq Q$). 
\enddefinition

The definitions of linking system and of $p$-local compact group may now be given as follows, by 
[BLO4, Corollary A.5].

\definition {Definition A.9} The transporter system 
$$ 
\ca T_{\D}(S)@>\e>> \ca T @>\r>>\ca F 
$$  
is a {\it linking system} associated with $\ca F$ if the following conditions hold. 
\roster 

\item "{(1)}" $\ca F$ is order-saturated. 

\item "{(2)}" We have $P\in\D$ for each $\ca F$-centric subgroup $P\leq S$ such that 
$O_p(Out_{\ca F}(P))=1$. 

\item "{(3)}" For each $P$ in $\D$ the kernel of the homomorphism 
$\r_P:Aut_{\ca T}(P)\to Aut_{\ca F}(P)$ is discrete $p$-toral. 

\endroster 
In the special case that $\D$ is the set $\ca F^c$ of all $\ca F$-centric subgroups of $S$ we say 
that $(\e,\r)$ is a {\it $p$-local compact group}. 
\enddefinition

\proclaim {Proposition A.10} Let $(\ca L,\D,S)$ be a compact locality on $\ca F$, with 
$\D=\ca F^c$. Then the transporter system $(\e,\r)$ given by proposition A.7 is a $p$-local compact group. 
\endproclaim 

\demo {Proof} We need only verify the conditions (1) and (3) in the preceding definition, since the 
hypothesis that $\D=\ca F^c$ yields the remaining requirements. Condition (1) is given by A.6. 
Now let $P\in\D$. There is then an isomorphism $\a:N_{\ca L}(P)\to Aut_{\ca T}(P)$ given by 
$g\maps (g,P,P)$, and the kernel of $\r_P$ is then the image under $\a$ of $C_{\ca L}(P)$. 
As $\ca L$ is proper and $P\in\ca F^c$, 6.9 yields $C_{\ca L}(P)=Z(P)$, and so (3) holds. 
\qed 
\enddemo

Our goal now is to proceed in the opposite direction from that of the preceding result. Thus, starting with 
a $p$-local compact group, we aim now to construct a compact locality $(\ca L,\D,S)$.

\vskip .1in 
In what follows we fix the $p$-local finite group $(\ca T_{\D}(S)@>\e>>\ca T\@>\r>>\ca F)$, with the 
abbreviations $\r_P$ and $\iota_{P,Q}$ as earlier. Write $\iota_P$ for the identity morphism $\iota_{P,P}$ 
in $Aut_{\ca T}(P)$. Condition A.8(B) implies that the image of $\iota_{P,P'}$ under $\r$ is the 
inclusion map $P\to P'$, so $\iota_{P,P'}$ is referred to as an {\it inclusion morphism} of $\ca T$. 
This leads to the following definition. 

\definition {Definition A.11} Let $P,Q,P',Q'\in\D$ with $P\leq P'$ and $Q\leq Q'$, and further let 
$\phi\in Mor_{\ca T}(P,Q)$ and $\phi'\in Mor_{\ca T}(P',Q')$. Then $\phi$ is an 
{\it extension} of $\phi'$, and $\phi'$ is a {\it restriction} of $\phi$ if  
$$ 
\iota_{P,P'}\circ\phi=\phi\circ\iota_{Q,Q'}. 
$$
\enddefinition

The following result collects the basic properties concerning the transporter system $(\e,\r)$.

\proclaim {Lemma A.12} 
\roster 

\item "{(a)}" Let $P,Q,R\in\D$, and let 
$$ 
P@>\bar{\phi}>>Q\quad\text{and}\quad Q@>\bar{\psi}>> R
$$ 
be $\ca F$-homomorphisms. Further, let $\psi\in Mor_{\ca T}(Q,R)$ with $\r(\psi)=\bar\psi$, and let 
$\l\in Mor_{\ca T}(P,R)$ with $\r(\l)=\bar\phi\circ\bar\psi$. Then there exists a unique 
$\phi\in Mor_{\ca T}(P,R)$ such that $\r(\phi)=\bar\phi$ and such that $\l=\phi\circ\psi$. 

\item "{(b)}" Let $\psi:P\to Q$ be a $\ca T$-morphism and let $P_0,Q_0\in\D$ with $P_0\leq P$ and 
with $Q_0\leq Q$. Suppose that $\r(\psi)$ maps $P_0$ into $Q_0$. There is then a unique 
$\ca T$-morphism $\psi_0:P_0\to Q_0$ such that $\psi$ is an extension of $\psi_0$. 

\item "{(c)}" A $\ca T$-homomorphism $\phi$ is a $\ca T$-isomorphism if and only if 
$\r(\phi)$ is an $\ca F$-isomorphism. 

\item "{(d)}" All morphisms of $\ca T$ are both monomorphisms and epimorphisms in the categorical 
sense. That is, we have left and right cancellation for morphisms in $\ca T$. 

\item "{(e)}" Let $\phi_0:P_0\to Q_0$ be a $\ca T$-morphism and let $P_0\leq P\leq S$ and 
$Q_0\leq Q\leq S$. Then there exists at most one extension of $\phi_0$ to a $\ca T$-homomorphism 
$P\to Q$. 

\item "{(f)}" Let $P,\bar P,Q,\bar Q$ be objects of $\ca T$, with $P\norm\bar P$ and with $Q\norm\bar Q$. 
Suppose that we are given a $\ca T$-isomorphism $\phi:P\to Q$ and an extension of $\phi$ to a 
$\ca T$-homomorphism $\bar\phi:\bar P\to \bar Q$. Then for each $x\in\bar P$ there is a 
commutative square: 
$$ 
\CD 
P                  @>\phi>>          Q \\ 
@V{x\d_{P,P}}VV                  @VV{y\d_{Q,Q}}V  \\ 
P                  @>\phi>>          Q 
\endCD 
$$ 
where $y$ is the image of $x$ under $\r(\bar\phi)$. 

\item "{(g)}" Every $\ca T$-morphism $\psi:P\to Q$ is the composite of a $\ca T$-isomorphism 
$\phi:P\to Q_0$ followed by an inclusion morphism $\iota_{Q_0,Q}$, where $Q_0$ is the image of $P$ 
under $\r(\psi)$.  

\endroster 
\endproclaim 

\demo {Proof} Points (a) through (d) constitute [BLO4, Proposition A.2], and (e) follows from the 
left cancellation in (d). For (f) one may appeal to the proof of [OV, Lemma 3.3(d)], as that proof 
does not depend on the finiteness of $S$. Finally, let $\phi:P\to Q$ be a $\ca T$-homomorphism and 
let $Q_0$ be the image of $P$ under $\r(\phi)$. Then there exists a restriction of $\phi$ to 
a $\ca T$-homomorphism $\phi_0:P\to Q_0$ by (b), and $\phi_0$ is then a $\ca T$-isomorphism by (c). 
This yields (g). 
\qed
\enddemo  

(A.13) By [BLO3, Section 3] it is a feature of an order-saturated fusion system $\ca F$ over a discrete 
$p$-toral group $S$ that there is a mapping $P\maps P^{\bull}$ from $Sub(S)$ into $Sub(S)$
having the following properties. 
\roster 

\item "{(1)}" $\{P^{\bull}\mid P\leq S\}$ is $\ca F$-invariant, and is 
the union of a finite number of $S$-conjugacy classes of subgroups of $S$. 

\item "{(2)}" For subgroups $P\leq Q\leq S$ we have $P^{\bull}\leq Q^{\bull}$ and 
$(P^{\bull})^{\bull}=P^{\bull}$. 

\item "{(3)}" For all $P,Q\leq S$ we have $N_S(P,Q)\sub N_S(P^{\bull},Q^{\bull})$. 

\item "{(4)}" For all $P,Q\leq S$, each $\ca F$-homomorphism $\a:P\to Q$ extends to an 
$\ca F$-homomorphism $\a^{\bull}:P^{\bull}\to Q^{\bull}$. 

\endroster 
In fact, we will not need (4) here. Rather, what we require is the following result concerning 
$p$-local finite groups.

\proclaim {Lemma A.14} Let $\ca T^{\bull}$ be the full subcategory of $\ca T$ whose set of objects is 
$\{P^{\bull}\mid P\in\D\}$. Then there is a functor 
$$ 
(-)^{\bull}:\ca T\to\ca T^{\bull}, 
$$ 
having the following properties. 
\roster 

\item "{(a)}" $(-)^{\bull}$ is the mapping $P\maps P^{\bull}$ on objects $P\in\D$. 

\item "{(b)}" $(-)^{\bull}$ restricts to the identity functor on $\ca T^{\bull}$. 

\item "{(c)}" For all $P,Q\in\D$ and all $\phi\in Mor_{\ca T}(P,Q)$, the image $\phi^{\bull}$ of 
$\phi$ under $(-)^{\bull}$ is an extension of $\phi$. 

\item "{(d)}" If $\a:X\to Y$ and $\phi:P\to Q$ are $\ca T$-morphisms such that $\phi$ is an 
extension of $\a$, then $\phi^{\bull}$ is an extension of $\a^{\bull}$. 

\endroster 
\endproclaim 

\demo {Proof} Points (a) through (c) are given by [JLL, Proposition 1.12]. By the same reference we 
have also the result that for all $X,P\in\D$ and all $g\in N_S(X,P)$ we have (in accord with A.15(3)) 
$(g)\e_{X,P})^{\bull}=(g)\e_{X^{\bull},P^{\bull}}$. In particular, by taking $X\leq P$ and $g=\1$ we obtain 
$(\iota_{X,P})^{\bull}=\iota_{X^{\bull},P^{\bull}}$. If $\a$ and $\phi$ are given as in (d), so that 
$\iota_{X,P}\circ\phi=\a\circ\iota_{Y,Q}$, the functoriality of $(-)^{\bull}$ now yields 
$\iota_{X^{\bull},P^{\bull}}\circ\phi^{\bull}=\a^{\bull}\circ\iota_{Y^{\bull},Q^{\bull}}$. 
Thus (d) holds. 
\qed 
\enddemo

\proclaim {Lemma A.15} Let $\phi_0:P_0\to Q_0$, $\phi:P\to Q$, and $\phi':P'\to Q'$ be 
$\ca T$-isomorphisms, and suppose that both $\phi$ and $\phi'$ are extensions of $\phi_0$. 
\roster 

\item "{(a)}" If $P=P'$ or if $Q=Q'$, then $\phi=\phi'$. 

\item "{(b)}" There is a unique extension of $\phi_0$ to an isomorphism $\psi:P\cap P'\to Q\cap Q'$, 
and both $\phi$ and $\phi'$ are extensions of $\psi$. 

\endroster 
\endproclaim 

\demo {Proof} Assume that (a) is false. We may take $P=P'$, since the case where $Q=Q'$ will then follow 
by considering the inverses of the given $\ca T$-isomorphisms. Note that if $\phi^{\bull}=(\phi')^{\bull}$ 
then $\phi=\phi'$ by restriction. Since $\phi^{\bull}$ and $(\phi')^{\bull}$ are extensions of 
$(\phi_0)^{\bull}$ by A.14(b), it therefore suffices to consider the case where $\phi_0$, $\phi$, and 
$\psi$ are $\ca T^{\bull}$-isomorphisms. The finiteness condition A.13(1) then yields the existence 
of a counter-example $(\phi_0,\phi,\phi')$ to (a) in which $|P_0|$ is maximal. 

Let $x\in N_P(P_0)$, let $y$ be the image of $x$ under $\r(\phi)$, and let $y'$ be the image of $x$ 
under $\r(\phi')$. We appeal to A.13(f) with $(P_0,Q_0,N_P(P_0),N_Q(Q_0))$ in the role of  
$(P,Q,\bar P,\bar Q)$, and obtain 
$$ 
\phi_0\i\circ(x)\e_{P_0,P_0}\circ\phi_0=(y)\e_{Q_0,Q_0}=(y')\e_{Q_0,Q_0}.  
$$ 
As $\e_{Q_0,Q_0}$ is injective (by A.8(B)) we get $y=y'$, and thus $\r(\phi)$ and $\r(\phi')$ 
agree on $P_1:=N_P(P_0)$. Let $Q_1$ be the image of $P_1$ under $\r(\phi)$. By A.12(b) there is a 
restriction $\phi_1:P_1\to Q_1$ of $\phi$ and a restriction $\phi_1':P_1\to Q_1$ of 
$\phi'$, and then $\phi_1=\phi_1'$ by A.12(e). Now $(\phi_1,\phi,\phi')$ is a counter-example to (a)   
with $|P_1|>|P_0|$, in violation of the maximality of $|P_0|$. This contradiction 
completes the proof of (a). 

Set $X=P\cap P'$ and $Y=Q\cap Q'$. Then $\phi$ and $\phi'$ have 
restrictions $\psi:X\to(X)(\r(\phi))$ and $\psi':X\to(X)(\r(\phi'))$ which, in turn, restrict to 
$\phi_0$. Then (a) yields $\psi=\psi'$, and this establishes (b). 
\qed 
\enddemo 

Define a relation $\up$ on the set $Iso(\ca T)$ $\ca T$-isomorphisms by $\phi\up\phi'$ if 
$\phi'$ is an extension of $\phi$. We may also write $\phi'\down\phi$ to indicate that $\phi$ is a 
restriction of $\phi'$. 

\proclaim {Lemma A.16} The following hold.
\roster

\item "{(a)}" The relation $\up$ induces a partial order on $Iso(\ca T)$.  

\item "{(b)}" The relation $\up$ respects composition of morphisms. That is, if $\phi\up\phi'$ and 
$\psi\up\psi'$, and the compositions $\phi\circ\psi$ and $\phi'\circ\psi'$ are defined, then 
$(\phi\circ\psi)\up(\phi'\circ\psi')$. 

\item "{(c)}" For each $\ca T$-isomorphism $\a$ there exists a unique $\ca T$-isomorphism $\phi$ 
such that $\phi$ is maximal with respect to $\up$ and such that $\a\up\phi$. 

\endroster 
\endproclaim 

\demo {Proof} For points (a) and (b) we repeat the proof of [Che1, Lemma X.7]. The transitivity of  
$\up$ is immediate. Suppose that both $\phi\up\phi'$ and $\phi\down\phi'$, where $\phi\in Iso_{\ca T}(P,Q)$ 
and $\phi'\in Iso_{\ca T}(P',Q')$. Then $P=P'$, $Q=Q'$, $\iota_{P,P'}=\iota_P$, and $\iota_{Q,Q'}=\iota_Q$. 
Further, $\iota_P\phi'=\phi\circ\iota_{Q}$ and then $\phi'=\phi$ since $\iota_P$ and $\iota_Q$ are identity 
morphisms in $\ca T$. Thus (a) holds. 

Suppose that we are given $\phi\up\phi'$ and $\psi\up\psi'$, with $\phi\circ\psi$ and $\phi'\circ\psi'$ 
defined on objects $P$ and $P'$ respectively. Set $Q=P\phi$ and $R=Q\psi$, and set $Q'=P'\phi'$ and 
$R'=Q'\psi'$. The following diagram, in which the vertical arrows are inclusion morphisms, demonstrates that 
$\phi\circ\psi\up\phi'\circ\psi'$. 
$$
\CD
P'         @>\phi'>>    Q'       @>\psi'>>       R' \\ 
@AAA                   @AAA                     @AAA  \\ 
P          @>>\phi>     Q        @>>\psi>        R 
\endCD 
$$ 
This yields (b). 

Let $\a\in Iso(\ca T)$. The finiteness condition A.13(1), together with A.14(c,d), yields the existence 
of at least one $\ca T$-isomorphism $\phi$ such that $\a\up\phi$ and such that $\phi$ is maximal with 
respect to $\up$. Assuming now that $\a$ is a counter-example to (c), there then exists an 
$\up$-maximal $\ca T$-isomorphism $\phi'$ with $\a\up\phi'$ and with $\phi\neq\phi'$. Write 
$\a:X\to Y$, $\phi:P\to Q$, and $\phi':P'\to Q'$. We may again apply A.13(1) in conjunction with 
A.14(c,d) in order to obtain such a triple $(\a,\phi,\phi')$ in which $|X|$ has been maximized. 

Set $P_1=N_P(X)$ and $Q_1=N_{Q}(Y)$, and similarly define $P_1'$ and  $Q_1'$. Set 
$X_1=\<P_1,P_1'\>$ and $Y_1=\<Q_1,Q_1'\>$. Let $\l:Aut_{\ca T}(P_0)\to Aut_{\ca T}(Q_0)$ be the 
isomorphism induced by conjugation by $\a$. Then A.12(f) implies that $\l$ maps 
$(X_1)\e_{P_0}$ to $(Y_1)\e_{Q_0}$. Condition (II) in the definition of transporter system 
then yields the existence of an extension of $\a$ to an isomorphism $\a_1:X_1\to Y_1$. 
Let $\phi_1$ be the restriction of $\phi$ to an isomorphism $P_1\to Q_1$. Then 
$\phi\down\phi_1\up\a_1$. Let $\psi$ be an $\up$-maximal extension of $\a_1$. If $\phi=\psi$ 
then $X_1=X$, whence $P=X=P'$, and then A.15(a) yields $\phi=\phi'$. Thus $\phi\neq\psi$, so 
$(\phi_1,\phi,\psi)$ provides a counter-example to (c). The maximality of $|X|$ then yields 
$X=N_P(X)$, and we again obtain $X=P$. Then $\phi=\a$, $\a$ is $\up$-maximal, and again $\a=\phi'$. 
Thus we have a contradiction, proving (c). 
\qed 
\enddemo

Let $\equiv$ be the equivalence relation on $Iso(\ca T)$ generated by $\up$, and let 
$$ 
\ca L=Iso(\ca T)/\equiv 
$$ 
be the set of equivalence classes. For $\phi\in Iso(\ca T)$ we write $[\phi]$ for the equivalence 
class containing $\phi$.

\proclaim {Lemma A.17} Let $f\in\ca L$. 
\roster 

\item "{(a)}" There is a unique $\phi\in f$ such that $\psi$ is $\up$-maximal in the poset $Iso(\ca T)$. 
Moreover, we then have $\a\up\phi$ for all $\a\in f$, and $\phi\i$ is the unique $\up$-maximal member 
of $[\phi\i]$. 

\item "{(b)}" The unique maximal $\phi\in f$ is a $\ca T^{\bull}$-isomorphism. 

\item "{(c)}" $f\cap Iso_{\ca T}(P,Q)$ has cardinality at most 1 for any 
$P,Q\in Ob(\ca T)$. 

\endroster 
\endproclaim 

\demo {Proof} There exists at least one $\up$-maximal member $\phi\in f$ by A.13(1) with A.14(c,d). Suppose 
that $\phi$ and $\phi'$ are distinct $\up$-maximal members of $f$. 
As $\phi\equiv\phi'$ there is a sequence $\s=(\psi_0,\cdots,\psi_n)$ of members of $f$ with $\phi=\psi_0$, 
$\phi'=\psi_n$, and such that for all $i$ with $1\leq i\leq n$ we have either $\psi_{i-1}\up\psi_i$ or 
$\psi_{i-1}\down\psi_i$. Assume that the pair $(\phi,\phi')$ has been chosen so that $n$ is as small as 
possible. The maximality of $\phi$ implies $\phi\down\psi_1$ and that $\psi_2$ is not a restriction of 
$\phi$. Since $\psi_2$ is the restriction of some maximal isomorphism (necessarily in $f$), we obtain 
$n=2$. Thus $\phi\down\psi_1\up\phi'$. As this violates A.16(c) we obtain the uniqueness asserted in (a). 
Moreover, A.16(c) then implies that each $\a\in f$ extends to $\phi$. The inverse of any extension of 
$\a$ is an extension of $\a\i$, so (a) holds. Point (b) then follows from A.14(c).  

In order to prove (c), let $\psi,\psi'\in f\cap Iso_{\ca T}(P,Q)$. Then both $\psi$ and $\psi'$ are 
restrictions of a single $\phi\in f$, by (a). Now A.12(b) yields $\psi=\psi'$. 
\qed 
\enddemo

Define $\bold D$ to be the set of words $w=(f_1,\cdots,f_n)\in\bold W(\ca L)$ such that there exists a 
sequence $(\phi_1,\cdots,\phi_n)$ of $\ca T$-isomorphisms with $\phi_i\in f_i$, and a sequence 
$(P_0,\cdots,P_n)$ of members of $\D$, such that each $\phi_i$ is a $\ca T$-isomorphism $P_{i-1}\to P_i$.  
As in section I.2 we may say also that $w\in\bold D$ via $(P_0,\cdots,P_n)$, or via $P_0$. Define 
$$
\Pi:\bold D\to\ca L 
$$
by $\Pi(w)=f$, where $f$ is the unique maximal element of $[\phi_1\circ\cdots\circ\phi_n]$ given by 
A.17(a). That $\Pi$ is well-defined follows from A.16(b). Set 
$\1=[\iota_S]$, and for any $f\in\ca L$ let $f\i$ be the equivalence class of 
$\phi\i$, where $\phi$ is the unique maximal member of $f$.

\proclaim {Proposition A.18} $\ca L$ with the above structures is a partial group. Moreover, the 
following hold. 
\roster 

\item "{(a)}" For any $x\in S$, the $\equiv$-class $[(x)\e_S]$ is the set of all $(x)\e_{P,Q}$ such that 
$P^x=Q$, and $(x)\e_S$ is the maximal member of $[(x)\e_S]$. 

\item "{(b)}" $[\iota_S]=\{\iota_P\mid P\in Ob(\ca T)\}$, and $\iota_S$ is the maximal member of its class. 

\item "{(c)}" For any $\phi\in Iso(\ca T)$, $[\phi\i]$ is the set of 
inverses of the members of $[\phi]$. 

\endroster 
\endproclaim 

\demo {Proof} We first check via definition I.1.1 that $\ca L$ is a partial group. For any $f\in\ca L$ 
the word $(f)$ of length 1 is in $\bold D$ since $f$ is represented by a $\ca T$-isomorphism. If $w\in\bold D$ 
and $w=u\circ v$ then it is immediate from the definition of $\bold D$ that both $u$ and $v$ are in 
$\bold D$. Thus the condition I.1.1(1) in the definition of partial group is satisfied. By definition of 
$\Pi$ we have $\Pi(f)=f$ for $f\in\ca L$, so I.1.1(2) holds. Condition I.1.1(3) is a straightforward 
consequence of associativity of composition of isomorphisms in $\ca T$. 

That the inversion map $f\maps f\i$ is an involutory bijection follows from 
A.17(a). Now let $u=(f_1,\cdots,f_n)\in\bold D$ via $(P_0,\cdots,P_n)$, and 
set  $u\i=(f_n\i,\cdots,f_1\i)$. Then $u\i\in\bold D$ via $(P_n,\cdots,P_0)$, 
so $u\i\circ u\in\bold D$. One obtains a representative in the class 
$\Pi(u\i\circ u)$ via a sequence of cancellations $\phi_k\i\circ\phi=\iota_{P_k}$, for 
representatives $\phi_k\in f_i$, so $\Pi(u\i\circ u)$ is the equivalence 
class containing $\iota_{P_0}$. Since $\iota_{P_0}\up\iota_S$, and since 
$\1=[\iota_S]$ by definition, we get $\Pi(u\i\circ u)=\1$. Thus I.1.1(4) holds 
in $\ca L$, and $\ca L$ is a partial group. 

We now prove (a). Let $P\leq P'$ and $Q\leq Q'$ in $\D$, and let $x$ be an element of $S$ such that 
$P^x=Q$ and $(P')^x=Q'$. The functoriality of $\e$ yields 
$$  
(\1)\e_{P,P'}\circ(x)\e_{P',Q'}=(x)\e_{P,Q'}=(x)\e_{P,Q}\circ(\1)\e_{Q,Q'}, 
$$ 
which means that $(x)\e_{P,Q}(g)\up(x)\e_{P',Q'}$. In particular, we get 
$(x)\e_{P,Q}(g)\up(x)\e_S$. In order to complete the proof of (a), it now 
suffices to show that for any $\phi\in Iso_{\ca T}(P,Q)$ with $(x)\e_S\equiv\phi$, we have  
$\phi=(x)\e_{P,Q}$. 
Suppose false, and let $\s=(\phi_1,\cdots,\phi_n)$ be a sequence of 
$\ca T$-isomorphisms with $\phi=\phi_1$, $(x)\e_S=\phi_n$, and with either 
$\phi_i\up\phi_{i+1}$ or $\phi_i\downarrow\phi_{i+1}$ for all $i$ with 
$1\leq i<n$. Among all $(\phi,P,Q)$ with $\phi\neq(x)\e_{P,Q}$ and 
$(x)\e_S\equiv\phi$, choose $(\phi,P,Q)$ so that the length of such a chain 
$\s$ is as small as possible. Set $\psi=\phi_2$. Then $\psi=(x)\e_{X,Y}$, 
where $X$ and $Y$ are objects of $\ca T$ with $X^x=Y$. Suppose 
$\phi\uparrow\psi$. Applying the functor $\r$ to the commutative diagram 
$$ 
\CD 
     X     @>(x)\e_{X,Y}>>    Y  \\
@V\iota_{P,X}VV               @VV\iota_{Q,Y}V \\ 
     P    @>\phi>>            Q
\endCD,
$$
and applying condition (B) in the definition of transporter system to $\r(\e_{X,Y}(g))$, we conclude that 
$\r(\phi)$ is the restriction of $c_g$ to the homomorphism $\r(\phi):P\to Q$. In particular, we get $P^x=Q$, 
so that also $(x)\e_{P,Q}$ is a restriction of $(x)\e_{X,Y}$. Then $\phi=(x)\e_{P,Q}$ by A.12(d), and 
contrary to hypothesis. On the other hand, if $\phi\down\psi$, then $\phi=(x)\e_{P,Q}$ by A.15(a), again 
contrary to hypothesis. This completes the proof of (a), and then (b) is the special case of (a) given
by $x=1$. 

Let $f=[\phi]$ be an equivalence class, with $\phi$ maximal in $f$. One 
checks (by reversing pairs of arrows in the appropriate diagrams) that if 
$\psi$ is a $\ca T$-isomorphism, and $\psi$ is a restriction of $\phi$, then 
the $\ca T$-isomorphism $\psi\i$ is a restriction of $\phi\i$. Point (c) 
follows from this observation. 
\qed 
\enddemo

\definition {Remark} In view of A.18(a), there will be no harm in writing $x$ to denote the equivalence 
class $[(x)\e_S]$, for $x\in S$. That is to say that from now on we shall identify $S$ with the image 
of $S$ under the composition of $\e_S$ with the projection $Iso(\ca T)\to\ca L$. 
\enddefinition

\proclaim {Lemma A.19} Let $\phi:Z\to W$ be a $\ca T$-isomorphism, maximal in 
its $\equiv$-class. Let $X$ and $Y$ be objects of $\ca T$ contained in $Z$, 
and let $U$ and $V$ be the images of $X$ and $Y$, respectively, under 
$\r(\phi)$. Suppose that there exist elements $x$ and $x'$ in $S$ such that 
the following diagram commutes. 
$$
\CD 
X                  @>{\phi\mid_{X,U}}>>            U   \\ 
@V{(x)\e_{X,Y}}VV                         @VV{(x')\e_{U,V}}V   \\ 
Y                  @>{\phi\mid_{Y,V}}>>            V 
\endCD\tag* 
$$ 
Then $x\in Z$, and $x'$ is the image of $x$ under $\r(\phi)$. 
\endproclaim

\demo {Proof} Let $\phi'$ be the composition (in right-hand notation) 
$$
\phi'=(x\i)\e_{Z^x,Z}\circ\phi\circ(x')\e_{W,W^{x'}}.  
$$
Thus, $\phi'\in Iso_{\ca T}(Z^x,W^{x'})$, and the commutativity of (*) 
yields $\phi\equiv\phi'$. The maximality 
of $\phi:Z\to W$ implies that $Z^x\leq Z$ and $W^{x'}\leq W$. That is, 
$x\in N_S(Z)$ and $x'\in N_S(W)$. There is then a commutative diagram as follows. 
$$
\CD
Z                 @>{\phi}>>            W   \\
@V{(x)\e_{Z}}VV                   @VV{(x')\e_{W}}V   \\
Z                 @>{\phi}>>            W
\endCD
$$
Condition (II) in the definition of transporter system implies that there is an extension of $\phi$ to a 
$\ca T$-isomorphism $\<Z,x\>\to\<W,x'\>$, and the maximality of $\phi$ then yields $x\in Z$ and $x'\in W$. 
Condition (C) in the definition of transporter system implies that $x'$ is the image under $\r(\phi)$ of $g$. 
\qed
\enddemo

\proclaim {Corollary A.20} Let $f\in\ca L$ and let $P\in\D$ with the property 
that, for all $x\in P$, $(f\i,x,f)\in\bold D$ and $\Pi(f\i,x,f)\in S$. Let 
$Q$ be the set of all such products $\Pi(f\i,x,f)$. Then $Q\in\D$ and there 
exists $\psi\in f$ such that $\psi\in Iso_{\ca T}(P,Q)$. 
\endproclaim 

\demo {Proof} As $(f\i,x,f)\in\bold D$ there exist $U,X,Y,V\in\D$ and 
representatives $\psi$ and $\bar\psi$ of $f$ such that 
$$ 
\CD 
U@>\bar\psi\i>>X@>(x)\e_{X,Y}>>Y@>\psi>>V 
\endCD 
$$ 
is a chain of $\ca T$-isomorphisms, and where the middle arrow in the 
diagram is given by A.18(a). As $\Pi(f\i,x,f)\in S$ there exists $x'\in S$ such that 
$\bar\psi\i\circ(x)\e_{X,Y}\circ\psi=(x')\e_{U,V}$. Let $\phi:Z\to W$ be the 
maximal element of $f$. Then A.19 implies that $x\in Z$, and 
$x'$ is the image of $x$ under $\r(\phi)$. In particular, we 
have $P\leq Z$ and $Q\leq W$, and we may therefore take $X=Y=P$ 
and $U=V=Q$, obtaining $\psi\in Iso_{\ca T}(P,Q)$. 
\qed  
\enddemo 

\proclaim {Lemma A.21} Let $\psi:P\to Q$ be a $\ca T$-isomorphism, and let $f=[\psi]$ be the 
equivalence class of $\psi$. Then $(f\i,x,f)\in\bold D$ for all $x\in P$, and $P^f=Q$ in the partial 
group $\ca L$. Moreover, the conjugation map $(f\i,x,f)\maps\Pi(f\i,x,f)$ is equal to $\r(\psi)$. 
\endproclaim

\demo {Proof} For any $x\in P$, we have the composable sequence 
$$
Q@>\psi\i>>P@>(x)\e_P>>P@>\psi>>Q  
$$
of $\ca T$-isomorphisms, so $(f\i,x,f)$ is in $\bold D$. By Condition (C) in the definition of transporter 
system then yields $\psi\i\circ(x)\e_P\circ\psi=(x')\e_Q$, where $x'=(x)(\r(\psi))\in Q$. 
The class $[(x')\e_Q]$ is the same as $[(x')\e_S]$ by 
A.18(a); and we recall that we have introduced the convention to denote this 
class simply as $x'$. Thus $x^f=x'$, and so $P^f\sub Q$. Similarly $Q^{f\i}\sub P$, from which one 
deduces that the conjugation map $c_f:P\to Q$ is surjective. Injectivity of $c_f$ follows from 
left and right cancellation in the partial group $\ca L$, so $P^f=Q$. The final assertion of the 
lemma is given by the observation, made above, that $x'=(x)(\r(\psi))$. 
\qed 
\enddemo

\proclaim {Theorem A.22} Let $(\ca T_{\D}(S)@>\e>> \ca T @>\r>>\ca F)$ be a $p$-local compact group, 
and let $\ca L=Iso(\ca T)/\equiv$ be the partial group given by A.18. For $x\in S$, identify $x$ with 
the $\equiv$-class of the $\ca T$-automorphism $(x)\e_S$ of $S$. Then $(\ca L,\D,S)$ is a 
compact locality on $\ca F$. 
\endproclaim 

\demo {Proof} We have seen in A.18 that $\ca L$ is a partial group, and the remark following A.18 
shows how to identify $S$ with its image under the composition of $\e_S$ with the quotient map 
$Iso(\ca T)\to\ca L$. In order to show that $(\ca L,\D)$ 
is objective, begin with $w=(f_1,\cdots,f_n)\in\bold D$. By definition, there exist representatives 
$\psi_i$ of the classes $f_i$, and a sequence $(P_0,\cdots,P_n)$ of objects in $\D$, such that each 
$\psi_i$ is a $\ca T$-isomorphism $P_{i-1}\to P_i$. Then $P_{i-1}^{f_i}=P_i$ for all $i$, 
by A.21. Conversely, given $w=(f_1,\cdots,f_n)\in\bold W(\ca L)$, and given 
$(P_0,\cdots,P_n)\in\bold W(\D)$ with $P_{i-1}^{f_i}=P_i$ for all $i$, it follows from A.20 that 
$w\in\bold D$. Thus, $(\ca L,\D)$ satisfies the condition (O1) in the definition I.2.1 of 
objective partial group. Since $\D=\ca F^c$ is $\ca F$-closed we also have (O2), and thus $(\ca L,\D)$ 
is objective. As $\D$ is a set of subgroups of $S$, $(\ca L,\D,S)$ is then a pre-locality (as defined 
in I.2.5). 

The conjugation maps $c_g:S_g\to S$ for $g\in\ca L$ are $\ca F$-homomorphisms, by A.20 and A.8(C). Thus 
$\ca F_S(\ca L)$ is a subsystem of $\ca F$. Assuming now that $\ca F\neq\ca F_S(\ca L)$, there exists 
an $\ca F$-isomorphism $\b:X\to Y$ such that $\b$ is not an $\ca F_S(\ca L)$-homomorphism. By A.13(4) 
we may take $X=X^{\bull}$ and $Y=Y^{\bull}$, and by the finiteness condition in A.13(1) we may then 
assume that from among all $\ca F$-isomorphisms which are not $\ca F_S(\ca L)$-homomorphisms, 
$\b$ has been chosen so that $|X|$ is as large as possible. If $X\in\D$ then $Y\in\D$ and the surjectivity 
of $\r_{X,Y}$ (condition (A2) in definition A.8) implies that $\b$ is an $\ca F_S(\ca L)$-homomorphism, 
so in fact $X\notin\D$. In particular $X<S$, and so $X<N_S(X)$. Similarly $Y<N_S(Y)$.  

As $\ca F$ is order-saturated there exists a fully order-normalized $\ca F$-conjugate $Z$ of $X$, 
and there then exist $\ca F$-homomorphisms $\eta_1:N_S(X)\to N_S(Z)$ 
and $\eta_2:N_S(Y)\to N_S(Z)$ such that $X\eta_1=Z=Y\eta_2$. Each $\eta_i$ is an 
$\ca F_S(\ca L)$-homomorphism by the maximality in the choice of $X$, and it then suffices to show 
that the $\ca F$-automorphism $\a=\eta_1\i\circ\b\circ\eta_2$ of $Z$ is an $\ca F_S(\ca L)$-homomorphism. 
As $\ca F$ is order-receptive, $\a$ extends to an $\ca F$-automorphism $\bar{\a}$ of 
$C_S(Z)Z$. But $C_S(Z)Z$ is centric in $\ca F$, so $C_S(Z)Z\in\D$, and $\bar{\a}$ is then 
an $\ca F_S(\ca L)$-homomorphism. The same is then true of $\a$, and so we have shown that 
$\ca F=\ca F_S(\ca L)$. 

Set $\G=\{P^{\bull}\mid P\leq S\}$. For $P,Q\in\G$ we have (by A.13(2)) 
$$ 
(P\cap Q)^{\bull}\leq P\cap Q^{\bull}=P\cap Q, 
$$ 
and thus $\G$ is closed under finite intersections. Let $g\in\ca L$ and let $\psi:P\to Q$ be the unique 
$\up$-maximal representative of $g$. Then $P=S_g$ by A.21, and then $S_g\in\G$ by A.14(c). We may 
prove by induction on the length of $w\in\bold W(\ca L)$ that $S_w\in\G$. Namely, write 
$w=(g)\circ v$ where $g\in\ca L$ and $v\in\bold W(\ca L)$ with $S_v\in\G$. Then 
$S_w=(S_{g\i}\cap S_v)^{g\i}\in\G$ since, as we have seen, $\G$ is closed with respect to finite 
intersections and since (by A.13(1)) $\G$ is $\ca F$-invariant. The finiteness condition in 
A.13(1) now implies that $\bigcap\{S_w\mid w\in \bold X\}\in\G$ for each non-empty subset $\bold X$ 
of $\bold W(\ca L)$. This shows that the poset $\Omega_S(\ca L)$ defined in I.2.10 is finite-dimensional. 

Let $P\in\D$ and let $\a_P:Aut_{\ca T}(P)\to N_{\ca L}(P)$ be the mapping $\psi\maps[\psi]$. Then 
$\a_P$ is a homomorphism by A.16(b), $\a_P$ is injective by A.18(c), and $\a_P$ is surjective by A.20. 
 Thus $\a_P$ is an isomorphism, and then condition (I) in definition A.8 implies that $N_{\ca L}(P)$ is 
virtually $p$-toral. As $\Omega_S(\ca L)$ is finite-dimensional, all subgroups of $\ca L$ are then virtually 
$p$-toral by I.2.16. In particular, if $\bar S$ is a $p$-subgroup of $\ca L$ containing $S$ then $\bar S$ 
is discrete $p$-toral. If $S<\bar S$ then A.2(c) yields $S<N_{\bar S}(S)$, which is contrary to condition (I) 
in A.4. Thus $S$ is a maximal $p$-subgroup of $\ca L$, and we have established that $\ca L$ is a  
locality on $\ca F$. Notice that A.21 implies that $\a_P$ restricts to an isomorphism 
$Ker(\r_P)\to C_{\ca L}(P)$. As $(\e,\r)$ is a $p$-local compact group, $Ker(\r_P)$ is a $p$-group, 
and thus $N_{\ca L}(P)$ is of characteristic $p$. That is, $\ca L$ satisfies the condition (PL2) in 
the definition (6.7) of proper locality. Condition (PL1), that $\ca F^{cr}$ be contained in $\D$, is 
given by $\D=\ca F^c$. Condition (PL3), that $S$ has the normalizer-increasing property, is given by 
A.2(c). Thus $\ca L$ is proper, and the proof is complete. 
\qed 
\enddemo

\proclaim {Theorem A.23} Let $(\ca L,\D,S)$ be a compact locality on $\ca F$, such that $\D$ 
is the set $\ca F^c$ of $\ca F$-centric subgroups of $S$. Let $(\ca T_{\D}(S)@>\e>> \ca T @>\r>>\ca F)$ 
be the $p$-local compact group constructed from $\ca L$ as in A.7 and A.10, and let $(\ca L',\D,S)$ be 
the compact locality constructed from $(\e,\r)$ as in A.18 and A.22. Then the mapping 
$$ 
\Phi:\ca L\to\ca L', 
$$  
which sends $g\in\ca L$ to the $\equiv$-class of the $\ca T$-isomorphism $(g,S_g,S_{g\i})$ (and with 
the identifications given by the remark following A.18) is an isomorphism of partial groups which restricts 
to the identity map on $S$. 
\endproclaim 

\demo {Proof} Let $\Phi^*:\bold W(\ca L)\to\bold W(\ca L')$ be the mapping induced by $\Phi$, and let 
$w=(g_1,\cdots,g_n)\in\bold D(\ca L)$ via $(P_0,\cdots,P_n)$, with $P_0=S_w$. We shall denote the 
$\equiv$-class of a $\ca T$-isomorphism $(g,P,Q)$ by $[g,P,Q]$. Then 
$w\Phi^*=([g_1,P_0,P_1],\cdots,[g_n,P_{n-1},P_n])$, and $w\Phi^*\in\bold D(\ca L')$ via $(P_0,\cdots,P_n)$ 
by A.21. The definition of the product $\Pi'$ in $\ca L'$ then yields 
$$ 
\Pi'(w\Phi^*)=[\pi(w),P_0,P_n]=(\Pi(w))\Phi, 
$$ 
and thus $\Phi$ is a homomorphism of partial groups. 

Recall that for $P,P'\in\D$ with $P\leq P'$, we have $\iota_{P,P'}=(\1,P,P')$. It follows that the extensions 
of a $\ca T$-isomorphism $(f,P,Q)$ are of the form $(f,P',Q')$, and this implies that $\Phi$ is injective. 
Since $\ca L'=Im(\Phi)$ by definition, $\Phi$ is a bijection, and we now leave it to the reader to verify 
that $\Phi\i$ is a homomorphism. For $x\in S$ we have identified $x$ with $[x,S,S]$, so $\Phi$ restricts to 
the identity map on $S$. 
\qed
\enddemo

\proclaim {Theorem A.24} Let $(\ca L,\D,S)$ and $(\ca L',\D,S)$ be compact localities on $\ca F$, having 
the same set of objects. Then there exists an isomorphism $\a:\ca L\to\ca L'$ of partial groups, 
such that $P\a=P$ for all $P\in\D\cap\ca F^c$. In particular, $\a$ restricts to an automorphism of $S$. 
\endproclaim 

\demo {Proof} By Theorem A1 we may assume without loss of generality that $\D$ is equal to the set $\ca F^c$ 
of $\ca F$-centric subgroups of $S$. Let $(\ca T_{\D}(S)@>\e>>\ca T\@>\r>>\ca F)$ be the $p$-local 
compact group constructed from $(\ca L,\D,S)$ via A.7 and A.10, and let 
$(\ca T_{\D}(S)@>\e'>>\ca T'\@>\r'>>\ca F)$ be the $p$-local compact group similarly constructed from 
$(\ca L',\D,S)$. 

Following [BLO1] (but with notation which reflects our preference for right-hand composition of morphisms), 
we define the {\it orbit category} $\ca O=\ca O^c(\ca F)$ to be the category whose set of 
objects is $\D$, with 
$$ 
Mor_{\ca O}(P,Q)=Hom_{\ca F}(P,Q)/Inn(Q). 
$$ 
That is, the $\ca O$-morphisms $P\to Q$ are the sets 
$$ 
[\phi]=\{\phi\circ c_x\mid x\in Q\}, 
$$ 
where $\phi$ is an $\ca F$-homomorphism $P\to Q$. If also $\psi:Q\to R$ is an $\ca F$-homomorphism then 
one has the well-defined composition $[\phi]\circ[\psi]=[\phi\circ\psi]$. 

There is a contravariant functor 
$$ 
\ca Z:\ca O^{op}\to Ab
$$ 
(where $Ab$ is the category of abelian groups), given by $\ca Z(P)=Z(P)$ on objects, and defined in the 
following way on $\ca O$-morphisms. If  $\phi:P\to Q$ is an $\ca F$-homomorphism (with $P$ and $Q$ 
centric in $\ca F$), then $\ca Z$ sends the $\ca O$-homomorphism $[\phi]$ to the homomorphism 
$Z(Q)\to Z(P)$ obtained as the composition of the inclusion map $Z(Q)\to Z(P\phi)$ followed by 
the map $\phi_0\i:Z(P\phi)\to Z(P)$, where $\phi_0$ is the $\ca F$-isomorphism $P\to P\phi$ induced 
by $\phi$. The main result of [LL] is: If $\ca F$ is saturated then the higher limit functors 
${\underset\leftarrow\to{lim}}^k(\ca Z)$ are trivial for all $k\geq 2$. 

There is a direct analogy with the theory of group extensions, which estabishes that  
the vanishing of ${\underset\leftarrow\to{lim}}^3(\ca Z)$ implies the existence of a 
$p$-local compact group $(\ca T_{\D}(S)@>\e>> \ca T @>\r>>\ca F)$. Such a $p$-local compact group 
can be viewed as an ``extension" of $\ca O$, in the sense that there is a functor 
$$ 
\ca T@>\s>>\ca O, 
$$ 
such that $\s$ induces the identity map $Ob(\ca T)\to Ob(\ca O)$ (i.e. the identity map on $\D$), and 
such that the image of a $\ca T$-morphism $\psi:P\to Q$ under $\s$ is equal to $[\phi]$, where 
$\phi$ is the $\ca F$-homomorphism $(\psi)\r$. The vanishing of ${\underset\leftarrow\to{lim}}^2(\ca Z)$ 
yields the uniqueness of this extension, up to isomorphism. That is, if 
$(\ca T_{\D}(S)@>\e'>> \ca T' @>\r'>>\ca F)$ is another $p$-local compact group, then there is 
an isomorphism $\b:\ca T\to\ca T'$ of categories, such that the following diagram commutes:  
$$ 
\CD
\ca T    @>\s>>  \ca O \\
@V{\b}VV           @|    \\
\ca T'   @>\s'>>  \ca O
\endCD 
$$ 
(and where $\s'$ is the functor defined by obvious analogy with $\s$). 

It is immediate from the commutativity of the diagram that $\b$ is the identity map on objects. Then 
for each $P\in\D$, $\b$ restricts to a group isomorphism 
$$ 
\b_P:Aut_{\ca T}(P)\to Aut_{\ca T'}(P).  
$$ 
In particular, $\b$ maps the identity element $\iota_P$ to $\iota'_P$, where 
$\iota_P=(\1)\e_P$ and where $\iota'_P=(\1)\e'$. Now let $P\leq Q$ in $\D$. Then 
$\iota_{P,Q}$ is the unique (by A.12(d)) $\ca T$-morphism $\g$ having the property that 
$\iota_P\circ\g=\iota_Q$. Similarly, letting $\iota'_{P,Q}$ denote the image of $\1$ under $\e_{P,Q}$, 
then $\iota'_{P,Q}$ is the unique $\ca T'$-morphism $\g'$ such that $\iota'_P\circ\g'=\iota'_Q$.
Thus $(\iota_{P,Q})\b=\iota'_{P,Q}$. Now let  
$\psi:P\to Q$ and $\bar\psi:\bar P\to\bar Q$ be $\ca T$-isomorphisms such that $\psi\up\bar\psi$. 
That it, assume that $\iota_{P,\bar P}\circ\bar\psi=\psi\circ\iota_{Q,\bar Q}$. Applying $\b$ then 
yields $(\psi)\b\up(\bar\psi)\b$, and thus $\b$ induces a mapping $\a:\ca L\to\ca L'$ (on equivalence 
classes of isomorphisms). The product $\Pi(g_1,\cdots,g_n)$ in $\ca L$ is the equivalence class 
of a composite of a sequence of representatives for $(g_1,\cdots,g_n)$, so $\a$ is a homomorphism 
of partial groups. As $\b$ is invertible, $\a$ is then an isomorphism, as required. 
\qed 
\enddemo

\proclaim {Corollary A.25} Let $\ca F$ be an order-saturated fusion system on a discrete $p$-toral group 
$S$Then $\ca F$. Then $\ca F$ is saturated. 
\endproclaim 

\demo {Proof} By [LL] (and [M]), there exists a $p$-local compact group on $\ca F$. By A.22 there 
then exists a compact locality on $\ca F$, and so $\ca F$ is saturated by 8.3. 
\qed 
\enddemo

\definition {Remark} If $(\ca L,\D,S)$ and $(\ca L',\D,S)$ are proper, {\it finite} localities on $\ca F$, 
with $\D=\ca F^c$, then the uniqueness result in [Ch1] yields an isomorphism $\a:\ca L\to\ca L'$ such 
that $\a$ restricts to the identity map on $S$; and we say that $\a$ is {\it rigid}. 
By Theorem A1 above, one has the same result (i.e. a rigid isomorphism $\ca L\to\ca L'$ of proper, finite  
localities) for any $\ca F$-closed set $\D$ of subgroups of $S$ with $\ca F^{cr}\sub\D\sub\ca F^s$. 

We ask whether the analogous result obtains for compact localities. That is, in Theorem A.25, can $\a$ be 
chosen so that the restriction of $\a$ to $S$ is the identity map on $S$ ? We expect such a result to be 
the case, as we believe that [LL] shows how to extend the existence/uniqueness proof in [Ch1] 
(in which the locality $\ca L$ is constructed from $\ca F$ by an iterative procedure) from the finite 
case to the case of compact localities. Such an exercise would be an enormous undertaking, however, since 
the proof in [Ch1] is far longer than the proof via obstruction theory in [LL]. So, we wonder if there 
is not a better approach to the question. 
\enddefinition

\enddocument 

\ref \key O \paper Existence and uniqueness of linking systems: Chermak's proof via obstruction theory 
\by Bob Oliver  \jour Acta Math. \vol 211 \yr 2013  \pages 141-175 
\endref

\vskip .2in 
\noindent 
{\bf Section 9: The fusion system of a proper locality} 
\vskip .1in 

--(some words. Mention [He].) 

\proclaim {Theorem 9.1} Let $\ca F=\ca F_S(\ca L)$ be the fusion system of a proper locality $(\ca L,\D,S)$. 
Then $\ca F$ is saturated. That is, $\ca F$ is inductive; and for each $V\leq S$ with $V$ fully normalized 
in $\ca F$ the fusion systems $N_{\ca F}(V)$ and $C_{\ca F}(V)$ are $(cr)$-generated. 
\endproclaim 

\demo {Proof} Let $Q\in\D$. Then each $\ca F$-automorphism of $Q$ is given by 
conjugation by an element of $N_{\ca L}(Q)$, and so $Aut_{\ca F}(Q)\cong N_{\ca L}(Q)/C_{\ca L}(Q)$. 
Assume that $Q$ is fully normalized in $\ca F$. Then $N_S(Q)\in Syl_p(N_{\ca L}(Q)$ by 2.1, and hence 
$Aut_S(Q)\in Syl_p(Aut_{\ca F}(Q))$. That is, $Q$ is fully automized in $\ca F$. Now let $P$ be an 
$\ca F$-conjugate of $Q$, let $\phi:P\to Q$ be an $\ca F$-isomorphism, and let $R$ be a subgroup of 
$N_S(P)$ such that 
$$ 
\phi\i Aut_R(P)\phi\leq Aut_S(Q). \tag*
$$ 
Then $\phi$ is conjugation by an 
element $g\in\ca L$, and I.2.3(b) shows that $g$-conjugation is an isomorphism from $N_{\ca L}(P)$ to 
$N_{\ca L}(Q)$. Thus $R^g$ is a $p$-subgroup of $N_{\ca L}(Q)$, and the condition (*) implies that 
$R^g\leq C_{\ca F}(Q)N_S(Q)$. As $N_S(Q)\in Syl_p(C_{\ca F}(Q)N_S(Q))$ there exists $h\in C_{\ca F}(Q)$ 
such that $R^{gh}\leq N_S(Q)$. Then conjugation by $gh$ is an extension of $\phi$ to an 
$\ca F$-homomorphism $\bar\phi:R\to N_S(Q)$; which is to say that $Q$ is receptive in $\ca F$. Thus 
$\ca F$ is $\D$-saturated. By Theorem A1 we may assume $\ca F^c\sub\D$, and then 

[Theorem 2.2 in 5a] 

implies that $\ca F$ is saturated. 

Let $Y\leq S$ be fully normalized in $\ca F$, and let $X$ be an $\ca F$-conjugate of $Y$. Any 
$\ca F$-isomorphism $X\to Y$ conjugates $Aut_{\ca F}(X)$ to $Aut_{\ca F}(Y)$. As $\ca F$ is saturated, 
$Y$ is fully automized, and so there exists an $\ca F$-isomorphism $\psi:X\to Y$ such that $\psi$ 
conjugates $Aut_S(X)$ into $Aut_S(Y)$. As also $Y$ is receptive, $\psi$ extends to an 
$\ca F$-homomorphism $N_S(X)\to N_S(Y)$, and thus $\ca F$ is inductive. 
By [Theorem I.5.5 in AKO], $N_{\ca L}(X)$ and $C_{\ca L}(X)$ are saturated, and these fusion systems are 
then $(cr)$-generated as an immediate consequence of 
Alperin's fusion theorem [Theorem I.3.5 in AKO]. 
\qed 
\enddemo 

We wish to expand on 6.1, by showing that if $V$ is fully normalized in $\ca F$ then $N_{\ca F}(V)$ and 
$C_{\ca F}(V)$ are fusion systems of proper localities. The following lemma will be needed.  

\proclaim {Lemma 6.2} Let $\ca F$ be the fusion system of a proper locality, let $V$ be fully 
normalized in $\ca F$, and let $P$ be a subgroup of $N_S(V)$ containing $V$. Then $P$ is centric in 
$\ca F$ if and only if $P$ is centric in $N_{\ca F}(V)$. 
\endproclaim 

\demo {Proof} Set $\ca F_V=N_{\ca F}(V)$. We are free to replace $P$ by any $\ca F_V$-conjugate of $P$, 
so we may assume that $P$ is fully normalized in $\ca F_V$. As $\ca F$ is inductive by 6.1, $P$ is then 
fully centralized in $\ca F$ by 1.16. Then $P\in\ca F^c$ if and only if $C_S(P)\leq P$ by 1.10. As 
$V\leq P$ we have $C_S(P)=C_{N_S(V)}(P)$, and the lemma follows. 
\qed 
\enddemo

\proclaim {Proposition 6.3} Let $(\ca L,\D,S)$ be a proper locality on $\ca F$, and let $V\leq S$ be a 
subgroup of $S$ such that $V$ is fully normalized in $\ca F$. Then there exists a proper locality 
$(\ca L_V,\D_V,N_S(V))$ on $N_{\ca F}(V)$, and a proper locality $(\ca C_V,\S_V,N_S(V))$ on $C_{\ca F}(V)$. 
\endproclaim 

\demo {Proof} By Theorem A1 we may take $\D=\ca F^c$. Set $\ca F_V=N_{\ca F}(V)$, and set $\D_V=(\ca F_V)^c$. 
Further, set 
$$ 
\ca L_V=\{g\in N_{\ca L}(V)\mid N_{S_g}(V)\in\D_V\}, 
$$ 
and write  
$$ 
\bold D_V=\bold D_{\D_V}=\{w\in\bold W(N_{\ca L}(V))\mid N_{S_w}(V)\in\D_V\}. 
$$ 
Then 6.2 yields $\ca D_V\sub\ca F^c$, and $\D_V=\{P\in\D\mid V\norm P\}$. Notice that $\Pi(w)\in\ca L_V$ 
for any $w\in\bold D_V$, by I.2.5(c). It is then a straightforward exercise with definition I.1.1 
to verify that $\ca L_V$ is a partial group with respect to the restriction $\Pi_V:\bold D_V\to\ca L_V$ of 
$\Pi$, and with respect to the restriction to $\ca L_V$ of the inversion in $\ca L$. Since $(\ca F_V)^c$ is 
$\ca F_V$-closed, it is immediate from the above definition of $\D_V$ and from definition I.2.1 that 
$(\ca L_V,\D_V)$ is objective. As $V$ is fully normalized in $\ca F$, $N_S(V)$ is a maximal $p$-subgroup 
of $\ca L_V$ by I.2.10(b). Thus $(\ca L_V,\D_V,N_S(V)$ is a locality. 
Set $\ca E_V=\ca F_{N_S(V)}(N_{\ca L}(V))$. Then $\ca E_V$ is a fusion subsystem of $\ca F_V$. 
Since $(\ca F_V)^{cr}\sub\D_V$, and since $\ca F_V$ is $(cr)$-generated by 6.1, it follows that 
$\ca E_V=\ca F_V$. Let $P\in\D_V$. Then  
$$ 
N_{\ca L_V}(P)=N_{N_{\ca L}(V)}(P)=N_{N_{\ca L}(P)}(V),  
$$
and hence $N_{\ca L_V}(P)$ is a group of characteristic $p$ by 2.7(b). Thus 
$(\ca L_V,\D_V,N_S(V)$ is a proper locality on $N_{\ca F}(V)$.

We may assume henceforth that $\ca L=\ca L_V$. Set $\S=C_{\ca F}(V)^c$ and define $C_{\ca L}(V)$ to be 
the set of all $g\in\ca L$ such that $g^x=g$ for all $x\in V$. One observes that $C_{\ca L}(V)$ is a 
partial normal subgroup of $\ca L$. Set $H=N_{\ca L}(C_S(V))$. Then 
$\ca L=C_{\ca L}(V)H$ by the Frattini Lemma (I.3.11 -- ?), and 1.5 -- shows that $H$ acts on $C_{\ca F}(V)$. 
Then $\S$ is $\ca F$-invariant. Now let $X\in\S$. Then each $\ca F$-conjugate of $VX$ is of the form 
$VY$ where $Y\in X^{\ca F}$. Then $Y\in\S$, and so 
$$ 
C_S(VY)=C_{C_S(V)}(Y)\leq Y\leq VY.  
$$ 
This shows that $VX\in\ca F^c$ for $X\in\S$. 

Set $\ca C_V=\{g\in C_{\ca L}(V)\mid C_{S_g}(V)\in\S\}$, and write 
$$ 
\bold E=\bold D_{\S}=\{w\in\bold W(\ca C_V)\mid C_{S_w}(V)\in\S\}. 
$$ 
Then $\bold E\sub\bold D$, and $\Pi$ restricts to a mapping $\bold E\to\ca C_V$ which, together with 
the restriction of the inversion in $\ca L$, makes $\ca C_V$ into a partial group, and which makes 
$(\ca C_V,\S)$ into an objective partial group. As $C_S(V)$ is a maximal $p$-subgroup of $C_{\ca L}(V)$  
by I.3.1(c), $C_S(V)$ is also a maximal $p$-subgroup of $C_V$. Thus $(\ca C_V,\S,C_S(V))$ is a locality. 
As $C_{\ca F}(V)^{cr}\sub\S$ 

we find that

$C_{\ca F}(V)=\ca F_{C_S(V)}(\ca C_V)$. 
For $X\in\S$, $N_{\ca C_V}(X)$ is a normal subgroup of the group $N_{\ca L}(VX)$, and then 2.7(a)    
shows that $N_{\ca C_V}(X)$ is of characteristic $p$. Thus $(\ca C_V,\S,C_S(V))$ is a proper locality 
on $C_{\ca F}(V)$. 
\qed 
\enddemo

DOES THIS NEXT GUY REPLACE 8.-- ??

\proclaim {Lemma 6.10} Let $\ca L$ be a proper locality on $\ca F$, let $P\in\ca F^{cr}$, let $T\leq S$ be 
strongly closed in $\ca F$, and let $\ca E$ be an inductive fusion system on $T$ such that $\ca E$ is a 
subsystem of $\ca F$. Then there exists $Q\in P^{\ca F}$ with $Q\cap T\in\ca E^c$. 
\endproclaim 

\demo {Proof} Set $U=P\cap T$ and set $A=N_{C_T(U)}(P)$. Then $[P,A]\leq U$, and thus $A$ centralizes the 
chain $(P\geq U\geq 1)$ of normal subgroups of the group $N_{\ca L}(P)$. Then $A\leq O_p(N_{\ca L}(P))$ 
by 2.7(c), and then $A\leq P$ by an application of 2.3 to the fusion systen $N_{\ca F}(P)$. Thus 
$C_T(U)\leq P$, and so $C_T(U)\leq U$. 

Let $V\in U^{\ca F}$ be fully normalized in $\ca F$. As $\ca F$ is inductive by 6.1, there exists an 
$\ca F$-homomorphism $\phi:N_S(U)\to N_S(V)$ with $U\phi=V$. Set $Q=P\phi$. Then $Q\in\ca F^{cr}$, and 
so $C_T(V)\leq V$ by the result of the preceding paragraph. Let $V'\in V^{\ca E}$. Then $V'\in V^{\ca F}$, 
and so there exists an $\ca F$-homomorphism $\psi:N_S(V')\to N_S(V)$ with $V'=V\psi$. Then 
$N_T(V')\psi\leq N_T(V)$ as $T$ strongly closed in $\ca F$, and thus $|N_T(V')|\leq |N_T(V)|$. This shows 
that $V$ is fully normalized in $\ca E$. As $\ca E$ is inductive, $V$ is fully centralized in $\ca E$ by 
1.13. As $C_T(V)\leq V$, $V$ is then centric in $\ca E$ by 1.10. 
\qed 
\enddemo

\enddocument